\RequirePackage[l2tabu,orthodox]{nag}		
\documentclass[reqno]{amsart}		
\usepackage[margin=1.5in,bottom=1.25in]{geometry}		

\usepackage{amsmath}		
\usepackage{amssymb}		
\usepackage{amsfonts}		
\usepackage{amsthm}		
\usepackage[foot]{amsaddr}		

\usepackage{mathtools}		
\mathtoolsset{%
}

\usepackage[utf8]{inputenc}		
\usepackage[T1]{fontenc}		

\usepackage[proportional,tabular,lining,sf,mono=false]{libertine}

\usepackage{dsfont}		

\usepackage[
cal=cm,
]
{mathalfa}

\usepackage{acronym}		
\renewcommand*{\aclabelfont}[1]{\acsfont{#1}}		
\newcommand{\acli}[1]{\textit{\acl{#1}}}		
\newcommand{\acdef}[1]{\textit{\acl{#1}} \textup{(\acs{#1})}\acused{#1}}		

\usepackage[labelfont={bf,small},labelsep=colon,font=small]{caption}	
\usepackage[dvipsnames,svgnames]{xcolor}		
\colorlet{MyRed}{Crimson!75!Black}
\colorlet{MyBlue}{MediumBlue!90!Black}
\colorlet{MyGreen}{DarkGreen!80!Black}

\newcommand{\afterhead}{.}
\newcommand{\ackperiod}{}		
\newcommand{\para}[1]{\smallskip\paragraph{\textbf{#1\afterhead}}}

\usepackage{quoting}			
\quotingsetup{vskip=\medskipamount}

\usepackage{booktabs}		
\usepackage[inline,shortlabels]{enumitem}		
\setenumerate{itemsep=\smallskipamount,topsep=1ex,left=\parindent}

\usepackage{xspace}		

\usepackage[numbers,sort&compress]{natbib}		

\usepackage{hyperref}
\hypersetup{
colorlinks=true,
linktocpage=true,
pdfstartview=FitH,
breaklinks=true,
pdfpagemode=UseNone,
pageanchor=true,
pdfpagemode=UseOutlines,
plainpages=false,
bookmarksnumbered,
bookmarksopen=false,
bookmarksopenlevel=1,
hypertexnames=true,
pdfhighlight=/O,
urlcolor=MyRed,linkcolor=MyBlue,citecolor=MyBlue,	
pdftitle={},
pdfauthor={},
pdfsubject={},
pdfkeywords={},
pdfcreator={pdfLaTeX},
pdfproducer={LaTeX with hyperref}
}

\usepackage[sort&compress,capitalize,nameinlink]{cleveref}		
\crefname{assumption}{Assumption}{Assumptions}
\crefname{algo}{Algorithm}{Algorithms}
\crefname{example}{Example}{Examples}
\crefname{method}{Method}{Methods}
\creflabelformat{assumption}{\upshape(#2#1#3\upshape)}

\crefname{assumptionenum}{Assumption}{Assumptions}
\creflabelformat{assumptionenum}{#2#1#3}

\crefname{item}{}{}
\creflabelformat{item}{#2#1#3}

\crefname{eq}{}{}
\creflabelformat{eq}{\upshape(#2#1#3\upshape)}

\def\endenv{\hfill{\small$\blacktriangle$}}		

\usepackage{thmtools}		
\usepackage{thm-restate}		

\theoremstyle{plain}
\newtheorem{corollary}{Corollary}		
\newtheorem{lemma}{Lemma}		
\newtheorem{proposition}{Proposition}		

\newtheorem*{corollary*}{Corollary}		

\theoremstyle{definition}
\newtheorem{definition}{Definition}		
\newtheorem{example}{{\small$\blacktriangledown$} Example}		
\newtheorem{algo}{{\small$\blacktriangledown$} Algorithm}		

\newtheorem*{definition*}{Definition}		
\newtheorem*{assumption*}{Assumptions}		
\newtheorem*{example*}{Example}		

\theoremstyle{remark}
\newtheorem{remark}{Remark}		

\newtheorem*{remark*}{Remark}		

\newcounter{proofpart}

\numberwithin{remark}{section}		
\numberwithin{example}{section}		

\usepackage[textwidth=20mm]{todonotes}		

\newcommand{\debug}[1]{#1}		

\newcommand{\commtag}[1]{\tag*{\small\{#1\}}}

\newcommand{\newmacro}[2]{\newcommand{#1}{\debug{#2}}}		
\newcommand{\newop}[2]{\DeclareMathOperator{#1}{\debug{#2}}}		

\DeclarePairedDelimiter{\braces}{\{}{\}}		
\DeclarePairedDelimiter{\bracks}{[}{]}		
\DeclarePairedDelimiter{\parens}{(}{)}		

\DeclarePairedDelimiter{\abs}{\lvert}{\rvert}		

\DeclarePairedDelimiterX{\inner}[2]{\langle}{\rangle}{#1,#2}		
\DeclarePairedDelimiter{\norm}{\lVert}{\rVert}		
\DeclarePairedDelimiterXPP{\twonorm}[1]{}{\lVert}{\rVert}{}{#1}		
\DeclarePairedDelimiterXPP{\dnorm}[1]{}{\lVert}{\rVert}{_{\ast}}{#1}		

\DeclarePairedDelimiterX{\braket}[2]{\langle}{\rangle}{#1,#2}		

\DeclarePairedDelimiterX{\setdef}[2]{\{}{\}}{#1:#2}		
\DeclarePairedDelimiterXPP{\exclude}[1]{\mathopen{}\setminus}{\{}{\}}{}{#1}

\newcommand{\alt}[1]{#1'}		

\newcommand{\N}{\mathbb{N}}		
\newcommand{\R}{\mathbb{R}}		

\DeclareMathOperator*{\intersect}{\bigcap}		
\DeclareMathOperator*{\union}{\bigcup}		

\DeclareMathOperator{\bigoh}{\mathcal{O}}		
\DeclareMathOperator{\cl}{cl}		
\DeclareMathOperator{\crit}{crit}		
\DeclareMathOperator{\dist}{dist}		
\DeclareMathOperator{\grad}{\nabla}		
\DeclareMathOperator{\Jac}{D}		
\DeclareMathOperator{\one}{\mathds{1}}		

\newmacro{\coef}{\lambda}		
\newmacro{\dd}{\:d}		

\newcommand{\eps}{\varepsilon}		
\newcommand{\pd}{\partial}		

\newcommand{\insum}{\sum\nolimits}		

\newmacro{\pexp}{p}		
\newmacro{\qexp}{q}		
\newmacro{\rexp}{r}		

\newcommand{\cf}{cf.\xspace}		
\newcommand{\eg}{e.g.,\xspace}		
\newcommand{\ie}{i.e.,\xspace}		
\newcommand{\vs}{vs.\xspace}		

\newcommand{\textpar}[1]{\textup(#1\textup)}		

\newcommand{\txs}{\textstyle}		

\newcommand{\from}{\colon}		

\newcommand{\defeq}{\coloneqq}		

\newmacro{\set}{\mathcal{S}}		

\newmacro{\points}{\mathcal{Z}}		
\newmacro{\intpoints}{\points^{\circ}}		
\newmacro{\point}{z}		
\newmacro{\pointalt}{\alt\point}		

\newmacro{\dpoints}{\mathcal{W}}		
\newmacro{\dpoint}{w}		
\newmacro{\dpointalt}{\alt\dpoint}		

\newmacro{\base}{p}		
\newmacro{\basealt}{q}		

\newmacro{\open}{\mathcal{U}}		
\newmacro{\closed}{\mathcal{C}}		
\newmacro{\cpt}{\mathcal{K}}		
\newmacro{\nhd}{\mathcal{U}}		

\newmacro{\start}{1}		
\newmacro{\halfafterstart}{3/2}		
\newmacro{\afterstart}{2}		
\newmacro{\running}{\start,\afterstart,\dotsc}		
\newmacro{\halfrunning}{\start,\halfafterstart,\dotsc}

\newmacro{\runalt}{k}		
\newmacro{\run}{n}		
\newmacro{\nRuns}{T}		
\newmacro{\runs}{\mathcal{\nRuns}}		

\newmacro{\state}{Z}		
\newmacro{\dstate}{Y}		

\newcommand{\init}[1][\state]{\debug{#1}_{\start}}		
\newcommand{\preiter}[1][\state]{\debug{#1}_{\runalt-1}}		
\newcommand{\iter}[1][\state]{\debug{#1}_{\runalt}}		
\newcommand{\afteriter}[1][\state]{\debug{#1}_{\runalt+1}}		
\newcommand{\preprev}[1][\state]{\debug{#1}_{\run-2}}		
\newcommand{\prev}[1][\state]{\debug{#1}_{\run-1}}		
\newcommand{\curr}[1][\state]{\debug{#1}_{\run}}		
\newcommand{\prelead}[1][\state]{\debug{#1}_{\run-1}^{+}}		
\newcommand{\lead}[1][\state]{\debug{#1}_{\run}^{+}}		
\renewcommand{\next}[1][\state]{\debug{#1}_{\run+1}}		

\newmacro{\ctime}{t}		
\newmacro{\ctimealt}{s}		
\newmacro{\cstart}{0}		

\newmacro{\horizon}{T}		

\newmacro{\vecspace}{\R^{\vdim}}		

\newmacro{\coord}{i}		
\newmacro{\vdim}{d}		
\newmacro{\vvec}{x}		
\newmacro{\bvec}{e}		
\newmacro{\bvecs}{\mathcal{E}}		

\newmacro{\subspace}{\mathcal{W}}		
\newmacro{\wvec}{w}		
\newmacro{\subdim}{m}		

\newmacro{\tanhull}{\mathcal{Z}}		
\newmacro{\tanvec}{z}		

\newcommand{\dual}[1]{#1^{\ast}}		
\newmacro{\dspace}{\dual\vecspace}		
\newmacro{\dvec}{v}		
\newmacro{\dbvec}{\eps}		

\newmacro{\ones}{\mathbf{1}}		
\newmacro{\mat}{M}		
\newmacro{\eye}{I}		

\newcommand{\tspace}[1][\point]{ T_{#1}  }		

\newop{\tcone}{TC}		
\newop{\dcone}{\dual\tcone}		
\newop{\ncone}{NC}		
\newop{\pcone}{PC}		

\newmacro{\cvx}{\mathcal{C}}		
\newmacro{\subd}{\partial}		
\newmacro{\jmat}{J\vecfield}		
\newmacro{\hmat}{H}		

\newop{\Opt}{Opt}		
\newop{\Sol}{Sol}		

\newmacro{\obj}{f}		
\newmacro{\objalt}{g}		
\newmacro{\sobj}{F}		

\newmacro{\param}{\theta}		
\newmacro{\params}{\Theta}		

\newmacro{\gvec}{g}		
\newmacro{\vecfield}{V}		

\newmacro{\gbound}{G}		
\newmacro{\vbound}{M}		

\newcommand{\sol}[1][\point]{#1^{\ast}}		
\newcommand{\sols}{\sol[\points]}		

\newmacro{\strong}{\ell}		
\newmacro{\smooth}{\beta}		
\newmacro{\lips}{L}		

\newmacro{\minmax}{\Phi}		

\newmacro{\minvar}{x}		
\newmacro{\minvaralt}{\alt x}		
\newmacro{\minvars}{\mathcal{X}}		

\newmacro{\maxvar}{y}		
\newmacro{\maxvaralt}{\alt y}		
\newmacro{\maxvars}{\mathcal{Y}}		

\newmacro{\minsol}{\sol[\minvar]}		
\newmacro{\maxsol}{\sol[\maxvar]}		

\newmacro{\mindim}{\vdim_{\minvars}}		
\newmacro{\maxdim}{\vdim_{\maxvars}}		

\newmacro{\minstate}{X}		
\newmacro{\maxstate}{Y}		

\newmacro{\payoffmat}{A}		
\newmacro{\perturb}{\phi}	   

\newop{\NE}{NE}		
\newop{\CE}{CE}		
\newop{\CCE}{CCE}		

\newop{\brep}{br}		
\newop{\reg}{Reg}		
\newop{\preg}{\overline{Reg}}		
\newop{\val}{val}		

\newmacro{\play}{i}		
\newmacro{\playalt}{j}		
\newmacro{\nPlayers}{N}		
\newmacro{\players}{\mathcal{\nPlayers}}		

\newmacro{\pure}{a}		
\newmacro{\purealt}{a'}		
\newmacro{\nPures}{A}		
\newmacro{\pures}{\mathcal{\nPures}}		

\newmacro{\cost}{c}		
\newmacro{\loss}{\ell}		
\newmacro{\pay}{u}		
\newmacro{\payv}{v}		
\newmacro{\pot}{\obj}		

\newmacro{\game}{\mathcal{G}}		
\newmacro{\gamefull}{\game(\players,\points,\pay)}		

\newmacro{\fingame}{\Gamma}		
\newmacro{\fingamefull}{\Gamma(\players,\pures,\pay)}		

\newop{\Eucl}{\Pi}		
\newop{\logit}{\Lambda}		

\newmacro{\hreg}{h}		
\newmacro{\breg}{D}		
\newmacro{\pmap}{P}		
\newmacro{\mirror}{Q}		
\newmacro{\fench}{F}		
\newmacro{\hstr}{K}		
\newmacro{\depth}{H}		
\newmacro{\zone}{\mathbb{D}}		
\newmacro{\subpoints}{\points^{\circ}}		

\DeclareMathOperator{\ex}{\mathbb{E}}		
\DeclareMathOperator{\prob}{\mathbb{P}}		

\newmacro{\seed}{\omega}		
\newmacro{\seeds}{\Omega}		
\newmacro{\history}{\mathcal{H}}		

\newmacro{\sample}{\omega}		
\newmacro{\samples}{\Omega}		
\newmacro{\filter}{\mathcal{F}}		
\newmacro{\probspace}{(\samples,\filter,\prob)}		

\newmacro{\event}{\mathcal{E}}       
\newmacro{\eventalt}{\mathcal{H}}       
\newcommand{\comp}[1]{#1^{\mathtt{c}}}		

\newmacro{\mean}{\mu}		
\newmacro{\sdev}{\sigma}		
\newmacro{\variance}{\sdev^{2}}		

\newmacro{\dkl}{D_{\mathrm{KL}}}		
\newcommand{\as}{\debug{\textpar{a.s.}}\xspace}		

\providecommand\given{}		

\DeclarePairedDelimiterXPP{\exof}[1]{\ex}{[}{]}{}{
\renewcommand\given{\nonscript\,\delimsize\vert\nonscript\,\mathopen{}} #1}

\DeclarePairedDelimiterXPP{\probof}[1]{\prob}{(}{)}{}{
\renewcommand\given{\nonscript\:\delimsize\vert\nonscript\:\mathopen{}} #1}

\newcommand{\oneof}[1]{\one_{\{#1\}}}

\newmacro{\step}{\gamma}		
\newmacro{\temp}{\eta}		

\newmacro{\efftime}{\tau}		
\newmacro{\tinv}{M}		
\newcommand{\apt}[2][]{\state_{#1}(#2)}		

\newop{\orcl}{\mathsf{V}}		
\newop{\err}{\mathsf{U}}		

\newmacro{\signal}{V}		
\newmacro{\error}{W}		
\newmacro{\noise}{U}		
\newmacro{\bias}{b}		
\newmacro{\brown}{W}		

\newmacro{\bbound}{B}		

\newmacro{\totbound}{S}		
\newmacro{\noisepar}{\sdev}		
\newmacro{\noisevar}{\variance}		

\newmacro{\snoise}{\xi}		
\newmacro{\sbias}{\psi}		
\newmacro{\scorr}{\theta}		

\newmacro{\mix}{\delta}		
\newmacro{\unitvar}{E}		
\newmacro{\pertvar}{W}		

\newmacro{\radius}{r}		

\usepackage{twoopt}

\newmacro{\limset}{\mathcal{L}}		

\newcommand{\orbit}[2][]{\point_{#1}(#2)}		
\newcommand{\dotorbit}[2][]{\dot\point_{#1}(#2)}		

\newmacro{\flowmap}{\Theta}		
\newcommandtwoopt{\flow}[2][\ctime][\point]{\flowmap_{#1}(#2)}

\newcommand{\umfd}[1][\point]{\mathcal{E}^u_{{#1}}} 	

\newmacro{\graph}{\mathcal{G}}
\newmacro{\vertices}{\mathcal{V}}
\newmacro{\edges}{\mathcal{E}}

\newmacro{\gmat}{g}		
\newmacro{\gdist}{\dist_{\gmat}}
\newmacro{\ball}{\mathbb{B}}		
\newmacro{\sphere}{\mathbb{S}}		

\newcommand{\Grass}[1][\subdim]{
G({#1},\GrassRelay)
}      

\newcommand{\GrassRelay}[1][\vdim]{#1}		

\newmacro{\const}{C}
\newmacro{\conf}{\alpha}
\newmacro{\lyap}{E}

\newmacro{\aux}{\tilde\lyap}		

\begin{document}


\newcommand{\longtitle}{\uppercase{The Limits of Min-Max Optimization Algorithms:\\
Convergence to Spurious Non-Critical Sets}}
\newcommand{\runtitle}{\uppercase{The Spurious Limits of Min-Max Optimization Algorithms}}		

\title
[\runtitle]
{\longtitle}		

\author
[Y.~P.~Hsieh]
{Ya-Ping Hsieh$^{\ast}$}
\address{$^{\ast}$\,%
LIONS, École Polytechnique Fédérale de Lausanne (EPFL).}
\email{ya-ping.hsieh@epfl.ch}

\author
[P.~Mertikopoulos]
{Panayotis Mertikopoulos$^{\diamond,\sharp}$}
\address{$^{\diamond}$\,%
Univ. Grenoble Alpes, CNRS, Inria, LIG, 38000, Grenoble, France.}
\address{$^{\sharp}$\,%
Criteo AI Lab.}
\email{panayotis.mertikopoulos@imag.fr}

\author
[V.~Cevher]
{Volkan Cevher$^{\ast}$}
\email{volkan.cevher@epfl.ch}



\thanks{
%
%
The authors are grateful to Thomas Pethick for his help in the numerical simulation of adaptive methods.
This research was partially supported by the COST Action CA16228 ``European Network for Game Theory'' (GAMENET),
the Army Research Office under grant number W911NF-19-1-0404, the Swiss National Science Foundation (SNSF) under  grant number 200021\_178865 / 1, the European Research Council (ERC) under the European Union's Horizon 2020 research and innovation programme (grant agreement $n^\circ$ 725594 - time-data), and 2019 Google Faculty Research Award.
P.~Mertikopoulos is also grateful for financial support by
the French National Research Agency (ANR) under grant no.~ANR\textendash 16\textendash CE33\textendash 0004\textendash 01 (ORACLESS)\ackperiod}

\newacro{LHS}{left-hand side}
\newacro{RHS}{right-hand side}
\newacro{iid}[i.i.d.]{independent and identically distributed}
\newacro{lsc}[l.s.c.]{lower semi-continuous}

\newacro{APT}{asymptotic pseudotrajectory}
\newacroplural{APT}[APTs]{asymptotic pseudotrajectories}
\newacro{GD}{gradient dynamics}
\newacro{GF}{gradient flow}
\newacro{ICT}{internally chain-transitive}
\newacro{MDS}{martingale difference sequence}
\newacro{NE}{Nash equilibrium}
\newacroplural{NE}[NE]{Nash equilibria}
\newacro{ODE}{ordinary differential equation}
\newacro{SA}{stochastic approximation}
\newacro{SFO}{stochastic first-order oracle}
\newacro{SG}{stochastic gradient}
\newacro{SP}{saddle-point}
\newacro{WAC}{weak asymptotic coercivity}

\newacro{AH}{Arrow\textendash Hurwicz}
\newacro{BDG}{Burkholder\textendash Davis\textendash Gundy}
\newacro{ConO}{consensus optimization}
\newacro{RM}{Robbins\textendash Monro}
\newacro{KW}{Kiefer\textendash Wolfowitz}
\newacro{GDA}{gradient descent/ascent}
\newacro{SGA}{symplectic gradient adjustment}
\newacro{SGD}{stochastic gradient descent}
\newacro{SGDA}{stochastic gradient descent/ascent}
\newacro{SPSA}{simultaneous perturbation stochastic approximation}
\newacro{ASGDA}[alt-SGDA]{alternating stochastic gradient descent/ascent}
\newacro{SEG}{stochastic extra-gradient}
\newacro{EG}{extra-gradient}
\newacro{PEG}{Popov's extra-gradient}
\newacro{RG}{reflected gradient}
\newacro{OG}{optimistic gradient}
\newacro{PPM}{proximal point method}

\newacro{GAN}{generative adversarial network}
\newacro{NN}{neural network}
\newacro{FTRL}{``follow the regularized leader''}
\newacro{CGD}{Competitive Gradient Descent}
\newacro{wp1}[w.p.$1$]{with probability $1$}

\begin{abstract}
%
%
Compared to ordinary function minimization problems, min-max optimization algorithms encounter far greater challenges because of the existence of periodic cycles and similar phenomena.
Even though some of these behaviors can be overcome in the convex-concave regime, the general case is considerably more difficult.
On that account, we take an in-depth look at a comprehensive class of state-of-the art algorithms and prevalent heuristics in \emph{non-convex / non-concave} problems, and we establish the following general results:
\begin{enumerate*}
[\itshape a\upshape)]
\item
generically, the algorithms' limit points are contained in the \acdef{ICT} sets of a common, mean-field system;
\item
the attractors of this system also attract the algorithms in question with arbitrarily high probability;
and
\item
all algorithms avoid the system's unstable sets with probability $1$.
\end{enumerate*}
On the surface, this provides a highly optimistic outlook for min-max algorithms;
however, we show that there exist \emph{spurious attractors} that do not contain \emph{any} stationary points of the problem under study.
In this regard, our work suggests that existing min-max algorithms may be subject to inescapable convergence failures.
We complement our theoretical analysis by illustrating such attractors in simple, two-dimensional, almost bilinear problems.
\end{abstract}
\acresetall

\maketitle
\acresetall

\section{Introduction}
\label{sec:introduction}

Consider a min-max optimization \textendash\ or \acli{SP} \textendash\ problem of the form
\begin{equation}
\tag{SP}
\label{eq:minmax}
\min_{\minvar\in\minvars} \max_{\maxvar\in\maxvars} \minmax(\minvar,\maxvar).
\end{equation}
Given an algorithm for solving \eqref{eq:minmax}, it is then natural to ask:
\begin{equation}
\label{eq:question}
\tag{$\star$}
\textit{Where does the algorithm converge to?}
\end{equation}
The goal of our paper is
to treat \eqref{eq:question} in a general non-convex\,/\,non-concave setting
and
to provide answers for a comprehensive array of state-of-the-art algorithms.

\para{Related work}

This question has attracted significant interest in the machine learning literature because of its potential implications to
\aclp{GAN} \citep{GPAM+14},
robust reinforcement learning \citep{PDSG17},
and
other models of adversarial training \citep{MMST+18}.
In this broad setting, it has become empirically clear that the joint training of two \acp{NN} is fundamentally more difficult than that of a \emph{single} \ac{NN} of similar size and architecture.
The latter task boils down to successfully finding a (good) local minimum of a non-convex function, so it is instructive to revisit \eqref{eq:question} in the context of \emph{non-convex minimization}.

In this case, the existing convergence theory for \acdef{SGD} \textendash\ the ``gold standard'' for deep \ac{NN} training \textendash\ can be informally summed up as follows:
\begin{enumerate}
[left=\parindent]
\item
\ac{SGD} always converges to critical points.
\item
\ac{SGD} does not converge to strict saddle points or other spurious solutions.
\end{enumerate}
These results could be seen as plausible expectations for algorithmic proposals to solve \eqref{eq:minmax}.
Unfortunately however, there are well-known examples of simple \emph{bilinear} min-max games where \ac{SGDA}, the min-max analogue of \ac{SGD}, leads to recurrent orbits that do not contain \emph{any} critical point of $\minmax$.
Such \emph{spurious convergence} phenomena arise from the min-max structure of \eqref{eq:minmax} and have no counterpart in minimization problems.

This well-documented failure of \ac{SGDA} has led to an extensive literature that is impossible to survey here.
As a purely indicative \textendash\ and highly incomplete \textendash\ list, we mention the works of \citet{DISZ18}, \citet{GBVV+19}, \citet{MLZF+19} and \citet{MOP19a}, who studied how these failures can be overcome in \emph{deterministic} bilinear problems by means of an \acli{EG} step (or an optimistic proxy thereof).
By contrast, in \emph{stochastic} problems, the convergence of optimistic\,/\,\acl{EG} methods is compromised unless additional, tailor-made mitigation mechanisms are put in place \textendash\ such as variance reduction \cite{IJOT17,CGFLJ19} or variable step-size schedules \cite{HIMM20}.
This shows that the convergence of min-max training methods can be particularly fragile, even in simple, bilinear problems.

Beyond the class of convex-concave problems analyzed above, another vigorous thread of research has focused on the \emph{local analysis} of a min-max optimization algorithm close to the game's critical points \textendash\ typically subject to a second-order sufficient condition; \cf \citet{heusel2017gans, nagarajan2017gradient, daskalakis2018limit, adolphs2019local, mazumdar2020gradient, fiez2020gradient, grimmer2020landscape, grimmer2020limiting}.
The global analysis is much more challenging and requires strong structural assumptions such as variational coherence \citep{MLZF+19} and/or the existence of a Minty-type solution \citep{liu2019towards}.
In the absence of such conditions, \citet{FVGP19a,flokas2021solving} showed that periodic and/or Poincaré recurrent behavior may persist in deterministic, continuous-time min-max dynamics.

From a practical viewpoint, these studies have led to a broad array of sophisticated algorithmic proposals for solving min-max games;
we review many of these algorithms in \cref{sec:algorithms}.
However, a central question that remains unanswered is whether it is theoretically plausible to expect a qualitatively different behavior relative to \ac{SGDA} in the full spectrum of non-convex\,/\,non-concave games.
Our work aims to provide concrete answers to this question.

\para{Our contributions}

Our first contribution is to provide a unified framework for a comprehensive selection of first- and zeroth-order min-max optimization methods (including \ac{SGDA}, \aclp{PPM}, optimistic / extra-gradient schemes, their alternating variants, etc.).
The principal ingredients of our approach are twofold:
\begin{enumerate*}
[(\itshape i\hspace*{1pt}\upshape)]
\item
a generalized \ac{RM} template that is wide enough to include all the above algorithms;
and
\item
an analytic framework leveraging the \ac{ODE} method of stochastic approximation \citep{Ben99,KY97}.
\end{enumerate*}
Based on these two elements, we prove a precise version of the following general principle:
\emph{the long-run behavior of all generalized \ac{RM} methods can be mapped to the study of\:\:\textbf{the same}, mean-field dynamical system.}

In more detail, we show that the limit points of all generalized \ac{RM} schemes belong to an \acdef{ICT} set of these mean dynamics.
The notion of an \ac{ICT} set is central in the study of dynamical systems \citep{Bow75,Con78,BH96} and, in some cases, they are easy to characterize:
in minimization problems (and possibly up to a ``hidden'' transformation in the spirit of \citealp{FVGP19a}), the dynamics' \ac{ICT} sets are the function's critical points.
As such, in this case, we recover \emph{exactly} the min-min landscape of \ac{SGD} \textendash\ but for an \emph{entire family} of algorithms, not just \ac{SGD}.

Moving on to \emph{general} min-max problems, the structure of the dynamics' \ac{ICT} sets could be considerably more complicated, so we provide two further, complementing results:
\begin{enumerate}
[left=\parindent]
\item
\emph{With high probability, all generalized \ac{RM} methods converge locally to attractors of the mean dynamics}.
\item
\emph{With probability $1$, all generalized \ac{RM} methods avoid the mean dynamics' unstable invariant sets}.
\end{enumerate}

As far as we are aware, there are no results of comparable generality in the min-max optimization literature.
From a high level, these theoretical contributions would seem to be analogous to existing results for \ac{SGD} in minimization problems (\ie that \ac{SGD} converges to critical points while avoiding strict saddles).
However, this similarity is only skin-deep:
as we show by a range of concrete, \emph{almost bilinear} examples, min-max optimization algorithms may encounter a series of immovable roadblocks.
Specifically:

\begin{itemize}
[left=\parindent,itemsep=\smallskipamount]
\item
An \ac{ICT} set may contain a \emph{globally attracting limit cycle}, and the range of algorithms under consideration cannot escape it \textendash\ even though \acl{EG} methods escape recurrent orbits in exact bilinear problems.
This suggests that bilinear games may not be representative as a testbed for \acs{GAN} training algorithms and heuristics.
\item
There exist \emph{unstable} critical points whose neighborhood contains an  (almost) \emph{globally stable} \ac{ICT} set.
Therefore, in sharp contrast to minimization, ``avoiding unstable critical points'' \emph{does not imply} ``escaping unstable critical points'' in min-max problems. 
\item
There exist \emph{stable} min-max points whose basin of attraction is ``shielded'' by an \emph{unstable} \ac{ICT} set.
As a result, if run with non-negligible noise in the gradients, then, with high probability, existing algorithms are repelled away from the desirable solutions.
\end{itemize}

Our results indicate a steep, qualitative increase in difficulty when passing from min-min to min-max problems, in line with concurrent works by \citet{daskalakis2020complexity} and \citet{letcher2020impossibility}.
In plain terms,
\citet{daskalakis2020complexity} proved the impossibility of attaining a critical point in polynomial time in deterministic, constrained min-max games.
In a similar spirit, the concurrent work of \citet{letcher2020impossibility} showed that there are min-max games where all ``reasonable'' deterministic algorithms may fail to converge.
By contrast, our paper focuses on the occurrence of \emph{spurious convergence phenomena} with probability $1$ in \emph{stochastic} algorithms.
In addition, our avoidance result (\cref{thm:repeller}) can be seen as a stochastic counterpart of the ``reasonableness'' requirement of \citet{letcher2020impossibility}, thereby enriching the applicability of the results therein.
Taken together, these works and our own provide a complementing look into the fundamental limits of min-max optimization algorithms.
\acresetall

\section{Setup and preliminaries}
\label{sec:setup}

Throughout our paper, we focus on general unconstrained problems with $\minvars = \R^{\mindim}$, $\maxvars=\R^{\maxdim}$, and $\minmax$ assumed $C^{1}$ and Lipschitz.
To avoid unnecessary notation, we will let $\point = (\minvar,\maxvar)$, $\points = \minvars \times \maxvars$ and $\vdim = \mindim+\maxdim$.
In addition, we will write
\begin{equation}
\vecfield(\point)
	\equiv (\vecfield_{\minvar}(\minvar,\maxvar),\vecfield_{\maxvar}(\minvar,\maxvar))
	\defeq (-\nabla_{\minvar}\minmax(\minvar,\maxvar),\nabla_{\maxvar}\minmax(\minvar,\maxvar))
\end{equation}
for the (min-max) gradient field of $\minmax$,
assumed here to be Lipschitz;
in some cases we may also require $\vecfield$ to be $C^{1}$ and write $\jmat(\point)$ for its Jacobian.
Finally, we will assume that $\vecfield$ satisfies the \acl{WAC} condition
\begin{equation}
\label[eq]{asm:coercive}
\braket{\vecfield(\point)}{\point}
	\leq 0
	\quad
	\text{for all sufficiently large $\point$}.
\end{equation}
This condition is a weaker version of standard coercivity conditions in the literature \citep{BC17}, it is satisfied by all convex-concave problems (including bilinear ones) and, importantly, it does not impose any growth requirements on the elements of $\vecfield$ (as standard coercivity conditions do).
We discuss it further in \cref{app:coercive}.

A \emph{solution} of \eqref{eq:minmax} is a tuple $\sol = (\minsol,\maxsol)$ with $\minmax(\minsol,\maxvar) \leq \minmax(\minsol,\maxsol) \leq \minmax(\minvar,\maxsol)$ for all $\minvar\in\minvars$, $\maxvar\in\maxvars$;
likewise, a \emph{local solution} of \eqref{eq:minmax} is a tuple $(\minsol,\maxsol)$ that satisfies this inequality locally.
Finally, a state $\sol$ with $\vecfield(\sol) = 0$ is said to be a \emph{critical} (or \emph{stationary}) \emph{point} of $\minmax$.

From an algorithmic standpoint, we will focus exclusively on the black-box optimization paradigm \citep{Nes04} with \acdef{SFO} feedback.
Algorithms with a more complicated feedback structure, such as a best-response oracle \cite{jin2019local, naveiro2019gradient, fiez2019convergence} or based on mixed-strategy sampling \cite{hsieh2019finding,domingo2020mean, kamalaruban2020robust}, are not considered in this work.

Specifically, when called at $\point = (\minvar,\maxvar)$ with random seed $\sample\in\samples$, an \ac{SFO} returns a random vector $\orcl(\point;\seed) \equiv (\orcl_{\minvar}(\point;\seed),\orcl_{\maxvar}(\point;\seed))$ of the form
\begin{equation}
\label{eq:SFO}
\tag{SFO}
\orcl(\point;\seed)
	= \vecfield(\point)
		+ \err(\point;\seed)
\end{equation}
where the error term $\err(\point;\seed)$ captures all sources of uncertainty in the model (\eg the selection of a minibatch in \acs{GAN} training, system state observations in reinforcement learning, etc.).
As is standard in the literature, we require $\err(\point;\seed)$ to be zero-mean and finite-variance:
\begin{equation}
\label{eq:err-base}
\forall \point \in \points, \quad
\exof{\err(\point;\seed)}
	= 0 
	\text{ and }
\exof{\norm{\err(\point;\seed)}^{2}}
	\leq \sdev^{2}.
\end{equation}
These will be our blanket assumptions throughout.

\section{Core algorithmic framework}
\label{sec:algorithms}

\subsection{The \acl{RM} template}

Much of our analysis will revolve around iterative algorithms that can be cast as generalized \acl{RM} algorithms \citep{RM51} of the general form
\begin{equation}
\label{eq:RM}
\tag{RM}
\next
	= \curr
		+ \curr[\step] \bracks{ \vecfield(\curr) + \curr[\error]}
\end{equation}
where:
\begin{enumerate}
[]
\item
$\curr = (\curr[X],\curr[Y]) \in \points$ denotes the state of the algorithm at each stage $\run = \running$
\item
$\curr[\error]$ is an abstract error term described in detail below.
\item
$\curr[\step]$ is the method's step-size hyperparameter, and is typically of the form $\curr[\step] \propto 1/\run^{\pexp}$ for some $\pexp\geq0$. Throughout the paper, we will always assume $\sum_{\run} \curr[\step] = \infty$ and $\lim_{\run} \curr[\step] = 0$.
\end{enumerate}

In the above, the error term $\curr[\error]$ is generated \emph{after} $\curr$;
thus, by default, $\curr[\error]$ is not adapted to the history $\curr[\filter] \defeq \history(\init,\dotsc,\curr)$ of $\curr$.
For concision,
we will also write
\begin{equation}
\label{eq:signal}
\curr[\signal]
	= \vecfield(\curr)
		+ \curr[\error]
\end{equation}
so $\curr[\signal]$ can be seen as a noisy estimator of $\vecfield(\curr)$.
In more detail, to differentiate between ``random'' (zero-mean) and ``systematic'' (non-zero-mean) errors in $\curr[\signal]$ it will be convenient to further decompose the error process $\curr[\error]$ as
\begin{equation}
\label{eq:error}
\curr[\error]
	= \curr[\noise]
		+ \curr[\bias]
\end{equation}
where
$\curr[\bias] = \exof{\curr[\error] \given \curr[\filter]}$ represents the systematic component
and
$\curr[\noise] = \curr[\error] - \curr[\bias]$ captures the random, zero-mean part.
In view of all this, we will consider the following descriptors for $\curr[\error]$:
\begin{subequations}
\label{eq:signal-stats}
\begin{alignat}{3}
\label{eq:bias}
a)\quad
	&\textit{Bias:}
	&\quad
	\curr[\bbound]
		&= \exof{\norm{\curr[\bias]} \given \curr[\filter]}
	\\
b)\quad
\label{eq:variance}
	&\textit{Variance:}
	&\quad\curr[\sdev]^{2}
		&= \exof{\norm{\curr[\noise]}^{2} \vert \curr[\filter]}
	\hspace{20em}
\end{alignat}
\end{subequations}
Note that both $\curr[\bbound]$ and $\curr[\sdev]$ are random (conditioned on $\curr[\filter]$);
this will play an important part in the sequel.

\subsection{Specific algorithms} \label{subsec:specific}

In the rest of this section, we discuss how a wide range of algorithms used in the literature can be seen as special instances of our general \ac{RM} template.
\smallskip

\begin{algo}[\Acl{SGDA}]
\label{alg:SGDA}
The basic \ac{SGDA} algorithm \textendash\ also known as the \acli{AH} method \citep{AHU58} \textendash\ queries an \ac{SFO} and proceeds as:
\begin{equation}
\label{eq:SGDA}
\tag{SGDA}
\next
	= \curr
		+ \curr[\step] \orcl(\curr;\curr[\seed]),
\end{equation}
where $\curr[\seed] \in \seeds$ ($\run=\running$) is an \ac{iid} sequence of oracle seeds.
As such, \eqref{eq:SGDA} admits a straightforward \ac{RM} representation by taking $\curr[\error] = \curr[\noise] = \err(\curr;\curr[\seed])$ and $\curr[\bias] = 0$.
\endenv
\end{algo}
\smallskip

\begin{algo}[\Acl{PPM}]
\label{alg:PPM}
The (deterministic) \acdef{PPM} \cite{rockafellar1976monotone} is an implicit update rule of the form:
\begin{equation}
\label{eq:PPM}
\tag{PPM}
\next
	= \curr
		+ \curr[\step] \vecfield(\next).
\end{equation}
The \ac{RM} representation of \eqref{eq:PPM} is obtained by taking $\curr[\error] = \curr[\bias] =  \vecfield(\next) -  \vecfield(\curr)$ and $\curr[\noise] = 0$.
\endenv
\end{algo}
\smallskip

\begin{algo}[\Acl{SEG}]
\label{alg:SEG}
Since \eqref{eq:PPM} is only implicitly defined, one can rarely run it in practice. Nonetheless, it is possible to approximate \eqref{eq:PPM} by locally querying two (stochastic) gradients at each iteration \cite{Nem04}.
This can be achieved by the \acdef{SEG}:
\begin{equation}
\label{eq:SEG}
\tag{SEG}
\begin{alignedat}{2}
\lead
	&= \curr
		+ \curr[\step] \orcl(\curr;\curr[\seed]),
	\\
\next
	&= \curr
		+ \curr[\step] \orcl(\lead;\lead[\seed]).
\end{alignedat}
\end{equation}
To recast \eqref{eq:SEG} in the \acl{RM} framework, simply take
$\curr[\error] = \orcl(\lead;\lead[\seed]) - \vecfield(\curr)$,
\ie
$\curr[\noise] = \err(\lead;\lead[\seed])$
and
$\curr[\bias] = \vecfield(\lead) - \vecfield(\curr)$.
\endenv
\end{algo}
\smallskip

\begin{algo}[\Acl{OG} / \Acl{PEG}]
\label{alg:PEG}
Compared to \eqref{eq:SGDA}, the scheme \eqref{eq:SEG} involves two oracle queries per iteration, which is considerably more costly.
An alternative iterative method with a single oracle query per iteration was proposed by \citet{Pop80}:
\begin{equation}
\label{eq:PEG}
\tag{OG/PEG}
\begin{alignedat}{2}
\lead
	&= \curr
		+ \curr[\step] \orcl(\prelead;\prev[\seed]),
	\\
\next
	&= \curr
		+ \curr[\step] \orcl(\lead;\curr[\seed]).
\end{alignedat}
\end{equation}
\acl{PEG}\acused{PEG} has been rediscovered several times and is more widely known as the \acf{OG} method in the machine learning literature \cite{RS13-COLT,CYLM+12,DISZ18,HIMM19}.
In unconstrained problems, \eqref{eq:PEG} turns out to be equivalent to a number of other existing methods, including
``extrapolation from the past'' \cite{GBVV+19} and
reflected gradient \cite{malitsky2020forward}.
Its \acl{RM} representation is obtained by setting
$\curr[\error] = \orcl(\lead;\curr[\seed]) - \vecfield(\curr)$,
\ie
$\curr[\noise] = \err(\lead;\curr[\seed])$
and
$\curr[\bias] = \vecfield(\lead) - \vecfield(\curr)$.
\endenv
\end{algo}

\begin{algo}[\acl{KW}]
\label{alg:SPSA}
When first-order feedback is unavailable, a popular alternative is to obtain gradient information of $\minmax$ via zeroth-order observations \cite{liu2019min}.
This idea can be traced back to the seminal work of \citet{KW52} and the subsequent development of the \ac{SPSA} method by \citet{Spa92}.
In our setting, this leads to the recursion:
\begin{equation}
\label{eq:SPSA}
\tag{SPSA}
\begin{aligned}
\curr[\signal]
	&= \pm (\vdim/\curr[\mix])
		\, \minmax(\curr + \curr[\mix]\curr[\seed])
		\, \curr[\seed]
	\notag\\
\next
	&= \curr
		+ \curr[\step] \curr[\signal]
\end{aligned}
\end{equation}
where
$\curr[\mix] \searrow 0$ is a vanishing ``sampling radius'' parameter,
$\curr[\seed]$ is drawn uniformly at random from the composite basis $\seeds = \bvecs_{\minvars} \cup \bvecs_{\maxvars}$ of $\points = \minvars \times \maxvars$,
and
the ``$\pm$'' sign is equal to $-1$ if $\curr[\seed]\in\bvecs_{\minvars}$ and $+1$ if $\curr[\seed]\in\bvecs_{\maxvars}$.
Viewed this way, the interpretation of \eqref{eq:SPSA} as a \acl{RM} method is immediate;
furthermore, a straightforward calculation (that we defer to \cref{app:RM-APT}) shows that the sequence of gradient estimators $\curr[\signal]$ in \eqref{eq:SPSA} has $\curr[\bbound] = \bigoh(\curr[\mix])$ and $\curr[\sdev]^{2} = \bigoh(1/\curr[\mix]^{2})$.
\endenv
\end{algo}
\smallskip

Further examples that can be cast in the \ac{RM} framework include
the negative momentum method \cite{gidel2019negative},
generalized \ac{OG} schemes \cite{mokhtari2019unified},
the Chambolle-Pock algorithm \cite{chambolle2011first},
the ``prediction method'' of \citet{yadav2017stabilizing},
and
centripetal acceleration \cite{peng2020training};
the analysis is similar and we omit the details.
Certain scalable second-order methods can also be viewed as \ac{RM} schemes, but the driving vector field $\vecfield$ is no longer the gradient field of $\minmax$;
we discuss this in the supplement.

\subsection{Alternating updates and moving averages}

There are two extremely common heuristics for practitioners in applying min-max algorithms to real applications: alternating and averaging.
An \emph{alternating} algorithm for \eqref{eq:minmax} updates the $\minvar$ and $\maxvar$ variables sequentially (instead of simultaneously as in \cref{subsec:specific}). An \emph{averaged} algorithm takes the next state as a convex combination of $\curr$ and $\next$ in \eqref{eq:RM}, \cf \cite{karras2018progressive}.

An important feature of our framework is that it captures alternating and averaged algorithms in a seamless manner.
Indeed, introducing alternating updates or a moving average in \ac{RM} schemes results in another \ac{RM} scheme:

\begin{lemma}
\label{lem:averaging}
Let $
\next
	= \curr
		+ \curr[\step] \bracks{ \vecfield(\curr) +  \curr[\error] }$ be an \ac{RM} scheme where $\curr[\error] = \curr[\noise]+ \curr[\bias]$ as in \eqref{eq:error}. Then its $\alpha$-averaged version (where $0 < \alpha < 1$), defined as
\begin{equation}
\label{eq:avg-RM}
\tag{avg-RM}
\begin{alignedat}{2}
\next[\state']
	&=
		\curr[\state]
		+ \curr[\step] \bracks{ \vecfield(\curr[\state]) + \curr[\error]}, 
	\\
\next[\state]
	&=
		\alpha \next[\state'] + (1-\alpha)\curr
\end{alignedat}
\end{equation}
is also an \ac{RM} scheme: $
\next
	= \curr
		+ \alpha\curr[\step] \bracks{ \vecfield(\curr) +  \curr[\error] }$.
\end{lemma}

\begin{remark}
\cref{lem:averaging} can be easily adapted to the scenario where one only averages either the $\curr[\minstate]$ or $\curr[\maxstate]$ variable.
\end{remark}


\begin{figure*}[th!]
\centering
\footnotesize
\includegraphics[height=.3\textwidth]{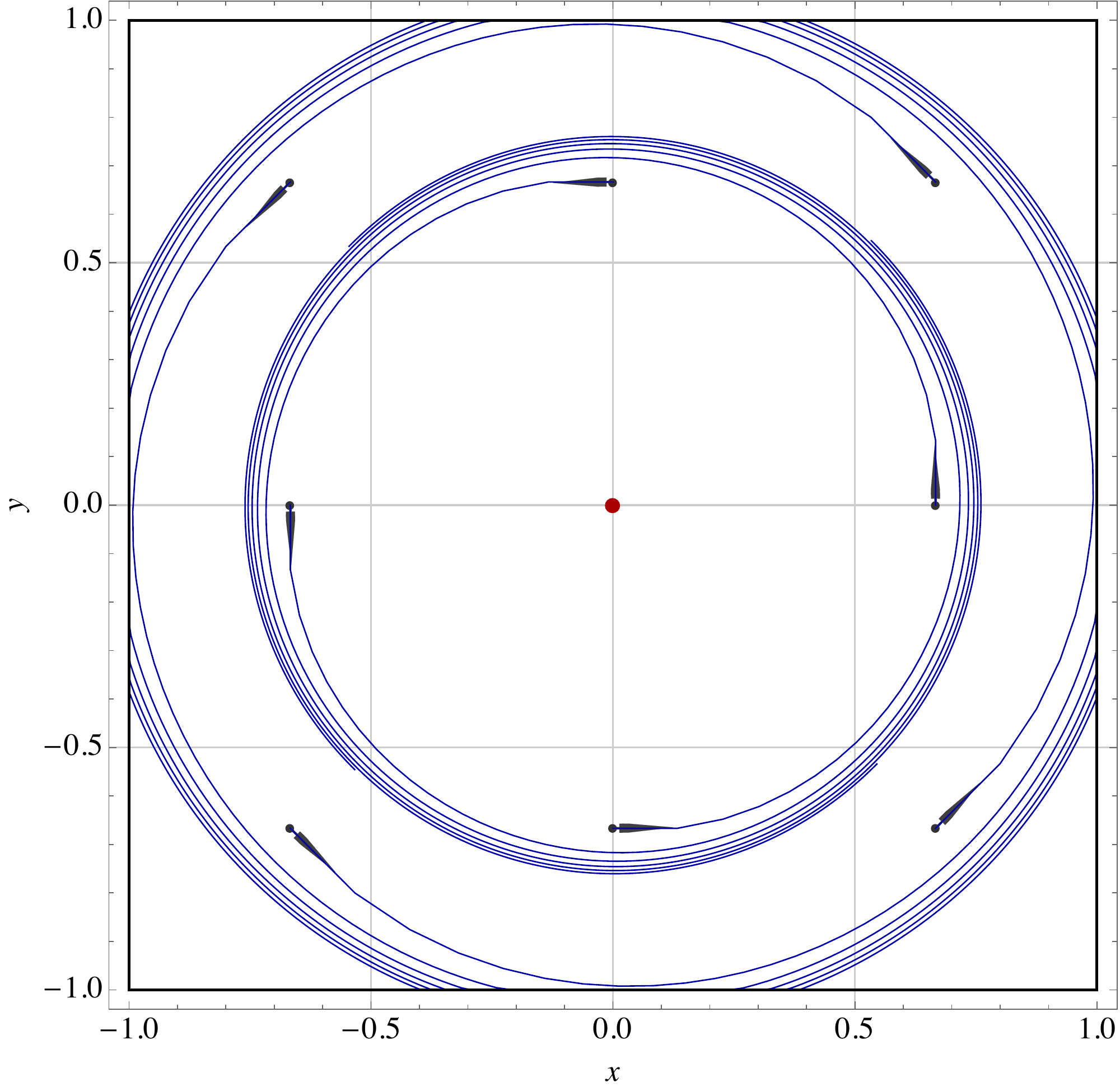}
\hfill
\includegraphics[height=.3\textwidth]{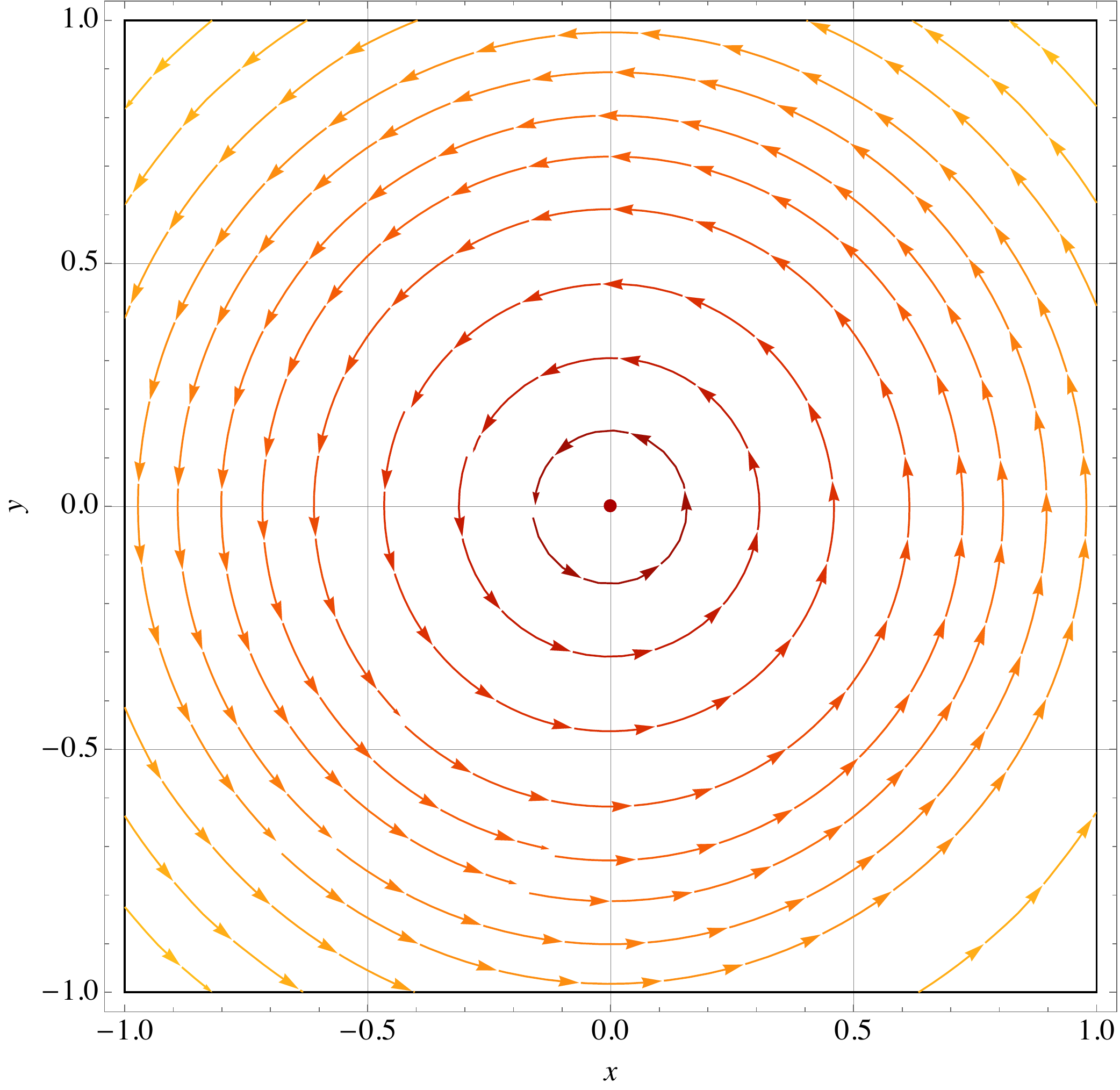}
\hfill
\includegraphics[height=.3\textwidth]{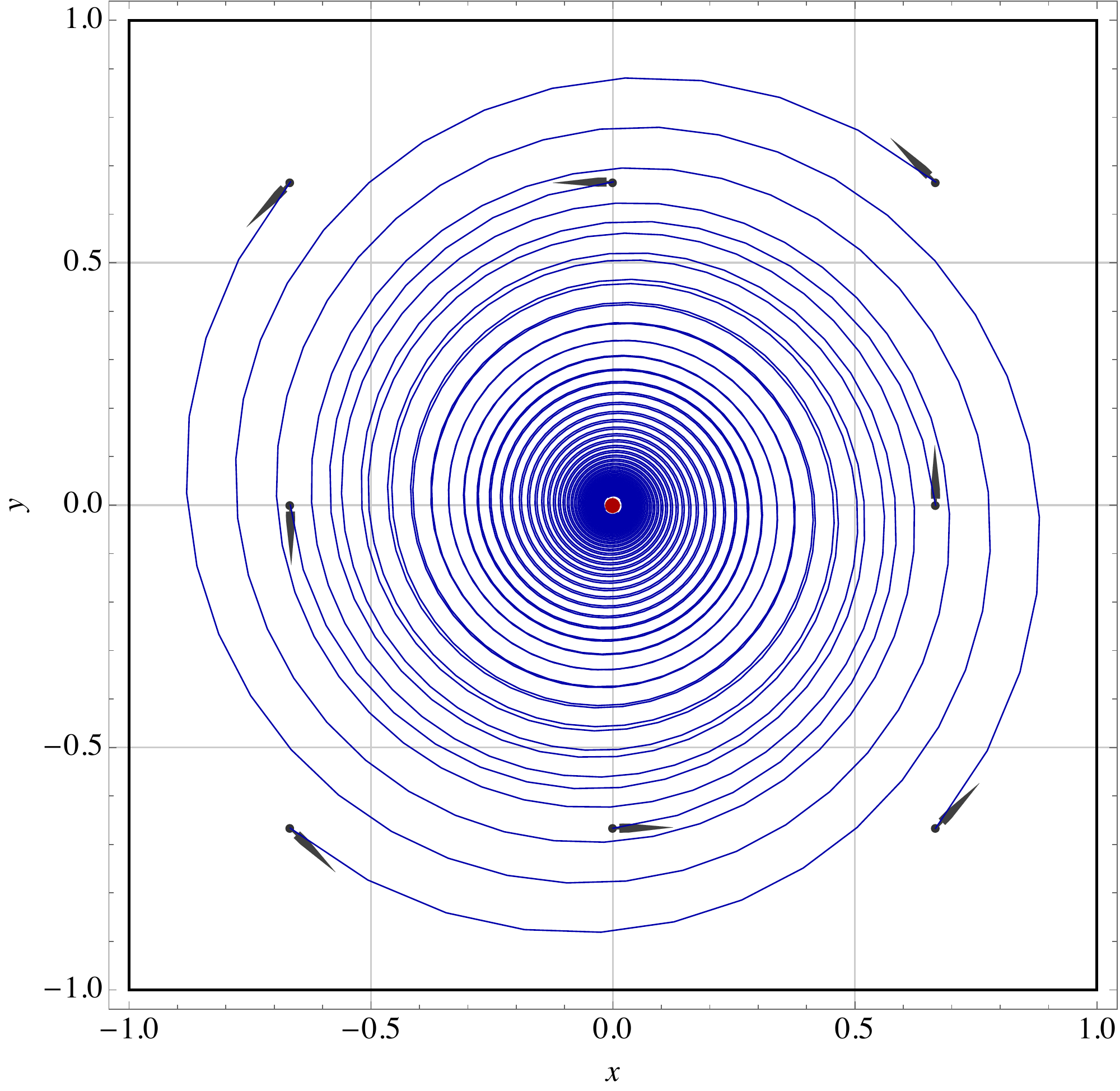}%
\caption{Comparison of different RM schemes for bilinear games $\minmax(\minvar,\maxvar) = \minvar\maxvar$, $\minvar,\maxvar\in\R$.
~From left to right:
(\itshape a\upshape)
\acl{GDA};
(\itshape b\upshape)
the mean dynamics \eqref{eq:MD};
(\itshape c\upshape)
\acl{EG}.}
\label{fig:bilinear}
\end{figure*}


\begin{lemma}
\label{lem:alternating}
Let $
\next
	= \curr
		+ \curr[\step] \bracks{ \vecfield(\curr) +  \curr[\error] }$ be an \ac{RM} scheme where $\curr[\error] = \curr[\noise]+ \curr[\bias]$ as in \eqref{eq:error}. Then its alternating version, defined as
\begin{equation}
\label{eq:alt-RM}
\tag{alt-RM}
\begin{alignedat}{2}
\next[\minstate]
	&=
		\curr[\minstate]
		+ \curr[\step] &&\bracks{ \vecfield_{\minvar}(\curr[\minstate],\curr[\maxstate]) + \error_{\minvar,\run}}, 
	\\
\next[\maxstate]
	&=
		\curr[\maxstate]
		+ \curr[\step] && \bracks{ \vecfield_{\maxvar}( \next[\minstate],\curr[\maxstate]) + \error_{\maxvar,\run}}, 
\end{alignedat}
\end{equation}
is also an \ac{RM} scheme: $
\next
	= \curr
		+ \curr[\step] \bracks{ \vecfield(\curr) +  \curr[\noise] + \curr[\bias'] }$ where
\begin{equation} \nonumber
\curr[\bias'] = \curr[\bias]+  \begin{bmatrix}
0 \\
\vecfield_{\maxvar}( \next[\minstate],\curr[\maxstate]) - \vecfield_{\maxvar}( \curr[\minstate],\curr[\maxstate])
\end{bmatrix}.
\end{equation}
\end{lemma}

\begin{remark}\label{rmk:k1k2}
One can easily generalize \cref{lem:alternating} to the ``$(k_1,k_2)$-\ac{RM} schemes'' where one performs $k_1$ updates for $\minvar$ and then $k_2$ updates for $\maxvar$ (here $k_1,k_2 \in \N$ are arbitrary but fixed). The resulting scheme will still be an \ac{RM} scheme.
In particular, our framework captures the popular $(k_1,k_2) = (1,5)$ variant of \eqref{eq:SGDA} used in the seminal works of \citet{GPAM+14} and \citet{arjovsky2017wasserstein}.
In view of \crefrange{lem:averaging}{lem:alternating}, \cref{rmk:k1k2}, and a simple calculation (see \eqref{eq:bias-order}), \emph{all of our results on first-order methods (\eg \crefrange{alg:SGDA}{alg:PEG}) apply also to their averaging/alternating and the more general $(k_1,k_2)$ versions.}
\end{remark}


\section{Convergence analysis}
\label{sec:analysis}

\subsection{Overview: Continuous \vs discrete time}
\label{sec:cont2disc}

The key in providing a unified treatment of all algorithms in \cref{sec:algorithms} is the reduction of \eqref{eq:RM} to the \emph{mean dynamics}
\begin{equation}
\label{eq:MD}
\tag{MD}
\dotorbit{\ctime}
	= \vecfield(\orbit{\ctime}).
\end{equation}
To see why \eqref{eq:MD} can capture the limiting behavior of a vast family of \ac{RM} schemes beyond GDA, let us illustrate the high-level intuition on the deterministic version of \cref{alg:SEG} ($\curr[\noise]=0$).

Since $\minmax$ and $\vecfield$ are assumed to be Lipschitz (say with constants $\vbound$ and $\lips$), we see that the bias term in \cref{alg:SEG} satisfies
\begin{multline}
\norm{\curr[\bias]}
	= \norm{ \vecfield(\lead) - \vecfield(\curr)}
	\leq \lips \norm{\lead - \curr}
	= \curr[\step] \lips \norm{\vecfield(\curr)} \leq  \curr[\step] \lips \vbound = \bigoh( \curr[\step]). \nonumber
\end{multline}
As a result, we can rewrite \cref{alg:SEG} as
\begin{equation}\label{eq:approx-PPM}
\frac{\next-\curr}{\curr[\step]} =  \vecfield(\curr) + \mathcal{O}(\curr[\step]).
\end{equation}If $\curr[\step]\searrow0$, 
we should then expect \eqref{eq:approx-PPM} to converge to \eqref{eq:MD}.
More generally, if the error term $\curr[\error]$ in \eqref{eq:RM} is sufficiently well-behaved, we should expect the iterates of \eqref{eq:RM} and the solutions of \eqref{eq:MD} to eventually come together.

Connecting \eqref{eq:RM} to \eqref{eq:MD} has proved very fruitful when the latter comprises a \emph{gradient system}, \ie $\vecfield = -\nabla\obj$ for some (possibly non-convex) $\obj\from\points\to\R$: Modulo mild assumptions, the systems \eqref{eq:RM} and \eqref{eq:MD} are known to both converge to the critical set of $\obj$ \citep{Lju77,KC78,BMP90,KY97,BT00}.

On the other hand, bona fide min-max problems are considerably more involved.
The most widely known illustration is given by the bilinear objective $\minmax(\minvar,\maxvar) = \minvar \maxvar$:
in this case (see \cref{fig:bilinear}), the trajectories \eqref{eq:MD} comprise periodic orbits of perfect circles centered at the origin (the unique critical point of $\minmax$).
However, the behavior of different \ac{RM} schemes can vary wildly, even in the absence of noise ($\sdev=0$):
trajectories of \eqref{eq:SGDA} spiral outwards, each converging to an initialization-dependent periodic orbit;
instead, \eqref{eq:SEG} trajectories spiral inwards, eventually converging to the solution $\sol=(0,0)$.

This particular difference between gradient and extra-gradient schemes has been well-documented in the literature \citep{DISZ18,GBVV+19,MLZF+19}.
More pertinent to our theory, it also raises several key questions:
\begin{enumerate}
\itshape
\item
What is the precise link between \ac{RM} methods and the mean dynamics \eqref{eq:MD}?
\item \label[item]{itm:discrepancy-RM}
When does \eqref{eq:MD} yield accurate predictions for the long-run behavior of an \ac{RM} method?
\end{enumerate}
Below, we devote \crefrange{sec:SA}{sec:applications} to the first question, and \cref{sec:convergence} to the second.

\subsection{Connecting \eqref{eq:RM} to \eqref{eq:MD}}
\label{sec:SA}

We begin by introducing a measure of ``closeness'' between the iterates of \eqref{eq:RM} and the solution orbits of \eqref{eq:MD}.
To do so, let $\curr[\efftime] = \sum\nolimits_{\runalt=\start}^{\run} \iter[\step]$ denote the ``effective time'' that has elapsed at the $\run$-th iteration of \eqref{eq:RM}, and define the continuous-time interpolation $\apt{\ctime}$ of $\curr$ as
\begin{equation}
\label{eq:interpolation}
\apt{\ctime}
	= \curr
		+ \frac{\ctime - \curr[\efftime]}{\next[\efftime] - \curr[\efftime]} (\next - \curr)
\end{equation}
for all $\ctime \in [\curr[\efftime],\next[\efftime]]$, $\run\geq\start$.
To compare $\apt{\ctime}$ to the solution orbits of \eqref{eq:MD},
we will further consider the \emph{flow} $\flowmap\from\R_{+}\times\points\to\points$ of \eqref{eq:MD}, 
which is simply the orbit of \eqref{eq:MD} at time $\ctime\in\R_{+}$ with an initial condition $\point(0) = \point \in \points$.
We then have the following notion of ``asymptotic closeness'':

\begin{definition}
\label{def:APT}
$\apt{\ctime}$ is an \acdef{APT} of \eqref{eq:MD} if, for all $\horizon > 0$, we have:
\begin{equation}
\label{eq:APT}
\txs
\lim_{\ctime\to\infty}
	\sup_{0 \leq h \leq \horizon}
		\norm{\apt{\ctime + h} - \flow[h][\apt{\ctime}]}
	= 0.
\end{equation}
\end{definition}
This comparison criterion is due to \citet{BH96} and it plays a central role in our analysis.
In words, it simply posits that $\apt{\ctime}$ eventually tracks the flow of \eqref{eq:MD} with arbitrary accuracy over windows of arbitrary length;
as a result, if $\curr$ is an \ac{APT} of \eqref{eq:MD}, it is reasonable to expect its behavior to be closely correlated to that of \eqref{eq:MD}.

Our first result below makes this link precise.
Consider an \ac{RM} scheme which satisfies
\begin{alignat}{2}
\label[assumption]{asm:bias}
\tag{A1}
\txs
\curr[\bbound] 
	\to 0
	\;\;
	\text{\as}
	\quad
	\text{and}
	\quad
\sum_{\run=\start}^\infty \exof{\curr[\step]\curr[\bbound]}  < \infty
	\\
\label[assumption]{asm:noise}
\tag{A2}
\txs
\sum\nolimits_{\run=\start}^{\infty} \exof{\curr[\step]^{2} (1+\curr[\bbound]^2+\curr[\sdev]^{2})}
	< \infty
	\hspace{1em}
\end{alignat}
We then have:


\begin{restatable}{theorem}{APT}
\label{thm:APT}
Suppose that 
\crefrange{asm:bias}{asm:noise} hold.
Then $\curr$ is an \ac{APT} of \eqref{eq:MD} \acs{wp1}.
\end{restatable}


\subsection{Applications and examples}
\label{sec:applications}

Of course, applying \cref{thm:APT} to a specific algorithm (\eg as in \cref{sec:algorithms}) would first require verifying \crefrange{asm:bias}{asm:noise}.
However, even though the noise $\err(\point;\seed)$ in \eqref{eq:SFO} is assumed zero-mean and finite-variance, this \emph{does not imply} that the error term $\curr[\error] = \curr[\noise] + \curr[\bias]$ in \crefrange{alg:PPM}{alg:SPSA} enjoys the same guarantees.
For example, the \ac{RM} representation of \crefrange{alg:PPM}{alg:PEG} has non-zero bias,
while \cref{alg:SPSA} has non-zero bias \emph{and} unbounded variance (the latter behaving as $\bigoh(1/\curr[\mix]^{2})$ with $\curr[\mix] \to 0$). 

In the following proposition we prove that \crefrange{alg:SGDA}{alg:SPSA} generate \aclp{APT} of \eqref{eq:MD} for the typical range of hyperparameters used to ensure almost sure convergence of stochastic first-order methods.


\begin{restatable}{proposition}{algoAPT}
\label{prop:APT}
Let $\curr$ be a sequence generated by any of the \crefrange{alg:SGDA}{alg:SPSA}.
Assume further that:
\begin{enumerate}
[\itshape a\upshape),leftmargin=.3in]
\item
For first-order methods \textpar{\crefrange{alg:SGDA}{alg:PEG}}, the algorithm is run with \ac{SFO} feedback satisfying \eqref{eq:err-base} and a step-size $\curr[\step]$ such that $A/\run \leq \curr[\step] \leq B/ \sqrt{\run(\log\run)^{1+\eps}}$ for some $A,B,\eps>0$.
\item
For zeroth-order methods \textpar{\cref{alg:SPSA}},
the algorithm is run with parameters $\curr[\step]$ and $\curr[\mix]$ such that
$\lim_{\run} (\curr[\step] + \curr[\mix]) = 0$,
$\sum_{\run} \curr[\step] = \infty$,
and
$\sum_{\run} \curr[\step]^{2} / \curr[\mix]^{2} < \infty$
\textpar{\eg $\curr[\step] = 1/\run$, $\curr[\mix] = 1/\run^{1/3}$}.
\end{enumerate}
Then $\curr$ is almost surely an \ac{APT} of \eqref{eq:MD}. 
\end{restatable}


\subsection{The limit sets of \ac{RM} schemes}
\label{sec:convergence}

The \ac{APT} results in \crefrange{sec:SA}{sec:applications} can be heuristically interpreted as:
``\ac{RM} schemes eventually behave as some orbits of \eqref{eq:MD}.''
We now further ask:
What are \emph{the} {candidate limit orbits} of \eqref{eq:MD} for \ac{RM} schemes?

To shed some light on the question, let us recall that, in non-convex \emph{minimization} problems,
\ac{SGD} enjoys the following properties:
\begin{enumerate}[left=\parindent,label={(\Roman*)}]
\item
\label[item]{itm:sgd-critical}
\ac{SGD} converges to the function's set of critical points \citep{Lju77,BT00}.
\item
\label[item]{itm:sgd-avoidance}
\ac{SGD} avoids unstable critical points \citep{Pem90,GHJY15, mertikopoulos2020almost}.
\end{enumerate}
This leads to the following ``law of the excluded middle'':
generically, the only solution candidates left for \ac{SGD} are stable critical points, \ie the local minimizers of the problem's minimization objective.

In the remaining of this section, we will assimilate \cref{itm:sgd-critical,itm:sgd-avoidance} in the context of \ac{RM} schemes applied to \eqref{eq:minmax}.

\subsubsection{The long-run limit of \ac{RM} schemes}

We first focus on generalizing \cref{itm:sgd-critical} for min-max optimization. To proceed, recall first that critical points alone cannot capture the broad spectrum of algorithmic behaviors when \eqref{eq:MD} is not a gradient system:
already in \cref{fig:bilinear} we see a critical point surrounded by \emph{spurious} periodic orbits.
In addition, in dynamical systems many other spurious convergence phenomena are known, such as homoclinic loops, limit cycles, or chaos.
To account for this considerably richer landscape, we will need some definitions from the theory of dynamical systems.

\begin{definition}[\citealp{Ben99}]
\label{def:ICT}
Let $\set$ be a nonempty compact subset of $\points$.
Then:
\begin{enumerate}
[\itshape a\upshape),leftmargin=.3in]
\item
$\set$ is \emph{invariant} if $\flow[\ctime][\set] = \set$ for all $\ctime\in\R$.
\item \label{item:b}
$\set$ is \emph{attracting} if it is invariant and there exists a compact neighborhood $\cpt$ of $\set$ such that $\lim_{\ctime\to\infty} \dist(\flow,\set) = 0$ uniformly in $\point\in\cpt$.
\item
$\set$ is \acdef{ICT} if it is invariant and $\flowmap\vert_{\set}$ admits no proper attractors in $\set$.
\end{enumerate}
\end{definition}

\begin{remark*}
Equivalently, \ac{ICT} sets can be viewed as ``minimal connected periodic orbits up to arbitrarily small numerical errors'', \cf \citet[Prop.~5.3]{Ben99}.
The definition above is more convenient to work with because it provides the key insights in \cref{sec:RMsame?} below.
\end{remark*}

Our next result shows that,
{with probability $1$, all limit points of \eqref{eq:RM} lie in these ``approximate periodic orbits'':

\begin{restatable}{theorem}{ICT}
\label{thm:ICT}
If \crefrange{asm:bias}{asm:noise} hold,
then $\curr$ converges almost surely to an \ac{ICT} set of $\minmax$.
\end{restatable}

\begin{corollary}
\label{cor:ICT}
Let $\curr$ be a sequence generated by any of the \crefrange{alg:SGDA}{alg:SPSA} with parameters as in \cref{prop:APT}.
Then $\curr$ converges almost surely to an \ac{ICT} set of $\minmax$.
\end{corollary}

\subsubsection{Avoidance of unstable points and sets}

Our next result provides an avoidance result for \ac{RM} schemes in min-max optimization.
In analogy with function minimization problems, we will focus on unstable \emph{invariant sets} of \eqref{eq:MD}, \ie invariant sets that admit a nontrivial unstable manifold (for an in-depth discussion and precise definition, see \citealp{Shu87} and \cref{app:avoid-repel}).

In generic minimization problems, these are precisely the sets of strict saddle points of the function being minimized.
However, since general min-max problems do \emph{not} comprise a gradient system, \eqref{eq:MD} could exhibit a plethora of unstable sets, not containing any stationary points of $\minmax$ (\eg periodic orbits, heteroclinic networks, etc.).
On account of the above, our result below is stated in terms of invariant \emph{sets} \textendash\ and not only points.
For convenience, we will assume that $\vecfield$ is $C^{2}$ and $\curr[\step]$ is as in \cref{prop:APT}.

\begin{restatable}{theorem}{repeller}
\label{thm:repeller}
Let $\cpt$ be an unstable invariant set of
\eqref{eq:MD}
\textpar{which trivially includes unstable periodic orbits and unstable critical points}.
Assume further that the noise in \eqref{eq:RM} satisfies: 
\begin{enumerate*}[label=\upshape(\itshape\roman*\hspace*{1pt}\upshape)]
\item
\label[assumptionenum]{asm:bounded-noise} 
$\sup_{\run} \norm{\curr[\noise]} < \infty$ \acs{wp1};
and
\item \label[assumptionenum]{asm:unif-exciting}
$\inf_{z:\norm{z}=1} \exof*{ \inner{\curr[\noise]}{z}_{+} \given \curr[\filter] } > 0$.
\end{enumerate*}
Then $\curr$ generated by any of the \crefrange{alg:SGDA}{alg:PEG} satisfies
\begin{equation}\nonumber
\txs
\prob \left({  \lim_{ \run \to \infty } \dist \parens{ \curr[\state], \cpt } = 0 } \right) = 0.
\end{equation}
\end{restatable}

\begin{remark}
We note that \cref{asm:bounded-noise,asm:unif-exciting} above are standard in the literature for avoidance results of \ac{SGD} \cite{Pem90,Ben99,mertikopoulos2020almost}, and are significantly lighter than other ``isotropic noise'' assumptions that are common in the literature \citep{GHJY15}.
Specifically, even though \cref{asm:unif-exciting} looks somewhat obscure, it only posits that the noise is not degeneratively equal to zero along certain directions in space;
for a more detailed discussion, see \cref{app:avoid-repel}.
We also stress that neither of these assumptions is required for the rest of our paper.
\end{remark}

\subsubsection{When do \ac{RM} schemes behave the same?} 
\label{sec:RMsame?}

So far, we have successfully generalized \cref{itm:sgd-critical,itm:sgd-avoidance} to the context of \eqref{eq:minmax} as follows:%
\footnote{To see why this is really a generalization, simply note that the only \ac{ICT} sets of $\vecfield = -\grad\obj$ are connected critical points of $\obj$; \cf \cref{prop:gradient}.} 
\begin{enumerate}[(I-SP)]
\item \label[item]{itm:minmax-critical}\ac{RM} schemes always converge to {\ac{ICT} sets}, and
\item \label[item]{itm:minmax-avoidance}\ac{RM} schemes always avoid invariant sets.
\end{enumerate}
 
Nonetheless, \cref{itm:minmax-critical,itm:minmax-avoidance} still fail to explain the distinct behaviors of \ac{RM} schemes in bilinear objectives: Why does \ac{SGDA} converge only to periodic orbits, while deterministic \ac{SEG} only to critical points?
Or, more generally, 
\begin{quoting}
\itshape
\centering
Are different \ac{RM} schemes more likely to exhibit different convergence topologies \textendash\ \eg cycles \vs critical points \textendash\ in generic min-max problems?
\end{quoting}

Our next result takes a closer look at \emph{attracting} \ac{ICT} sets and provides a generically negative answer to this question. 
To set the stage, suppose we want to apply \cref{itm:minmax-critical} to the bilinear objective $\minmax(\minvar,\maxvar)= \minvar\maxvar$.
{Stricto sensu, \cref{itm:minmax-critical} does not apply in this case since $\minmax$ is not Lipschitz.}
However, \cref{fig:bilinear}(\itshape b\upshape) shows (and we rigorously prove in \cref{app:proof-ICT}) that \emph{any} tuple $(\minvar,\maxvar) \in \R^{2}$ belongs to an \ac{ICT} set of $\minmax$, so \cref{thm:ICT} holds trivially.
This in turn implies that \emph{the only attractor for $\minmax$ is trivially the whole space $\R^2$}, since \cref{def:ICT}-\ref{item:b} is never satisfied for any $\set \subsetneq \R^2$.

Importantly, the celebrated \emph{Kupka-Smale theorem} \cite{kupka1963contributiona, smale1963stable} asserts that systems with degenerate periodic orbits (such as bilinear games) occur ``almost never'' in the Baire category sense.
More precisely, an arbitrarily small perturbation can fundamentally destroy the topological properties of their \ac{ICT} sets and give rise to proper, non-trivial attractors; \cf \cref{ex:bilinear}. 
In contrast, systems with nontrivial attractors are known to be robust under perturbations \cite{Shu87}, and our final result in the section shows that it is precisely the \emph{existence of nontrivial attractors} that makes the discrepancy of \ac{RM} schemes disappear, at least locally.

\begin{restatable}{theorem}{attract}
\label{thm:attract}
Let $\set$ be an attractor of \eqref{eq:MD} and fix some confidence level $\conf>0$.
If $\curr[\step]$ is small enough and \crefrange{asm:bias}{asm:noise} hold,
there exists a neighborhood $\nhd$ of $\set$, independent of $\conf$, such that
\(
\probof{ \text{$\curr$ converges to $\set$}}
	\geq 1-\conf
\)
if $\init\in\nhd$.
\end{restatable}

\begin{corollary}
\label{cor:attract}
Let $\curr$ be a sequence generated by any of the \crefrange{alg:SGDA}{alg:SPSA} with sufficiently small $\curr[\step]$ satisfying the conditions of \cref{prop:APT}.
If $\init\in\nhd$, then $\probof{ \text{$\curr$ converges to $\set$}} \geq 1-\conf$.
\end{corollary}

In short, \cref{thm:attract} asserts that any non-degenerate \ac{ICT} set dictates the local convergence of \emph{all} \ac{RM} schemes under the general \crefrange{asm:bias}{asm:noise}. 

On a positive note, since the Hartman-Grobman Theorem \cite{robinson1998dynamical} implies that all critical points of $\minmax$ with $\Re \braces{\lambda(\jmat(\sol)) } <0$ for all engenvalues $\lambda$ are attractors of \eqref{eq:MD}, \cref{thm:attract} immediately yields:
\begin{corollary}
\label{cor:stable}
Let $\sol$ be a critical point of $\minmax$ such that $\Re \braces{\lambda(\jmat(\sol))} <0$ for all engenvalues of $\jmat(\sol)$. Then all \ac{RM} schemes satisfying \crefrange{asm:bias}{asm:noise} locally converge to $\sol$ with high probability.
\end{corollary}

\cref{cor:stable} generalizes the local convergence of deterministic \ac{SGDA} and \ac{SEG} studied by \citet{daskalakis2018limit}. It also extends \citep[Theorem 5]{HIMM19} from \eqref{eq:PEG} to all generalized \ac{RM} schemes.

On the flip side, however, \cref{thm:attract} also bears an undesirable consequence:
it implies that many \ac{RM} schemes designed to improve \ac{SGDA} (\eg \crefrange{alg:PPM}{alg:PEG}) may in fact be trapped by \emph{spurious \ac{ICT} sets} in exactly the same way as \ac{SGDA}.
Thus, even though many of these algorithms have been motivated by their appealing properties in bilinear games, it is not clear whether they offer any significant advantages beyond the convex-concave case.
We examine this issue in detail in the next section.


\begin{figure}[t]
\centering
\footnotesize
\includegraphics[height=.48\textwidth]{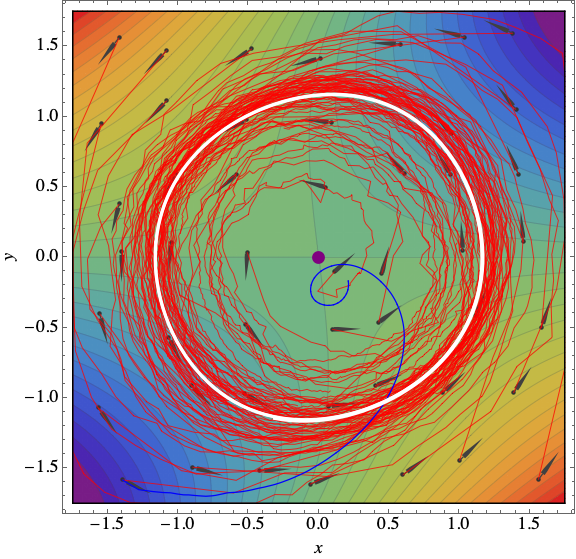}
\hfill
\includegraphics[height=.48\textwidth]{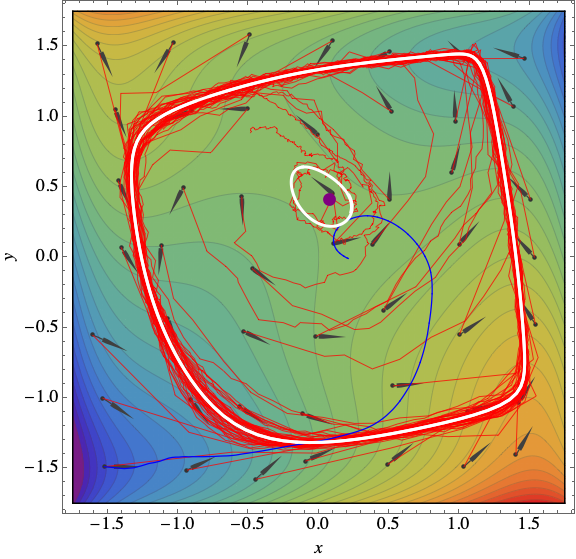}%
\caption{Spurious limits of min-max optimization algorithms. 
From left to right:
(\itshape a\upshape)
\eqref{eq:SEG} for \eqref{eq:perturbed-bilinear} with $\eps = 0.01$;
(\itshape b\upshape)
``forsaken solutions'' of \eqref{eq:SEG};
The red curves present trajectories with different initialization;
non-critical \ac{ICT} sets are depicted in white;
the blue curves represent an time-averaged sample orbit.
}
\label{fig:spurious}
\end{figure}

\section{Spurious attractors: Illustrations and examples}
\label{sec:spurious}

In this last section, we provide a range of simple examples that exhibit \emph{spurious attractors} \textendash\  \ie attractors that consist entirely of non-critical points.
For illustration purposes, we focus on the simple case $\minvars = \maxvars = \R$ with polynomial objectives.
In doing this, our goal is to highlight a number of issues that can arise in min-max optimization problems;
whether limit cycles of this type occur in actual large-scale experiments \textendash\ \eg in \acsp{GAN} \textendash\ is an open research question \citep{letcher2020impossibility}.


\begin{example}
[Almost bilinear $\not\approx$ bilinear, instability $\not\approx$ escape]
\label{ex:bilinear}
Consider an arbitrarily small perturbation of a bilinear game:
\begin{equation}
\label{eq:perturbed-bilinear}
\minmax(\minvar,\maxvar)
	= \minvar\maxvar
		+ \eps \perturb(\maxvar),
\end{equation}
where $\eps>0$ and $\perturb(\maxvar) = \frac{1}{2}\maxvar^2 - \frac{1}{4}\maxvar^4$.
There is an unstable critical point at the origin;
further, \cref{lem:cycles-in-almost-bilinear} asserts, for small $\eps$, the existence of an \emph{attracting} \ac{ICT} set $\set$ in a neighborhood of the circle $\setdef{\point}{\norm{\point}^2 = 4/3}$.
By \cref{cor:attract}, any \ac{RM} scheme of \cref{sec:algorithms} thus gets trapped by $\set$; see \cref{fig:spurious}(\itshape a\upshape) for an illustration for \eqref{eq:SEG}.

This example brings two issues of existing studies to light. First, it shows that ``almost bilinear games'' can still trap many methods for solving exact bilinear games.
Second, in contrast to minimization problems, the region around an unstable critical point can in fact be fully stable.
Thus, one has to be careful when interpreting algorithms that ``locally avoid unstable critical points'', since they might be incapable of escaping their neighborhoods.
\endenv
\end{example}

\begin{example}
[``Forsaken'' $\min\max$ solutions]\label{ex:forsaken-stable}
Suppose we apply \crefrange{alg:SGDA}{alg:SPSA} to the objective
\begin{equation}\label{eq:forsaken-stable}
\minmax(\minvar,\maxvar) 
	= \minvar(\maxvar-0.45)
		+ \perturb(\minvar)
		- \perturb(\maxvar)
\end{equation}
where $\perturb(\point) = \frac{1}{4} \point^2- \frac{1}{2}\point^4 + \frac{1}{6} \point^6$.
This problem has a desirable $(\minsol, \maxsol) \simeq (0.08, 0.4)$.
However, as we show in \cref{app:spurious-forsaken}, there exist \emph{two} spurious limit cycles that do not contain \emph{any} critical point of $\minmax$. Worse, the limit cycle closer to the solution is \emph{unstable} and repels any trajectory that comes close to the solution;
see \cref{fig:spurious}(\itshape b\upshape) for an illustration for \eqref{eq:SEG}.
As a result, the ``shielded'' solution is highly unlikely to be discovered by existing algorithms, even though it is perfectly stable.
\endenv
\end{example}

We conclude the paper by further examining several important settings that are not covered by our theory:




\begin{enumerate}
\item
Instead of the ``moving average'' in \cref{lem:averaging}, one can take the \emph{ergodic average} ($\curr[\state']= \frac{1}{n}\insum_{\runalt = 1}^\run \state_\runalt$) as is customary in convex-concave problems \citep{nemirovski2004prox,juditsky2011solving}. We plot one such trajectory in \cref{fig:spurious} (the blue curves). Evidently, we see that ergodic average can force the algorithms to halt at non-critical points, and this convergence is by no means min-max optimal.

\item
Many recent works attempt to address the cycling issues of min-max algorithms via incorporating \emph{second-order} oracles.
For completeness, we also study a range of popular {second-order} methods in \cref{app:perturb}.
Our analysis shows that these algorithms suffer similar symptoms as first-order schemes in our examples, \cf \crefrange{fig:CO}{fig:algs}.

\item
In addition to the diminishing step-size policies studied here, another common strategy in practice is to simply set $\curr[\step]$ to a \emph{constant step-size}. While our analysis does not cover this setting, there exist several techniques in stochastic approximation to boost from our ``almost surely'' statements for $\curr[\step]\searrow0$ to \emph{concentration} or \emph{high-probability} results when $\curr[\step] \equiv \step$ is small \cite{kushner1981asymptotic, KY97, Bor08}. 

For completeness, in \cref{app:const-step} we examine various constant step-size \ac{RM} schemes applied to \eqref{eq:perturbed-bilinear} and \eqref{eq:forsaken-stable}. The outcome coincides with our intuition that these schemes should concentrate around the spurious attracting \ac{ICT} sets, and hence exhibit similar behaviors as \ac{RM} schemes with $\curr[\step]\searrow0$; see \cref{fig:const-step}.

\item \emph{Adaptive methods} such as Adam \cite{kingma2014adam} are ubiquitous in \acs{GAN} training.
We study such methods in \cref{app:adaptive}:
our results show tha they fail solve the simple objectives \eqref{eq:perturbed-bilinear} and \eqref{eq:forsaken-stable}. Moreover, some methods even show a potentially detrimental tendency of converging to \emph{max-min points}, the exact opposite of desirable solutions; see \cref{fig:adaptive}.
\end{enumerate}

In closing, we should clarify that these illustrations are \emph{not} meant to suggest that the algorithms and practical tweaks discussed above are always doomed, or that they comprise the principal cause of failure in \acs{GAN} training.
However, we do believe that they constitute an important cautionary tale to the effect that, in min-max problems,
\emph{convergence does not imply optimality} \textendash\ or even \emph{stationarity}.

\appendix
\numberwithin{equation}{section}		
\numberwithin{lemma}{section}		
\numberwithin{proposition}{section}		
\numberwithin{theorem}{section}		

\section{Stabilization of \ac{RM} schemes}
\label{app:coercive}

Our aim in this appendix will be to prove the stability of generalized \ac{RM} schemes, namely that \aclp{APT} generated by \eqref{eq:RM} are bounded \acl{wp1}.
The key ingredient in our analysis is the \acf{WAC} condition \eqref{asm:coercive}, which, as we discussed in the main body of our paper, is a relaxation of the standard coercivity requirement
\begin{equation}
\label{eq:coercive}
\lim_{\norm{\point}\to\infty} \frac{\braket{\vecfield(\point)}{\point}}{\norm{\point}}
	= -\infty.
\end{equation}
The hypothesis \eqref{eq:coercive} is a mainstay in the analysis of monotone operators and variational inequalities \citep{BC17,FP03,Phe93}.
Roughly speaking, it states that the ``radial component''
\begin{equation}
\vecfield_{r}(\point)
	= \frac{\braket{\vecfield(\point)}{\point}}{\norm{\point}}
\end{equation}
of $\vecfield(\point)$ grows to $-\infty$ as $\norm{\point}\to\infty$.
In other words, any vector field that satisfies \eqref{eq:coercive} has an inward-pointing component that grows infinitely large for large $\norm{\point}$.

In view of the above, the coercivity assumption \eqref{eq:coercive} suggests that any process that takes successive steps along $\vecfield(\point)$ will be subject to an ``inwards drift'' towards regions with smaller norm, and this drift will be more and more pronounced the farther one moves away from the origin.
On that account, \eqref{eq:coercive} is a natural candidate for showing that \ac{RM} processes based on $\vecfield$ never escape to infinity.
On the other hand, vector fields that do not have a strong radial component \textendash\ such as the bilinear game field $\vecfield(\minvar,\maxvar) = (-\maxvar,\minvar)$ which has $\vecfield_{r}(\minvar,\maxvar) = 0$ \textendash\ are not covered by \eqref{eq:coercive}.
In this regard, the \ac{WAC} condition \eqref{asm:coercive} provides an important relaxation of \eqref{eq:coercive}, because it only posits that the radial component of $\vecfield(\point)$ is asymptotically non-positive \textendash\ or, more simply, that $\vecfield(\point)$ does not have a persistent outward-pointing component.

Before proving the stability of generalized \ac{RM} schemes under \eqref{asm:coercive}, we provide below a series of important examples that satisfy the \ac{WAC} condition \eqref{asm:coercive}:

\begin{enumerate}

\item
\emph{$\vecfield$ satisfies \eqref{eq:coercive}.}
Indeed, in this case, for all $M > 0$, there exists some $R\equiv R(M)$ such that
\begin{equation}
\frac{\braket{\vecfield(\point)}{\point}}{R}
	\leq \frac{\braket{\vecfield(\point)}{\point}}{\norm{\point}}
	\leq - M
	< 0
\end{equation}
whenever $\norm{z} \geq R$, \ie \eqref{asm:coercive} holds

\item
\emph{$\minmax$ is convex-concave and it admits a critical point.}
By shifting the problem's frame of reference if necessary, we can assume without loss of generality that $\sol = 0$ is a critical point of $\minmax$.
Then, with $\minmax$ assumed convex-concave, we readily get $\braket{\vecfield(\point)}{\point} \leq \braket{\vecfield(\sol)}{\point - \sol} = 0$, \ie \eqref{asm:coercive} holds
\end{enumerate}
The first item above justifies the terminology ``\acl{WAC}'';
for a geometric illustration, see \cref{fig:stability} below.

%
%


\begin{figure*}[t]
\centering
\footnotesize
\includegraphics[height=.4\textwidth]{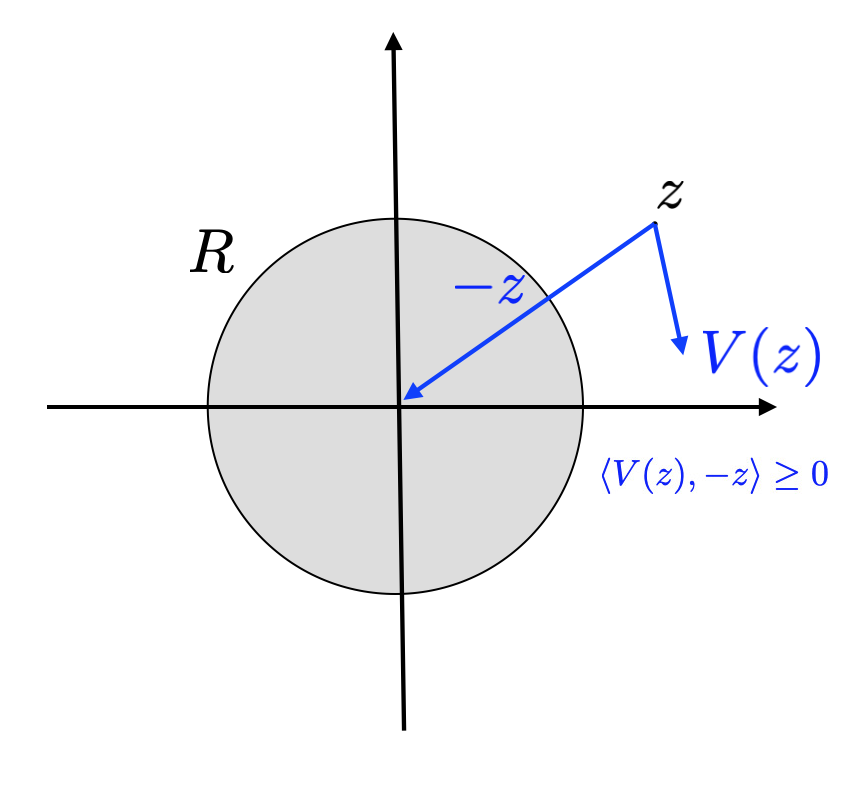}
\caption{Schematic illustration of the \acl{WAC} condition \cref{asm:coercive}. 
}
\label{fig:stability}
\end{figure*}


We now proceed to establish our main stability result for generalized \ac{RM} schemes under the \ac{WAC} condition \eqref{asm:coercive}:

\begin{proposition}
\label{prop:stability}
Suppose that $\vecfield$ satisfies \cref{asm:coercive}.
Then, under \crefrange{asm:bias}{asm:noise}, the sequence $\curr$ generated by \eqref{eq:RM} is bounded \as.
\end{proposition}

\newcommand{\coer}{R}

\begin{proof}[Proof of \cref{prop:stability}]
Our proof hinges on the introduction of a suitable ``energy function'' for \eqref{eq:MD}.
To define it, recall that that $\braket{\vecfield(\point)}{-\point} \geq 0$ whenever $\norm{\point} \geq \coer$.
~Then, with a fair amount of hindsight, fix some $\coef>0$ and let
\begin{equation}
\lyap(\point)
	= \begin{cases}
	0
		&\quad
		\text{if $\twonorm{\point}\leq\coer$},
		\\
	(\twonorm{\point} - \coer)^{2}/2
		&\quad
		\text{if $\coer \leq \twonorm{\point} \leq (1+\coef) \coer$},
		\\
	\coef \coer \twonorm{\point} - \coef(1 + \coef/2) \coer^{2}
		&\quad
		\text{if $(1+\coef)\coer \leq \twonorm{\point}$}.
	\end{cases}
\end{equation}
By a direct calculation, we can verify the following:
\begin{enumerate}
\item
$\lyap$ is continuously differentiable and its gradient is given by $\nabla\lyap(\point) = \phi(\twonorm{\point}/\coer) \, \point$
where
$\phi(u) = 0$ if $u \leq 1$,
$\phi(u) = 1-1/u$ if $1 \leq u \leq 1+\coef$,
and
$\phi(u) = \coef/u$ if $1+\coef \leq u$.
\item
$\lyap$ is negatively correlated to $\vecfield$, \ie $\braket{\vecfield(\point)}{-\nabla\lyap(\point)} \geq 0$ for all $\point\in\points$.
\item
$\lyap$ is $1$-smooth, \ie $\lyap(\pointalt) \leq \lyap(\point) + \braket{\nabla\lyap(\point)}{\pointalt - \point} + (1/2) \twonorm{\pointalt-\point}^{2}$ for all $\point,\pointalt\in\points$.
\end{enumerate}
Then, letting $\curr[\lyap] = \lyap(\curr)$ and $\curr[\phi] = \phi(\twonorm{\curr}/\coer)$, we get
\begin{align}
\label{eq:Lyap-bound1}
\next[\lyap]
	= \lyap(\curr -\curr[\step]\curr[\signal])
	&\leq \lyap(\curr)
		- \curr[\step] \braket{\nabla\lyap(\curr)}{\curr[\signal]}
		+ \frac{\curr[\step]^{2}}{2} \twonorm{\curr[\signal]}^{2}
	\notag\\
	&\leq \curr[\lyap]
		- \curr[\step] \curr[\phi] \braket{\curr[\noise] + \curr[\bias]}{\curr}
	+ \frac{3\curr[\step]^{2}}{2} \bracks[\big]{\twonorm{\vecfield(\curr)}^{2} + \twonorm{\curr[\noise]}^{2} + \twonorm{\curr[\bias]}^{2}},
\end{align}
where the second line follows from the properties of $\lyap$, the definition \eqref{eq:signal} of $\curr[\signal]$, and the Cauchy-Schwarz inequality.
Hence, conditioning on $\curr[\filter]$ and taking expectations, we obtain:
\begin{equation}
\label{eq:Lyap-bound2}
\exof{\next[\lyap] \given \curr[\filter]}
	\leq \curr[\lyap]
		+ \curr[\step] \curr[\phi] \twonorm{\curr} \curr[\bbound]
		+ \tfrac{3}{2} \curr[\step]^{2} \bracks[\big]{\vbound^{2} + \curr[\bbound]^{2} + \curr[\sdev]^{2}},
\end{equation}
where we made a second use of the Cauchy-Schwarz inequality in the term involving $\curr[\bbound]$ and $\vbound$ is the Lipschitz constant of $\minmax$.


To proceed, let $\curr[\eps] = \curr[\step] \curr[\phi] \twonorm{\curr} \curr[\bbound] + (3/2) \curr[\step]^{2} \bracks{\vbound^{2} + \curr[\bbound]^{2} + \curr[\sdev]^{2}}$ denote the ``residual'' term in \eqref{eq:Lyap-bound2}, and consider the auxiliary process $\curr[\aux] = \next[\lyap] + \sum_{\runalt=\run+1}^{\infty} \iter[\eps]$.
By \eqref{eq:Lyap-bound2}, we have $\exof{\curr[\aux] \given \curr[\filter]} \leq \curr[\lyap] + \sum_{\runalt=\run}^{\infty} \curr[\eps] = \prev[\aux]$, \ie $\curr[\aux]$ is a supermartingale relative to $\curr[\filter]$.
By the definition of $\phi$, we further have $\phi(u) \leq \coef/u$, so $\curr[\phi] \twonorm{\curr} \leq \coef\coer$ for all $\run$.
We thus get
\begin{equation}
\sum_{\run=\start}^{\infty} \curr[\eps]
	\leq \coef\coer \sum_{\run=\start}^{\infty} \curr[\step]\curr[\bbound]
		+ \frac{3}{2} \sum_{\run=\start}^{\infty} \curr[\step]^{2} (\vbound^{2} + \curr[\bbound]^{2} + \curr[\sdev]^{2})
\end{equation}
and hence, by \cref{asm:bias,asm:noise}, we conclude that $\exof{\sum_{\run} \curr[\eps]} < \infty$.
This shows that $\exof{\curr[\aux]} \leq \exof{\init[\aux]} < \infty$, \ie $\curr[\aux]$ is uniformly bounded in $L^{1}$.
Hence, by Doob's submartingale convergence theorem \citep[Theorem~2.5]{HH80}, it follows that $\curr[\aux]$ converges \acl{wp1} to some finite random limit $\aux_{\infty}$.
In turn, since $\sum_{\run} \curr[\eps] < \infty$, this implies that $\curr[\lyap] = \prev[\aux] - \sum_{\runalt=\run}^{\infty} \curr[\eps]$ also converges to some (random) finite limit \as.
From this we deduce that $\limsup_{\run} \twonorm{\curr} < \infty$, as claimed.
\end{proof}

\section{Proof of \cref{thm:APT,prop:APT}}
\label{app:APT}

In this appendix, we discuss how the algorithms in \cref{sec:algorithms} fit within the general \acl{SA} framework of \cref{sec:SA}.
Specifically, we prove the general conditions of \cref{thm:APT,prop:APT} which guarantee that \crefrange{alg:SGDA}{alg:SPSA} generate \aclp{APT} of the mean dynamics \eqref{eq:MD}.

\subsection{Generalities and preliminaries}

Before doing so, we will require some background material on \aclp{APT}.
Following \citet{BH96} and \citet{Ben99}, we first recall the definition of the ``effective time'' $\curr[\efftime] = \sum_{\runalt=\start}^{\run} \iter[\step]$ as the time that has elapsed at the $\run$-th iteration of the discrete-time process $\curr[\state]$;
recall also the definition \eqref{eq:interpolation} of the continuous-time interpolation $\apt{\ctime}$ of $\curr$ as
\begin{equation}
\tag{\ref{eq:interpolation}}
\apt{\ctime}
	= \curr
		+ \frac{\ctime - \curr[\efftime]}{\next[\efftime] - \curr[\efftime]} (\next - \curr)
\end{equation}
We will further require the ``continuous-to-discrete'' correspondence
\begin{equation}
\label{eq:tinv}
\tinv(\ctime)
	= \sup\setdef{\run\geq\start}{\ctime\geq\curr[\efftime]}
\end{equation}
which measures the number of iterations required for the effective time $\curr[\efftime]$ of the process to reach the timestamp $\ctime$;
for future use, we also define the quantity
\begin{equation}
\label{eq:tinv-horizon}
\curr[\tinv]
	\equiv \curr[\tinv](\horizon)
	= \tinv(\curr[\efftime] + \horizon).
\end{equation}
Finally, given an arbitrary sequence $\curr[A]$, we will denote its piecewise constant interpolation as
\begin{equation}
\overline A(\ctime)
	= \curr[A]
	\quad
	\text{for all $\ctime \in [\curr[\efftime],\next[\efftime]]$, $\run\geq\start$.}
\end{equation}
Using this notation, the (affinely) interpolated process $\apt{\ctime}$ can be expressed in integral form as
\begin{equation}
\apt{\ctime}
	= \apt{\cstart}
		+ \int_{\cstart}^{\ctime}
			\bracks{\vecfield(\overline\state(\ctimealt)) + \overline\error(\ctimealt)} \dd\ctimealt
\end{equation}
where $\curr[\error]$ denotes the generalized error term of \eqref{eq:RM}.

With all this in hand, \citet[Prop.~4.1]{Ben99} provides the following general condition for $\apt{\ctime}$ to be an \ac{APT} of the mean dynamics \eqref{eq:APT}:

\begin{proposition}
\label{prop:Benaim}
Suppose that $\apt{\ctime}$ is bounded and satisfies the general condition
\begin{equation}
\label{eq:Benaim}
\lim_{\ctime\to\infty} \Delta(\ctime;\horizon)
	= 0
	\quad
	\text{for all $\horizon>0$},
\end{equation}
where
\begin{equation}
\label{eq:delta}
\txs
\Delta(\ctime;\horizon)
	= \sup_{0\leq h\leq \horizon} \norm*{\int_{\ctime}^{\ctime+h} \overline\error(\ctimealt) \dd\ctimealt}.
\end{equation}
Then, $\apt{\ctime}$ is an \ac{APT} of \eqref{eq:MD}.
\end{proposition}

\subsection{Proof of \cref{thm:APT}} \label{app:proof-APT}

Our proof of \cref{thm:APT} revolves around the direct verification of the requirement \eqref{eq:Benaim} of \cref{prop:Benaim} via the use of maximal inequalities and martingale limit theory.\footnote{\citet{Ben99} provides a set of sufficient conditions for \eqref{eq:Benaim} to hold when $\apt{\ctime}$ is generated by a \ac{RM} scheme with $\bbound_{\run}=0$ and $\sup_{\run} \sdev_{\run} < \infty$;
however, our setting requires a more general treatment.}
For convenience, we restate the theorem below in full:

\APT*

\begin{proof}

Since we have shown that $\curr$ remains bounded in \cref{prop:stability}, it suffices to verify \eqref{eq:Benaim}.

Our proof relies on the \acf{BDG} inequality \citep{Bur73,HH80} which bounds the maximal value of a martingale $\curr[S]$ via its quadratic variation as
\begin{equation}
\label{eq:BDG}
\tag{BDG}
c_{2} \exof*{\sum_{\runalt=\start}^{\run} (\iter[S] - \preiter[S])^{2}}
	\leq \exof*{\max_{\runalt=\start,\dotsc,\run} \abs{\iter[S]}^{2}}
	\leq C_{2} \exof*{\sum_{\runalt=\start}^{\run} (\iter[S] - \preiter[S])^{2}},
\end{equation}
where $c_{2},C_{2}>0$ are universal constants.
As such, applying \eqref{eq:BDG} to the martingale $S_{m} = \sum_{\runalt=\run}^{m} \iter[\step] \iter[\noise]$ (after an appropriate shift of the starting time), we get
\begin{align}
\exof*{
	\sup_{\run\leq m \leq\curr[\tinv]}
		\norm*{\sum_{\runalt=\run}^{m} \iter[\step] \iter[\noise]}^{2}}
	&\leq C_{2} \exof*{\sum_{\runalt=\run}^{\curr[\tinv]} \iter[\step]^{2} \norm{\iter[\noise]}^{2}}
	\notag\\
	&= C_{2} \sum_{\runalt=\run}^{\curr[\tinv]} \iter[\step]^{2} \iter[\sdev]^{2}
	= C_{2} \int_{\curr[\efftime]}^{\curr[\efftime]+\horizon}
		\overline\step^{2}(\ctimealt) \overline\sdev^{2}(\ctimealt) \dd\ctimealt,
\end{align}
where $\curr[\tinv] = \curr[\tinv](\horizon) = \tinv(\curr[\efftime] + \horizon)$ is defined as in \eqref{eq:tinv-horizon}.
Now, mimicking \eqref{eq:delta}, let
\begin{equation}
\label{eq:delta-0}
\Delta_{0}(\ctime;\horizon)
	= \sup_{0\leq h\leq \horizon} \norm*{\int_{\ctime}^{\ctime+h} \overline\noise(\ctimealt) \dd\ctimealt}.
\end{equation}
so our previous bound shows that
\begin{align}
\exof{\Delta_{0}(\ctime;\horizon)^{2}}
	&\leq C_{2} \int_{\ctime}^{\ctime+\horizon}
		\overline\step^{2}(\ctimealt) \overline\sdev^{2}(\ctimealt) \dd\ctimealt.
\end{align}

We will proceed to show that $\lim_{\ctime\to\infty} \Delta_{0}(\ctime;\horizon) = 0$ for all $\horizon>0$ by considering the sequence of intervals $[\runalt\horizon,(\runalt+1)\horizon]$ and using the Borel-Cantelli lemma in order to show that $\Delta_{0}(\runalt\horizon;\horizon) \to 0$ as $\runalt\to\infty$.
Indeed, we have
\begin{align}
\sum_{\runalt=\start}^{\infty} \exof{\Delta_{0}(\runalt\horizon;\horizon)^{2}}
	&\leq C_{2} \int_{\cstart}^{\infty}
		\overline\step^{2}(\ctimealt) \overline\sdev^{2}(\ctimealt) \dd\ctimealt
	= C_{2} \sum_{\run=\start}^{\infty} \curr[\step]^{2} \curr[\sdev]^{2}
	< \infty
\end{align}
with the last step following from \cref{asm:noise}.
Then, if we consider the event $\event_{\runalt}(\eps) = \{\Delta_{0}(\runalt\horizon;\horizon) > \eps\}$,
Chebysev's inequality gives
\begin{equation}
\sum_{\runalt=\start}^{\infty} \probof{\event_{\runalt}(\eps)}
	\leq \frac{\sum_{\runalt=\start}^{\infty} \exof{\Delta_{0}(\runalt\horizon;\horizon)^{2}}}{\eps^{2}} < \infty,
\end{equation}
and hence, by the Borel-Cantelli lemma, we get
\begin{equation}
\probof*{\limsup_{\runalt\to\infty} \event_{\runalt}(\eps)}
	= 0.
\end{equation}
This shows that, with probability $1$, we have $\Delta_{0}(\runalt\horizon;\horizon) \leq \eps$ for all but a finite number of $\runalt$;
put differently, the event $\event(\eps) = \{\text{$\Delta_{0}(\runalt\horizon;\horizon)$ occurs infinitely often}\} = \intersect_{\run=\start}^{\infty} \union_{\runalt=\run}^{\infty} \event_{\runalt}(\eps)$ has $\probof{\event(\eps)} = 0$.
Therefore, as a union of probability zero events, we have
\begin{equation}
\probof*{\liminf_{\runalt\to\infty} \Delta_{0}(\runalt\horizon;\horizon) > 0}
	= \probof*{\union_{\run=\start}^{\infty} \event(1/\run)}
	\leq \sum_{\run=\start}^{\infty} \probof{\event(1/\run)}
	= 0,
\end{equation}
\ie $\Delta_{0}(\runalt\horizon;\horizon) \to 0$ with probability $1$.

Thus, going back to the requirements of \cref{prop:Benaim}, we get
\begin{align}
\Delta(\runalt\horizon;\horizon)
	&= \sup_{0\leq h\leq \horizon} \norm*{\int_{\runalt\horizon}^{\runalt\horizon+h} \overline\error(\ctime) \dd\ctime}
	= \sup_{0\leq h\leq \horizon} \norm*{\int_{\runalt\horizon}^{\runalt\horizon+h} [\overline\noise(\ctime) + \overline\bias(\ctime)]\dd\ctime}
	\notag\\
	&\leq \Delta_{0}(\runalt\horizon;\horizon)
		+ \sup_{0\leq h\leq \horizon} \int_{\runalt\horizon}^{\runalt\horizon+h} \overline\bbound(\ctime) \dd\ctime.
	\notag\\
	&\leq \Delta_{0}(\runalt\horizon;\horizon)
		+ \horizon \max_{0\leq h \leq \horizon} \overline\bbound(\runalt\horizon + h).
\end{align}
Given that $\lim_{\runalt\to\infty} \bbound_{\runalt} = 0$, the above shows that $\Delta(\runalt\horizon;\horizon) \to 0$ as $\runalt\to\infty$.
Moreover, for all $\ctime \in [\runalt\horizon,(\runalt+1)\horizon]$, we have $\Delta(\ctime;\horizon) \leq 2 \Delta(\runalt\horizon;\horizon) + \Delta((\runalt+1)\horizon;\horizon)$ so $\Delta(\ctime;\horizon) \to 0$ with probability $1$.
With $\horizon>0$ arbitrary, we conclude that \eqref{eq:Benaim} holds with probability $1$, so our claim follows from \cref{prop:Benaim}.
\end{proof}

\subsection{Proof of \cref{prop:APT}}\label{app:RM-APT}

We are now in a position to prove that the generalized \ac{RM} schemes presented in \cref{sec:algorithms} comprise \aclp{APT} of the mean dynamics \eqref{eq:MD}.
For convenience, we state the relevant result below:

\algoAPT*

\begin{proof}\label{pf:prop1}


We first note that $\sum_n \curr[\step]^2  < \infty$ by our choice of step-sizes. Thus, in order to prove that $\exof{\sum_n \curr[\step] \curr[\bbound]} < \infty$, $\exof{\sum_n \curr[\step]^2 \curr[\bbound]^2} < \infty$, and $\exof{\sum_n \curr[\step]^2 \curr[\sdev]^2} < \infty$, it suffices to show $\exof{\curr[\bbound]}= \exof{\norm{\curr[\bias]}} = \bigoh( \curr[\step])$ and $\exof{\curr[\sdev]^2} \leq \sdev^2$ for some constant $\sdev$.

We next proceed method-by-method:

\para{\cref{alg:SGDA}: \Acl{SGDA}}

For \eqref{eq:SGDA}, we have $\curr[\error] = \curr[\noise] = \err(\curr;\curr[\seed])$ and $\curr[\bias] = 0$, so \cref{asm:bias} is satisfied automatically (since $\curr[\bbound] = 0$).
Our claim then follows from the stated assumptions for \eqref{eq:SFO}.

\para{\cref{alg:PPM}: \Acl{PPM}}

For \eqref{eq:PPM}, we have $\curr[\noise] = 0$ and 
\begin{align*}
\norm{\curr[\bias]} &= \norm{  \vecfield(\curr) - \vecfield(\next)  } \\
&\leq \lips  \norm{  \curr - \next  } \\
&= \curr[\step] \lips \norm{  \vecfield(\curr) }  \\
&\leq \curr[\step] \lips M = \bigoh ( \curr[\step]).
\end{align*}where $\lips$ and $M$ are the Lipschitz constant of $\vecfield$ and $\minmax$, respectively. 

\para{\cref{alg:SEG}: \Acl{SEG}}

For \eqref{eq:SEG}, we have
$\curr[\noise] = \err(\lead;\lead[\seed])$
and
$\curr[\bias] = \vecfield(\lead) - \vecfield(\curr)$ so that $\exof{\curr[\sdev]^2} \leq \sdev^2$ by \eqref{eq:SFO}. 
To verify \cref{asm:bias}, by the definition of \eqref{eq:SEG}, we have  
\begin{align}
\label{eq:bias-SEG}
\norm{\curr[\bias]}
	= \norm{\vecfield(\lead) - \vecfield(\curr)}
	&\leq \lips \norm{\lead - \curr}
	\notag\\
	&= \curr[\step] \norm{\orcl(\curr;\curr[\seed])}
	= \curr[\step] \lips \norm{\vecfield(\curr) + \err(\curr;\curr[\seed])}
	\notag\\
	&\leq \curr[\step] \lips \norm{\vecfield(\curr)}
		+ \curr[\step] \lips \norm{\err(\curr;\curr[\seed])}.
\end{align}Since $\minmax$ is assumed to be Lipschitz and $\err(\curr;\curr[\seed])$ finite variance,
taking the expectation on both sides of the above shows $\exof{\curr[\bbound]} = \bigoh( \curr[\step])$.
\newcommand{\noiselevel}{\sqrt{\run \log^{1+\frac{\epsilon}{2}} \run}}
It remains to verify that $\lim_{\run\to\infty} \curr[\bbound] =0$ with probability 1. Now, by Chebyshev's inequality and \eqref{eq:SFO}, we have
\begin{equation}
\prob\left( {\norm{\err(\curr;\curr[\seed])} \geq \noiselevel}\right)
	\leq \frac{\sdev^2}{\run \log^{1+\frac{\epsilon}{2}} \run}
\end{equation}where $\eps$ is the same as in our choice of step-size in \cref{prop:APT}.
In turn, this implies that $$\sum_{\run = 2}^\infty\prob\left( {\norm{\err(\curr;\curr[\seed])} \geq \noiselevel}\right) < \infty$$so, by the Borel-Cantelli lemma, we have $\norm{\err(\curr;\curr[\seed])} = \bigoh\left( \noiselevel \right)$ with probability $1$.
Hence, by our assumptions for the method's step-size, we get
\begin{equation}
\curr[\step] \norm{\err(\curr;\curr[\seed])}
	= \bigoh\parens*{\frac{\noiselevel}{ \sqrt{\run \log^{1+\eps} \run} }}
	= \bigoh\parens*{\frac{1}{    \log^{\frac{\eps}{4}}\run  } }
\end{equation}
so that, in view of \eqref{eq:bias-SEG}, $\bbound_{\run}\to0$ with probability $1$.


\para{\cref{alg:PEG}: \Acl{OG}}

For \eqref{eq:PEG}, we have
$\curr[\noise] = \err(\curr;\lead[\seed])$
and
$\curr[\bias] = \vecfield(\lead) - \vecfield(\curr)$,
so $\exof{\curr[\sdev]^2} = \sdev^2$ again holds by \eqref{eq:SFO}.
The bias term can then be bounded exactly as in the case of \cref{alg:SEG}.

\para{\crefrange{alg:SGDA}{alg:PEG}: Alternating \ac{RM} schemes \eqref{eq:alt-RM}}

We now show that the alternating version of \crefrange{alg:SGDA}{alg:PEG} still constitute an \ac{APT} of \eqref{eq:MD}.

By \cref{lem:alternating}, we know that the alternating version of an \ac{RM} scheme is another \ac{RM} scheme with the same noise and new bias satisfying:
\begin{align} 
\norm{\curr[\bias']} 
&\leq \norm{\curr[\bias]} +  \norm{\vecfield_{\maxvar}(\next[\minstate],\curr[\maxstate]) - \vecfield_{\maxvar}(\curr[\minstate],\curr[\maxstate])}  \nonumber\\
	&\leq  \norm{\curr[\bias]} +  \lips \norm{\next[\minstate] - \curr[\minstate]} \nonumber\\
	&\leq \norm{\curr[\bias]} +  \curr[\step]\lips \left( \norm{ \vecfield(\curr)}  +\norm{ \curr[\bias]} + \norm{\curr[\noise]    } \right)  \label{eq:bias-order}
\end{align}by the definition of an \ac{RM} scheme. Since $\curr[\step]\lips\norm{ \curr[\bias]} = o \left( \norm{ \curr[\bias]} \right)$, the rest is the same as \cref{alg:SEG}.

We also note that \eqref{eq:bias-order} can be applied recursively to show that the bias term $\curr[\bias]^{(k_1,k_2)}$ of any $(k_1,k_2)$ version of \ac{RM} schemes satisfy
\begin{equation}
\norm{\curr[\bias]^{(k_1,k_2)}} \leq (k_1+k_2-1) \Big( \norm{\curr[\bias]} +  \curr[\step]\lips \left( \norm{ \vecfield(\curr)}  +\norm{ \curr[\bias]} + \norm{\curr[\noise]    } \right) \Big),
\end{equation}thus enjoying the same properties as the vanilla alternating $(1,1)$-\ac{RM} schemes in view of \eqref{eq:bias-order}.

\para{\cref{alg:SPSA}: \Acl{SPSA}}

Because of the algorithm's different oracle structure (zeroth- \vs first-order feedback), the analysis of \eqref{eq:SPSA} is different.
We begin with the algorithm's bias term, given here by
\begin{align}
\label{eq:bias-SPSA}
\curr[\bias]
	&= \exof{\curr[\signal] \given \curr[\filter]} - \vecfield(\curr)
\end{align}
with
\begin{equation}
\label{eq:signal-SPSA}
\curr[\signal]
	= \pm (\vdim/\curr[\mix])
		\, \minmax(\curr + \curr[\mix]\curr[\seed])
		\, \curr[\seed]
\end{equation}
denoting the method's one-shot \ac{SPSA} estimator.
To bound it, let
\begin{equation}
v_{\coord,\run}
	= \exof{\signal_{\coord,\run} \given \curr[\filter]}
\end{equation}
denote the $\coord$-th component of $\curr[\signal] \in \R^{\vdim}$ after having averaged out the choice of the random seed $\curr[\seed]$ (which, by default, is not $\curr[\filter]$-measurable).
We then have
\begin{align}
v_{\coord,\run}
	&= \pm \frac{\vdim}{\mix_{\run}}
		\cdot \frac{1}{2\vdim}
			\bracks[\big]{
			\minmax(\curr + \mix_{\run}\bvec_{\coord})
			- \minmax(\curr - \mix_{\run}\bvec_{\coord})}
\end{align}
where, as per our discussion in \cref{sec:algorithms}, the ``$\pm$'' sign is equal to $-1$ if $\bvec_{\coord}\in\bvecs_{\minvars}$ and $+1$ if $\bvec_{\coord}\in\bvecs_{\maxvars}$.
Then, by the mean value theorem, there exists some $\curr[\aux]$ in the line segment $\big[\curr - \mix_{\run}\bvec_{\coord},\curr + \mix_{\run}\bvec_{\coord}\big]$ such that
\begin{equation}
v_{\coord,\run}
	= \pm \pd_{\coord} \minmax(\curr[\aux])
	= \vecfield_{\coord,\run}(\curr[\aux]).
\end{equation}
Since $\vecfield$ is Lipschitz continuous, it follows that
\begin{align}
\label{eq:vecdiff}
\abs{v_{\coord,\run} - \vecfield_{\coord,\run}(\curr)}
	&= \abs[\big]{
		\vecfield_{\coord,\run}(\curr[\aux])
		- \vecfield_{\coord,\run}(\curr)
		} 
	\leq \lips \norm{\curr[\aux] - \curr}
	= \bigoh(\curr[\mix])
\end{align}
since $\curr[\aux] \in \big[\curr - \mix_{\run}\bvec_{\coord},\curr + \mix_{\run}\bvec_{\coord}\big]$.
Finally, for the oracle's variance, we have $\norm{\curr[\signal]}^{2} = \bigoh(1/\curr[\mix]^{2})$ by construction so, under the stated assumptions for $\curr[\step]$ and $\curr[\mix]$, \cref{asm:noise} is satisfied and our claim follows from \cref{thm:APT}.
\end{proof}

\section{Convergence analysis: Proof of \crefrange{thm:repeller}{thm:attract}}
\label{app:analysis}

With all this preliminary work in hand, we are finally in a position to prove \crefrange{thm:repeller}{thm:attract}.

The heavy lifting for \cref{thm:ICT} is already provided by the fact that, under the requirements of \cref{thm:APT} and/or \cref{prop:APT}, $\curr$ is an \ac{APT} of the mean dynamics \eqref{eq:MD}, so it inherits its limit structure.
\cref{thm:repeller,thm:attract} on the other hand require a completely different set of techniques and involve a much finer analysis of the process in hand.

\subsection{Avoidance of unstable periodic orbits}\label{app:avoid-repel}






%
%

While the proof of \cref{thm:repeller} is highly technical, the high-level intuition for its conclusion is crystal clear: Assume that we are given an \emph{unstable} critical point $\point^*$. Then, by the stable manifold theorem \cite{Shu87}, the set of all initializations such that the flow of \eqref{eq:MD} converges to $\point^*$ is of measure 0 in $\points$. Consequently, if the noise process $\{\curr[\noise]\}$ is such that it has ``non-negligible'' magnitude in the unstable directions near $\point^*$, then it is plausible that the \ac{RM} scheme should escape $\point^*$ along these directions. \cref{asm:unif-exciting} in \cref{thm:repeller} quantifies exactly the magnitude of noise for which we can formalize this heuristic argument.

Throughout this section we assume that we are given:
\begin{itemize}
\item A $(\vdim-\subdim)$-dimensional embedded submanifold $\set \subset \vecspace$ where $1 \leq \subdim \leq \vdim$ and $\vdim-\subdim$ is to be understood as the dimension of the unstable manifold.
\item A nonempty compact set $\cpt \subset \set$ invariant under $\flowmap \defeq \setdef{\flowmap_\ctime}{ \ctime \in \R_+}$. 
\item We also assume that $\set$ is $C^2$ is locally invariant: there exists a neighborhood $\nhd$ of $\cpt$ in $\vecspace$ and a positive time $\ctime_0$ such that
\begin{equation}
\flow[\ctime][\nhd \cap \set] \subset \set
\end{equation}for all $\abs{\ctime} \leq \ctime_0$.
\end{itemize}

We further assume that for every point $\point \in \cpt$, we have
\begin{equation}
\vecspace = \tspace\set \oplus \umfd
\end{equation}where
\begin{enumerate}[(i)]
\item $\point \to \umfd$ is a continuous map from $\cpt$ into the Grassmanian manifold $\Grass$ of $\subdim$ planes in $\vecspace$.

\item $\Jac \flow \umfd = \umfd[\flow]$ for all $\ctime \in \R, \point \in \umfd$. 

\item There exist $\lambda, C >0$ such that for all $\point \in \cpt, \wvec \in \umfd$ and $\ctime \geq 0$
\begin{equation}
\norm{  \Jac \flow \wvec } \geq C e^{\lambda \ctime} \norm{\wvec}.
\end{equation}
\end{enumerate}

We call any $\cpt$ satisfying the above an \emph{unstable invariant set}. As a simple illustration we show:
\begin{lemma}\label{lem:unstable-critical}
If $\sol$ is a critical point of $\minmax$ with any eigenvalue $\lambda$ of $\jmat(\sol)$ such that $\Re \braces{ \lambda(\jmat(\sol)) } >0$. Then $\sol$ verifies all the assumptions of an unstable invariant set.

As a corollary, $\curr$ generated by any of the \crefrange{alg:SGDA}{alg:PEG} in \cref{thm:repeller} avoids $\sol$ almost surely.
\end{lemma}
\begin{proof}
All the requirements for an unstable invariant set are readily verified by the Stable Manifold Theorem \cite{robinson1998dynamical}. The lemma follows by noting the the dimension for the unstable manifold is greater than or equal to 1; see \citep[Chap 5]{robinson1998dynamical}.
\end{proof}

A further justification of these technical assumptions is the following: Suppose $\cpt$ is a periodic orbit. We then say that $\cpt$ is (linearly) unstable if 1 is a Floquet multiplier of $\cpt$ and some multipliers have modulus strictly greater than 1 \cite{teschl2012ordinary}. If the vector field $\vecfield$ is assumed to be $C^2$, then a classical result in dynamical systems (see \eg \cite{Shu87}) states that $\cpt$ verifies all the above assumptions.


We now proceed to the proof of \cref{thm:repeller}. For ease of reading we reformulate its statements in the following more convenient form:

\newtheoremstyle{remboldstyle}
   {}{}{\itshape}{}{\bfseries}{.}{ }{}
\theoremstyle{remboldstyle}

\newtheorem{manualtheoreminner}{Theorem}
\newenvironment{manualtheorem}[1]{%
  \renewcommand\themanualtheoreminner{#1}%
  \manualtheoreminner
}{\endmanualtheoreminner}

\begin{manualtheorem}{2}
Let $\cpt$ be an unstable invariant set of $\vecfield$. Assume that
\begin{enumerate}[label=(\roman*)]
\item There exists $K >0$ such that $\norm{\curr[\noise]} \leq K$ for all $n$. 

\item $\curr[\step]$ is as in \cref{prop:APT}. 

\item There exists a neighborhood $\nhd(\cpt)$ of $\cpt$ and $b >0$ such that for all unit vector $v\in \vecspace$
\begin{equation}\nonumber
\ex \bracks{ \inner{\next[\noise] }{v}^+ \vert \curr[\filter]} \geq b \one_{  \braces{ \curr[\state] \in \nhd(\cpt) } }.
\end{equation}

\item The vector field $\vecfield$ is $C^2$.
\end{enumerate}
Then $\curr$ generated by any of the \crefrange{alg:SGDA}{alg:PEG} satisfies
\begin{equation}\nonumber
\prob \left({  \lim_{ \run \to \infty } \dist \parens{ \curr[\state], \cpt } = 0 } \right) = 0.
\end{equation}
\end{manualtheorem}

%
%
%

\begin{proof}[proof of \cref{thm:repeller}]

We first note that, for our choice of $\curr[\step]$, 
\begin{align}
\lim_{\run \to \infty}   \frac{\curr[\step]}{ \sqrt{   \insum_{\runalt=\run}^\infty  \step_{\runalt}^2  }} = 0
\end{align}so $\curr[\step] = o \parens*{ \sqrt{ \insum_{\runalt=\run}^\infty  \step_{\runalt}^2}  }$. This fact will be used in the proof when we invoke \cref{lem:prob-est} below with $\eps_\run = \bigoh(\curr[\step])$ and $\alpha_\run = \insum_{\runalt=\run}^\infty  \step_{\runalt}^2$ therein.

\para{Helper lemmas} 
We will need some technical lemmas. The first one is a deep result by \citet{BH95}, which asserts the existence of a local Lyapunov function near the unstable periodic orbits.

\newcommand{\nhdc}{\nhd(\cpt)}
\newcommand{\nhdalt}{\mathcal{W}}

For a right-differentiable function $\eta: \vecspace \to \R$ we define its right derivative $\Jac \eta$ applied to a vector $h \in \vecspace$ by 
\begin{equation}\label{eq:right-derivative}
\Jac \eta(\point)h = \lim_{\ctime \to 0^+} \frac{\eta( \point + \ctime h) - \eta(\point)}{\ctime}.
\end{equation}If $\eta$ is differentiable, then \eqref{eq:right-derivative} is simply $\inner{  \grad \eta(\point) }{ h  }$.

\begin{lemma}\label{lem:lyapunov}
There exists a compact neighborhood $\nhdc$ of $\cpt$, 
positive numbers $l, \beta >0$,
and a map $\eta: \nhd(\cpt) \to \R$ such that:
\begin{enumerate}[(i)]
\item \label[item]{itm:lyapunov1} $\eta$ is $C^2$ on $\ \nhdc \setminus \set$.

\item \label[item]{itm:lyapunov2} For all $\point \in \nhdc \cap  \set$, $\eta$ admits a right derivative $\Jac \eta (\point): \R^\vdim \to \R^ \vdim$ which is Lipschitz, convex and positively homogeneous.

\item \label[item]{itm:lyapunov3}
There exists $k>0$ and a neighborhood $\nhdalt \subset \R^\vdim$ of $0$ such that for all $\point \in \nhdc$ and $v\in \nhdalt$
\begin{equation}
\eta(\point+v) \geq \eta(\point) + {\Jac \eta(\point)}{ v} - k \norm{ v}^{2}.
\end{equation}

\item  \label[item]{itm:lyapunov4} There exists $c_1>0$ such that for all $\point \in \nhdc \setminus \set$
\begin{equation}
\norm{ \grad \eta(\point)} \geq c_1
\end{equation}and for all $\point \in \nhdc \cap \set$ and $v\in \R^\vdim$
\begin{equation}
\inner{\Jac \eta(\point)}{ v} \geq c_1 \norm{   v- \Jac \Eucl(\point) v }.
\end{equation}

\item  \label[item]{itm:lyapunov5} For all $\point \in \nhd(\gamma) \cap \set$, $u\in \tspace\set$ and $v\in \R^\vdim$
\begin{equation}
{\Jac \eta(\point)}{( u+v)} = {\Jac \eta(\point)}{v}.
\end{equation}

\item  \label[item]{itm:lyapunov6} For all $\point \in \nhdc$ we have 
\begin{equation}
{\Jac \eta(\point)}{ \vecfield(\point)} \geq \beta \eta(\point).
\end{equation}
\end{enumerate}
\end{lemma}

The second lemma we need is a probabilistic estimate from \cite{Pem92}.
\begin{lemma}\label{lem:prob-est}
Let $S_\run$ be a nonnegative stochastic process, $S_\run = S_0 + \insum_{\runalt = 1}^\run X_\runalt$ where $X_\run$ is $\curr[\filter]$-measurable.
~Let $\alpha_n \defeq \insum_{ \runalt = \run}^\infty \step_\runalt^2$.

Assume there exist a sequence $0 \leq \eps_\run = o( \sqrt{\alpha_\run})$, constants $a_1, a_2 > 0$ and an integer $N_0$ such that for all $\run \geq N_0$,
\begin{enumerate}[\upshape(\itshape i\hspace*{1pt}\upshape)]
\item \label[item]{itm:prob-est1}$\abs{X_\run}= o( \sqrt{\alpha_n} )$.
\item\label[item]{itm:prob-est2} $\one_{  \braces{S_\run > \eps_\run} } \ex \bracks{   X_{\run+1}   \vert \curr[\filter] }  \geq 0$. 
\item \label[item]{itm:prob-est3} $\ex  \bracks{  S^2_{\run+1} - S^2_\run \vert \curr[\filter]  }  \geq a_1 \step^2_{\run}$.
\item \label[item]{itm:prob-est4} $\ex \bracks{  X^2_{\run+1}  \vert \curr[\filter] }  \leq a_2 \step^2_{\run}.$
\end{enumerate}Then $\prob \parens{  \lim_{\run \to \infty}\nolimits   S_\run = 0 } = 0.$
\end{lemma}

Armed with \cref{lem:lyapunov,lem:prob-est} we are now ready to prove \cref{thm:repeller}.

Let $N \in \N$. Assume $\point_N \in \nhdc$ where $\nhdc$ is the neighborhood given by \cref{lem:lyapunov}. Define $T$ as
\begin{equation}
T \defeq \inf \braces{  \runalt \geq N: \state_\run \notin \nhdc }.
\end{equation}Evidently, $T$ is a stopping time adaptive to $\curr[\filter]$. Thus, proving \cref{thm:repeller} amounts to showing $\prob \parens{  T< \infty } = 1$.
~Without loss of generality we may assume $N=0$.

Define two sequences of random variables $\{X_\run\}_{\run \geq 2}$ and $\{S_\run\}$ as
\begin{subequations}
\begin{align}
\label{eq:Xdef}
\next[X]
	= &\parens{ \eta(\next[\state]) - \eta(  \curr[\state]  ) } \one_{\braces{\run \leq T}} + \curr[\step] \one_{ \braces{\run > T}},
\\
&S_0=\eta(\state_0),  \quad S_\run  = S_0 + \sum_{\runalt=2}^\run X_i. \label{eq:Sdef}
\end{align}
\end{subequations}
Note that $S_\run \geq 0$ \as for every $\run$. Our proof will revolve around verifying \cref{lem:prob-est}\crefrange{itm:prob-est1}{itm:prob-est4}.

\newcommand{\etacurr}{\eta(\curr[\state])}
\newcommand{\etanext}{\eta(\next[\state])}

\para{Verifying \cref{lem:prob-est}\cref{itm:prob-est1,itm:prob-est4} } By Lipschitz continuity of $\eta$ we know that
\begin{align}
\norm{\etacurr - \etanext} &\leq \lips'\norm{ \curr - \next} \nonumber\\
&= \curr[\step]\norm{  \vecfield(\curr) + \curr[\noise] + \curr[\bias]  } 
\end{align}where $\lips'$ is the Lipschitz constant of $\eta$. We have seen in the proof of \cref{prop:APT} that $\norm{ \curr[\bias] } = \bigoh \left(\step_\run (\norm{V(\curr)} + \norm{ \curr[\noise] }) \right)$. By \cref{prop:stability} and \cref{asm:bounded-noise} in \cref{thm:repeller}, we then have $\abs{\next[X]} = \bigoh({\curr[\step]}) = o ( \sqrt{\alpha_n}) $ which implies both \cref{lem:prob-est}\cref{itm:prob-est1,itm:prob-est4}.

\para{Verifying \cref{lem:prob-est}\cref{itm:prob-est2} } Let $k' = k \norm{\vecfield} + K$ where $k$ is given by \cref{lem:lyapunov}\cref{itm:lyapunov3} and $\norm{\vecfield} \defeq \sup  \setdef{ \vecfield(\point)}{ \point \in \nhdc}$ and $K$ is the uniform bound of $\curr[\noise]$. If $\run\leq T$, using \cref{lem:lyapunov}\cref{itm:lyapunov2,itm:lyapunov3,itm:lyapunov5,itm:lyapunov6} we have 
\begin{align}
\label{eq:eta-diff}
\eta(\next) - \eta(\curr) 
&\geq  \curr[\step]  \Jac \eta(\curr)  \parens*{ \curr[\noise] + \curr[\bias] + \vecfield(\curr) } - k \curr[\step]^2 \parens*{ \norm{\vecfield} + \norm{\curr[\noise]} + \norm{\curr[\bias]}   }^2 \nonumber \\
&\geq \curr[\step]  \beta \eta(\curr) +  \curr[\step]  \Jac \eta(\curr) \curr[\noise]+ \curr[\step]  \Jac \eta(\curr) \curr[\bias]- 2k' \curr[\step]^{2} - 2k \curr[\step]^{2}  \norm{\curr[\bias]}^2.
\end{align}By the same calculation leading up to \eqref{eq:bias-SEG}, \cref{asm:bounded-noise} in \cref{thm:repeller}, and the Lipschitz continuity of $\minmax$, there exists a constant $c'>0$ such that the bias sequence for \crefrange{alg:SGDA}{alg:PEG} can be bounded as $-\norm{\curr[\bias]} \geq - c' \curr[\step]$ \as. Combining this with the Lipschitz continuity of $\eta$, we can merge the last three terms in \eqref{eq:eta-diff} as 
\begin{align}
\label{eq:eta-diff-simplified}
\eta(\next) - \eta(\curr) 
&\geq \curr[\step]  \beta \eta(\curr) +  \curr[\step]  \Jac \eta(\curr) \curr[\noise]- 2k'' \curr[\step]^{2}
\end{align}for some constant $k'' > 0$. Thus
\begin{equation}\label{eq:noise-in}
\one_{\braces*{\run \leq T} } \exof{X_{\run+1} \vert \curr[\filter] } \geq \one_{\braces*{\run \leq T} } \bracks*{ \curr[\step]  \beta \eta(\curr)  - 2k'' \curr[\step]^{2} + \curr[\step] \exof{  \Jac \eta(\curr) \curr[\noise] \vert \curr[\filter] } }.
\end{equation}

By \cref{lem:lyapunov}\cref{itm:lyapunov2} again, we have
\begin{equation}\label{eq:noise-out}
\exof{  \Jac \eta(\curr) \curr[\noise] \vert \curr[\filter] }  \geq \Jac \eta(\curr)\exof{  \curr[\noise] \vert \curr[\filter] } =0
\end{equation}since we have assumed noise to be zero mean. Combining \eqref{eq:noise-in} and \eqref{eq:noise-out}, we then get
\begin{equation}\label{eq:nlT}
\one_{\braces*{\run \leq T} } \exof{X_{\run+1} \vert \curr[\filter] } \geq \one_{\braces*{\run \leq T} } \bracks*{ \curr[\step]  \beta \eta(\curr)  - 2k'' \curr[\step]^{2} } .
\end{equation}
If $\run  > T$, $X_{\run+1} = \curr[\step]$ so trivially 
\begin{equation}\label{eq:ngT}
\one_{\braces*{\run \leq T} } \exof{X_{\run+1} \vert \curr[\filter] } \geq 0.
\end{equation}Combining \eqref{eq:nlT} with \eqref{eq:ngT}, we see that  \cref{lem:prob-est}\cref{itm:prob-est2} is satisfied with $\eps_\run = \frac{k''}{\beta} \curr[\step]$.

\para{Verifying \cref{lem:prob-est}\cref{itm:prob-est3} }
We begin by observing that 
\begin{equation}\label{eq:ob}
\exof{  \next[S]^2 - \curr[S]^2 \vert \curr[\filter]  } = \exof{   \next[X]^2 \vert \curr[\filter] } + 2 \curr[S] \exof{   \next[X] \vert \curr[\filter] } .
\end{equation}
If $\curr[S] \geq \curr[\eps]$, then the right-hand side of \eqref{eq:ob} is non-negative by \cref{lem:prob-est}\cref{itm:prob-est2} that we just verified above. If  $\curr[S] < \curr[\eps]$, \eqref{eq:nlT} with \eqref{eq:ngT} imply that $\curr[S]\exof{\next[X] \vert \curr[\filter] }  \geq - \curr[\eps] k'' \curr[\step]^2  = - \bigoh ( \curr[\step]^3)$. In other words, \eqref{eq:ob} can be rewritten as
\begin{equation}\label{eq:ob2}
\exof{  \next[S]^2 - \curr[S]^2 \vert \curr[\filter]  } \geq \exof{   \next[X]^2 \vert \curr[\filter] } - \bigoh ( \curr[\step]^3) .
\end{equation}
Below, we shall prove that $\exof{   \next[X]^2 \vert \curr[\filter] }    \geq b_1 \curr[\step]^2  $ for some $b_1 >0$ and $\run$ large enough. Combining this with \eqref{eq:ob2} proves \cref{lem:prob-est}\cref{itm:prob-est3}.

\newcommand{\cfilter}{\curr[\filter]}

From \eqref{eq:eta-diff-simplified}, we deduce
\begin{equation}\label{eq:hard}
\one_{ \braces*{n\leq T} } \bracks*{   \exof{  \next[X^+] \vert \curr[\filter] }  -   \parens*{  \curr[\step] \exof{ (\Jac\eta(\curr)\curr[\noise])^+  \vert \curr[\filter] }  - k'' \curr[\step]^2 }   } \geq 0.
\end{equation}
Invoking \cref{lem:lyapunov}\cref{itm:lyapunov4} and \cref{asm:unif-exciting} in \cref{thm:repeller}, we see that 
\begin{equation}\label{eq:hard2}
\one_{ \braces{\run\leq T} \cap \braces{\curr \notin \set} } \Big(    \exof{ (\Jac \eta(\curr) \curr[\noise])^+  \vert \curr[\filter]  }  - c_1b  \Big) \geq 0
\end{equation}If $\curr \in \set$, we can choose a unit vector $v_\run \in \ker(\eye - \Jac \Eucl(\curr)  )^\perp$ where $\Eucl$ denotes the projection operator onto $\set$. By the definition of $v_\run$, we have $\inner{\curr[\noise]}{v_\run} = \inner{\curr[\noise] - \Jac \Eucl(\curr) \curr[\noise]}{v_\run} $. Let $\history=\braces{\run\leq T} \cap \braces{\curr \notin \set}.$ By \cref{lem:lyapunov}\cref{itm:lyapunov4}, Cauchy-Schwartz, and \cref{asm:unif-exciting} of \cref{thm:repeller} we get
\begin{align}
\one_{ \history }\exof{ (\Jac \eta(\curr) \curr[\noise])^+  \vert \curr[\filter]  }
&\geq c_1 \one_{ \history }\exof{ \norm{\curr[\noise] - \Jac \Eucl(\curr) \curr[\noise]}  \vert \curr[\filter]  }  \nonumber \\
&\geq c_1 \one_{ \history }\exof{ \inner{\curr[\noise] - \Jac \Eucl(\curr) \curr[\noise]}{v_\run}^+ \vert \curr[\filter]  } \nonumber \\
&= c_1 \one_{ \history }\exof{ \inner{\curr[\noise] }{v_\run}^+ \vert \curr[\filter]  } \nonumber\\
&\geq c_1b\one_\history. \label{eq:hard3}
\end{align}Combining \cref{eq:ngT,eq:hard,eq:hard2,eq:hard3} then gives 
\begin{equation}
\exof{   \next[X^+]  \vert \curr[\filter] }  \geq \curr[\step] c_1 b- k'' \curr[\step]^2.
\end{equation}On the other hand, we always have $\exof{   \next[X^2]  \vert \curr[\filter] }   \geq \exof{   \next[X^+]  \vert \curr[\filter] }$ by Jensen. It then follows that $\exof{   \next[X^2]  \vert \curr[\filter] }  \geq b_1 \curr[\step]^2$ for some $b_1 >0$ and large enough $\run$ as desired.

\para{Closing the gap} 
We have now verified \cref{lem:prob-est}\crefrange{itm:prob-est1}{itm:prob-est4}. Thus, \cref{lem:prob-est} concludes that
\begin{equation}\label{eq:final}
\probof{  \lim\nolimits_{\run \to \infty}S_\run =0  } = 0.
\end{equation}We will use \eqref{eq:final} to show that $T<\infty$ \as.

Suppose $T=\infty$. Then $\next[X] = \eta(\next[\state]) - \eta(  \curr[\state]  ) $ and $S_n = \eta(\curr)$ by \eqref{eq:Xdef}-\eqref{eq:Sdef}, and $\{ \curr\}$ remains in $\nhdc$ by definition of the stopping time $T$. \cref{thm:ICT} then asserts the the limit set $L(\{\curr\})$ of $\{\curr\}$ is a nonempty compact invariant subset of $\nhdc$, so that for all $\point' \in L(\{\curr\})$ and $\ctime\in \R $, $\flow[\ctime][\point'] \in \nhdc.$ But then \cref{lem:lyapunov}\cref{itm:lyapunov4} implies that $\eta(\flow[\ctime][\point']) \geq e^{\beta \ctime} \eta(\point')$ for all $t>0$, forcing $\eta(\point')$ to be zero. In other words, we have $L( \{ \curr\}) \subset \set$, which implies $S_\run = \eta(\curr) \to 0$ . By \eqref{eq:final}, this event occurs with probability 0, thus showing that $T$ is finite almost surely.
\end{proof}

\subsection{Convergence to \acp{ICT}} \label{app:proof-ICT}

We now prove \cref{thm:ICT}, which we restate below for convenience:

\ICT*

\begin{proof}
By \cref{thm:APT}, $\curr$ generates \ac{APT} of the mean dynamics \eqref{eq:MD}.
Now, let $\limset = \intersect_{\ctime\geq\cstart} \cl(\apt{t,\infty})$ be the limit set of $\apt{\ctime}$, \ie the set of limit points of convergent sequences $\apt{\curr[\ctime]}$ with $\lim_{\run} \curr[\ctime] = \infty$.
Our claim then follows by the limit set theorem of \citet[Theorem 8.2]{BH96}.
\end{proof}

As we discussed in the main part of our paper, the \ac{ICT} sets of $\minmax$ may exhibit a wide variety of structural properties (limit cycles, heteroclinic networks, etc.).
As a complement to this, we show below that, in \emph{gradient} systems ($\vecfield = -\nabla\obj$ for some $\obj\from\points\to\R$), \ac{ICT} sets can only be compoments of equilibria.
Specifically, building on a general result by \citet{Ben99}, we have:

\begin{proposition}
\label{prop:gradient}
Suppose that $\vecfield(\point) = -\nabla\obj(\point)$ for some $C^{\vdim}$-smooth potential function $\obj\from\points\to\R$ with a compact critical set $\crit(\obj) = \setdef{\sol}{\nabla\obj(\sol) = 0}$.
Then, every \ac{ICT} set $\set$ of \eqref{eq:MD} is contained in $\crit(\obj)$;
moreover, $\obj$ is constant on $\set$.
In particular, any \ac{ICT} set of \eqref{eq:MD} consists solely of critical points of $\obj$.
\end{proposition}

\begin{proof}
Under the stated conditions, the critical set $\sols \defeq \crit(\obj)$ of $\obj$ coincides with the set of rest points of \eqref{eq:MD}.
Moreover, by Sard's theorem \citep{Lee03}, $\obj(\sols)$ has zero Lebesgue measure and hence empty interior.
Our claim then follows from Proposition 6.4 of \citet{Ben99}.
\end{proof}

As another elementary illustration in addition to the gradient systems, one can show that for bilinear games $\minmax(\minvar,\maxvar)=\minvar\maxvar$, the \ac{ICT} sets the whole space $\R^2$.
~This can be easily seen by considering the widely known Hamiltonian function $H(\minvar,\maxvar)=\minvar^2+\maxvar^2$, which satisfies $\dot{H} = 0$ provided $(\minvar,\maxvar)$ follows \eqref{eq:MD}. An immediate consequence of this fact is that \emph{any} point on $\R^2$ lies in some ICT set of \eqref{eq:MD}, which further implies that there is no bounded attracting region, \ie attractors.

\subsection{Convergence to attractors}\label{app:proof-attract}

We now proceed with the analysis of \ac{RM} schemes in the presence of an attractor;
the relevant result is \cref{thm:attract}:

\attract*

Because of the generality of our assumptions, the proof of \cref{thm:attract} requires a range of completely different arguments and techniques.
We illustrate the main steps of our technical trajectory below:
\begin{enumerate}
\item
The first crucial component of our proof is to establish an energy function for \eqref{eq:RM} in a neighborhood of $\set$.
To do this, we rely on Conley's decomposition theorem (the so-called ``fundamental theorem of dynamical systems'') which states that the mean dynamics \eqref{eq:MD} are ``gradient-like'' in a neighborhood of an attractor, \ie they admit a (local) Lyapunov function.

\item
Because of the noise in \eqref{eq:RM}, the evolution of $\lyap$ along the trajectories of \eqref{eq:RM} could present \emph{signifcant} jumps:
in particular, a single ``bad'' realization of the noise could carry $\curr$ out of the basin of attraction of $\set$, possibly never to return.
A major difficulty here is that the driving vector field $\vecfield$ is \emph{not} assumed bounded, so it is not straightforward to establish proper control over the error terms of \eqref{eq:RM}.
However, we show that, with high probability (and, in particular, with probability at least $1-\conf$), the aggregation of these errors remains controllably small;
this is the most technically challenging part of our argument and it unfolds in a series of lemmas below.

\item
Conditioning on the above, we will show that, with probability at least $1-\conf$, the value of the trajectory's energy cannot grow more than a token threshold $\eps$;
as a result, if \eqref{eq:RM} is initialized close to $\set$, it will remain in a neighborhood thereof for all $\run$ (again, with probability at least $1-\conf$).

\item
Thanks to this ``stochastic Lyapunov stability'' result, we can regain control of the variance of the process and use martingale limit and maximal inequality arguments to show that $\curr$ converges to $\set$.
\end{enumerate}

In the rest of this section, we make this roadmap precise via a series of technical lemmas and intermediate results.

\para{A local energy function for \eqref{eq:RM}}

We begin by providing a suitable energy function for \eqref{eq:MD}.
Indeed, since $\set$ is an attractor of \eqref{eq:MD}, there exists a compact neighborhood $\cpt$ of $\set$, called the \emph{fundamental neighborhood} of $\set$, with the property that $\dist(\flow{\ctime}{\point},\set) \to 0$ as $\ctime\to\infty$ uniformly in $\point\in\cpt$.
Since all trajectories of \eqref{eq:MD} that start in $\cpt$ converge to $\set$, there are no other non-trivial invariant sets in $\cpt$ except $\set$.
Hence, with $\cpt$ compact, Conley's decomposition theorem \citep{Con78} shows that there exists a strongly smooth Lyapunov \textendash\ or ``energy'' \textendash\ function $\lyap\from\cpt\to\R$ such that
\begin{enumerate*}
[(\itshape i\hspace*{.5pt}\upshape)]
\item
$\lyap(\point) \geq 0$ with equality if and only if $\point\in\set$;
and
\item
$\dot\lyap(\point) \defeq \braket{\nabla\lyap(\point)}{\vecfield(\point)} < 0$ for all $\point\in\cpt\setminus\set$
(implying in particular that $\lyap(\flow{\ctime}{\point})$ is strictly decreasing in $\ctime$ whenever $\point\in\cpt\setminus\set$).
\end{enumerate*}

In the discrete-time context of \eqref{eq:RM}, the energy $\curr[\lyap] \defeq \lyap(\curr)$ of $\curr$ may fail to be decreasing (strictly or otherwise).
However, a simple Taylor expansion with Lagrange remainder yields the basic energy bound
\begin{align}
\label{eq:energy-bound}
\next[\lyap]
	&\leq \curr[\lyap]
		+ \curr[\step] \braket{\nabla\lyap(\curr)}{\vecfield(\curr)}
		+ \curr[\step] \curr[\snoise]
		+ \curr[\step] \curr[\sbias]
		+ \curr[\step]^{2} \curr[\scorr]^{2},
\end{align}
where the error terms $\curr[\snoise]$, $\curr[\sbias]$ and $\curr[\scorr]$ are defined as
\begin{subequations}
\label{eq:errors}
\begin{align}
\curr[\snoise]
	&= \braket{\nabla\lyap(\curr)}{\curr[\noise]}
	\\
\curr[\sbias]
	&= \curr[\bbound] \norm{\nabla\lyap(\curr)}
		+ \curr[\step] \smooth \curr[\bbound]^{2}
	\\
\curr[\scorr]^{2}
	&= \smooth \norm{\vecfield(\curr) + \curr[\noise]}^{2}
\end{align}
\end{subequations}
with $\smooth$ denoting the strong smoothness modulus of $\lyap$ over the compact set $\cpt$.
Clearly, each of these error terms can be positive, so $\curr[\lyap]$ may fail to be decreasing;
we discuss how these errors can be controlled below.

\para{Error control}

We begin by encoding the aggregation of the error terms in \eqref{eq:energy-bound} as
\begin{subequations}
\label{eq:err-agg}
\begin{align}
\label{eq:err-mart}
\curr[M]
	&= \sum_{\runalt=\start}^{\run} \iter[\step] \iter[\snoise]
\shortintertext{and}
\label{eq:err-subm}
\curr[S]
	&= \sum_{\runalt=\start}^{\run} \bracks{\iter[\step] \iter[\sbias] + \iter[\step]^{2} \iter[\scorr]^{2}}
\end{align}
\end{subequations}
Since $\exof{\curr[\snoise] \given \curr[\filter]} = 0$, we have $\exof{\curr[M] \given \curr[\filter]} = \prev[M]$, so $\curr[M]$ is a martingale;
likewise, $\exof{\curr[S] \given \curr[\filter]} \geq \prev[S]$, so $\curr[S]$ is a submartingale.
Interestingly, even though $\curr[M]$ appears more ``balanced'' as an error (because $\curr[\snoise]$ is zero-mean), it is more difficult to control because the variance of its increments is
\begin{equation}
\exof{\abs{\curr[\step] \curr[\snoise]}^{2} \given \curr[\filter]}
	= \curr[\step]^{2} \exof{\abs{\braket{\nabla\lyap(\curr)}{\curr[\noise]}}^{2} \given \curr[\filter]},
\end{equation}
so the jumps of $\curr[M]$ can become arbitrarily big if $\curr$ escapes $\cpt$ (which is the event we are trying to discount in the first place).
On that account, we will instead bound the total error increments by \emph{conditioning} everything on the event that $\curr$ remains within $\cpt$.

To make this precise, consider the ``mean square'' error process
\begin{align}
\label{eq:noise-tot}
\curr[R]
	&= \curr[M]^{2}
		+ \curr[S]
\end{align}
and the indicator events
\begin{align}
\label{eq:evt-stay}
\curr[\event]
	\equiv \curr[\event](\cpt)
	&= \braces*{\iter\in\cpt \; \text{for all $\runalt=\running,\run$}}
	\\[\smallskipamount]
\label{eq:evt-small}
\curr[\eventalt]
	\equiv \curr[\eventalt](\eps)
	&= \braces*{\iter[R] \leq \eps \; \text{for all $\runalt = \running,\run$}},
\end{align}
with the convention $\event_{0} = \eventalt_{0} = \samples$.
Moving forward, with significant hindsight, we will choose $\eps$ small enough so that
\begin{equation}
\label{eq:eps-bound}
	\setdef{\point\in\points}{\lyap(\point) \leq 2\eps + \sqrt{\eps}}
	\subseteq \cpt
\end{equation}
and we will assume that $\init$ is initialized in a neighborhood $\nhd \subseteq \cpt$ such that
\begin{equation}
\label{eq:nhd-init}
\nhd
	\subseteq \setdef{\point\in\points}{\lyap(\point) \leq \eps}.
\end{equation}
We then have the following estimates:

\begin{lemma}
\label{lem:events}
Suppose that $\init\in\nhd$ and \cref{asm:bias,asm:noise} hold.
Then
\begin{enumerate}
\item
$\next[\event] \subseteq \curr[\event]$
and
$\next[\eventalt] \subseteq \curr[\eventalt]$.
\item
$\prev[\eventalt] \subseteq \curr[\event]$.
\item
Consider the ``bad realization'' event
\begin{align}
\label{eq:evt-bad}
\curr[\tilde\eventalt]
	\defeq \prev[\eventalt] \setminus \curr[\eventalt]
	&= \prev[\eventalt] \cap \{\curr[R] > \eps\}
	\notag\\
	&= \braces*{\text{$\iter[R] \leq \eps$ for $\runalt=\running,\run-1$ and $\curr[R] > \eps$}},
\end{align}
and
let $\curr[\tilde R] = \curr[R] \one_{\prev[\eventalt]}$ denote the cumulative error subject to the noise being ``small'' until time $\run$.
Then:
\begin{equation}
\label{eq:noise-tot-cond}
\exof{\curr[\tilde R]}
	\leq \exof{\prev[\tilde R]
		+ \curr[\step] \gbound \curr[\bbound]
		+ \curr[\step]^{2} \bracks{2\smooth \gbound^{2} + (2\smooth + \gbound^{2})\curr[\sdev]^{2} + \smooth\curr[\bbound]^{2}} }
		- \eps \probof{\prev[\tilde\eventalt]},
\end{equation}
where $\gbound^{2} = \sup_{\point\in\cpt} \{ \norm{\nabla\lyap(\point)}^{2} + \norm{\vecfield(\point)}^{2} \}$
and, by convention,
$\tilde\eventalt_{0} = \varnothing$, $\tilde R_{0} = 0$.
\end{enumerate}
\end{lemma}


\begin{proof}
The first claim is obvious.
For the second, we proceed inductively:

\begin{enumerate}[leftmargin=2.5em]
\item
For the base case $\run=\start$, we have $\init[\event] = \{\init \in \cpt \} \supseteq \{\init \in \nhd \} = \samples$ (recall that $\init$ is initialized in $\nhd \subseteq \cpt$).
Since $\eventalt_{0} = \samples$, our claim follows.

\item
Inductively, suppose that $\prev[\eventalt] \subseteq \curr[\event]$ for some $\run\geq\start$.
To show that $\curr[\eventalt] \subseteq \next[\event]$,
suppose that $\iter[R] \leq \eps$ for all $\runalt=\running,\run$.
Since $\curr[\eventalt] \subseteq \prev[\eventalt]$, this implies that $\curr[\event]$ also occurs, \ie $\iter\in\cpt$ for all $\runalt=\running,\run$;
as such, it suffices to show that $\next\in\cpt$.

To do so, given that $\iter\in\nhd\subseteq\cpt$ for all $\runalt=\running\run$, the bound \eqref{eq:energy-bound} gives
\begin{equation}
\afteriter[\lyap]
	\leq \iter[\lyap]
		+ \curr[\step] \curr[\snoise]
		+ \curr[\step] \curr[\sbias]
		+ \curr[\step]^{2} \curr[\scorr]^{2},
	\quad
	\text{for all $\runalt = \running\run$},
\end{equation}
and hence, after telescoping over $\runalt = \running,\run$, we get
\begin{equation}
\next[\lyap]
	\leq \init[\lyap]
		+ \curr[M]
		+ \curr[S]
	\leq \init[\lyap]
		+ \sqrt{\curr[R]}
		+ \curr[R]
	\leq \eps
		+ \sqrt{\eps}
		+ \eps
	= 2\eps + \sqrt{\eps}.
\end{equation}
We conclude that $\lyap(\next) \leq 2\eps + \sqrt{\eps}$, \ie $\next\in\cpt$, as required for the induction.
\end{enumerate}

For our third claim, note first that
\begin{align}
\label{eq:noise-tot-upd}
\curr[R]
	&= (\prev[M] + \curr[\step]\curr[\snoise])^{2}
		+ \prev[S]
		+ \curr[\step] \curr[\sbias]
		+ \curr[\step]^{2} \curr[\scorr]^{2}
	\notag\\
	&= \prev[R]
		+ 2 \curr[\step] \curr[\snoise] \prev[M]
		+ \curr[\step]^{2} \curr[\snoise]^{2}
		+ \curr[\step] \curr[\sbias]
		+ \curr[\step]^{2} \curr[\scorr]^{2},
\end{align}
so, after taking expectations:
\begin{align}
\exof{\curr[R] \given \curr[\filter]}
	&= \prev[R]
		+ 2 \prev[M] \curr[\step]  \exof{\curr[\snoise] \given \curr[\filter]}
		+ \exof{
			\curr[\step]^{2} \curr[\snoise]^{2}
			+ \curr[\step] \curr[\sbias]
			+ \curr[\step]^{2} \curr[\scorr]^{2}
			\given \curr[\filter] }
	\geq \prev[R]
\end{align}
\ie $\curr[R]$ is a submartingale.
To proceed, let $\curr[\tilde R] = \curr[R] \one_{\prev[\eventalt]}
$ so
\begin{align}
\label{eq:noise-tot-cond1}
\curr[\tilde R]
	&= \prev[R] \one_{\prev[\eventalt]}
		+ (\curr[R] - \prev[R]) \one_{\prev[\eventalt]}
	\notag\\
	&= \prev[R] \one_{\preprev[\eventalt]}
		- \prev[R] \one_{\prev[\tilde\eventalt]}
		+ (\curr[R] - \prev[R]) \one_{\prev[\eventalt]},
	\notag\\
	&= \prev[\tilde R]
		+ (\curr[R] - \prev[R]) \one_{\prev[\eventalt]}
		- \prev[R] \one_{\prev[\tilde\eventalt]},
\end{align}
where we used the fact that $\prev[\eventalt] = \preprev[\eventalt] \setminus \prev[\tilde\eventalt]$ so $\one_{\prev[\eventalt]} = \one_{\preprev[\eventalt]} - \one_{\prev[\tilde\eventalt]}$.
Then, \eqref{eq:noise-tot-upd} yields
\begin{align}
\curr[R] - \prev[R]
		= 2 \prev[M] \curr[\step] \curr[\snoise]
		+ \curr[\step]^{2} \curr[\snoise]^{2}
		+ \curr[\step] \curr[\sbias]
		+ \curr[\step]^{2} \curr[\scorr]^{2}
\end{align}
so
\begin{subequations}
\begin{align}
\exof{(\curr[R] - \prev[R]) \one_{\prev[\eventalt]}}
	&\label{eq:noise-tot-zero}
		= 2 \exof{\curr[\step] \prev[M]\curr[\snoise] \one_{\prev[\eventalt]}}
	\\
	&\label{eq:noise-tot-noise}
		+ \exof{\curr[\step]^{2} \curr[\snoise]^{2} \one_{\prev[\eventalt]}}
	\\
	&\label{eq:noise-tot-signal}
		+ \exof{(\curr[\step]\curr[\sbias] + \curr[\step]^{2} \curr[\scorr]^{2}) \one_{\prev[\eventalt]}}.
\end{align}
\end{subequations}
However, since $\prev[\eventalt]$ and $\prev[M]$ are both $\curr[\filter]$-measurable, we have the following estimates:
\begin{enumerate}
\item
For the noise term in \eqref{eq:noise-tot-zero}, we have:
\begin{equation}
\exof{\prev[M] \curr[\snoise] \one_{\prev[\eventalt]}}
	= \exof{\prev[M] \one_{\prev[\eventalt]} \exof{\curr[\snoise] \given \curr[\filter]}}
	= 0.
\end{equation}

\item
The term \eqref{eq:noise-tot-noise} is where the reduction to $\prev[\eventalt]$ kicks in;
indeed:
\begin{align}
\exof{\curr[\snoise]^{2} \one_{\prev[\eventalt]}}
	&= \exof{\one_{\prev[\eventalt]} \exof{ \abs{\braket{\nabla\lyap(\curr)}{\curr[\noise]}}^{2}
		\given \curr[\filter]} }
	\notag\\
	&\leq \exof{\one_{\prev[\eventalt]} \norm{\nabla\lyap(\curr)}^{2} \exof{\norm{\curr[\noise]}^{2}
		\given \curr[\filter]}}
	\commtag{by Cauchy\textendash Schwarz}
	\\
	&\leq \exof{\one_{\curr[\event]} \norm{\nabla\lyap(\curr)}^{2} \exof{\norm{\curr[\noise]}^{2}
		\given \curr[\filter]}}
	\commtag{because $\prev[\eventalt] \subseteq \curr[\event]$}
	\\
	&\leq \gbound^{2} \curr[\sdev]^{2},
	\commtag{by \cref{eq:variance}}
\end{align}
where $\gbound^{2} = \sup_{\point\in\cpt} \{ \norm{\nabla\lyap(\point)}^{2} + \norm{\vecfield(\point)}^{2} \}$.
\usetagform{default}%

\item
Finally, for the term \eqref{eq:noise-tot-signal}, we have:
\begin{align}
\label{eq:noise-tot-term1}
\exof{\curr[\scorr]^{2} \one_{\prev[\eventalt]}}
	&\leq 2\smooth \exof{\norm{\vecfield(\curr)}^{2} \one_{\curr[\event]} + \norm{\curr[\noise]}^{2}}
	\leq 2\smooth (\gbound^{2} + \curr[\sdev]^{2}),
\end{align}
where we used the fact that $\one_{\prev[\eventalt]} \leq \one_{\curr[\event]} \leq 1$.
Likewise,
\begin{align}
\label{eq:noise-tot-term2}
\exof{\curr[\sbias] \one_{\prev[\eventalt]}}
	&\leq \exof{ \gbound \curr[\bbound]
		+ \curr[\step] \smooth \curr[\bbound]^{2}}.
\end{align}
\end{enumerate}
Thus, putting together all of the above, we obtain:
\begin{equation}
\exof{(\curr[R] - \prev[R]) \one_{\prev[\eventalt]}}
	\leq \exof{\curr[\step] \gbound \curr[\bbound]
		+ \curr[\step]^{2} \bracks{2\smooth \gbound^{2} + (2\smooth + \gbound^{2})\curr[\sdev]^{2} + \smooth\curr[\bbound]^{2}}}.
\end{equation}


Going back to \eqref{eq:noise-tot-cond1}, we have $\prev[R] > \eps$ if $\prev[\tilde\eventalt]$ occurs, so the last term becomes
\begin{equation}
\label{eq:noise-tot-term3}
\exof{\prev[R] \one_{\prev[\tilde\eventalt]}}
	\geq \eps \exof{\one_{\prev[\tilde\eventalt]}}
	= \eps \probof{\prev[\tilde\eventalt]}.
\end{equation}
Our claim then follows by combining \cref{eq:noise-tot-cond1,eq:noise-tot-term1,,eq:noise-tot-term2,eq:noise-tot-term3}.
\end{proof}

\para{Containment probability}

\Cref{lem:events} is the key to showing that $\curr$ remains close to $\set$ with high probability:
we formalize this in a final intermediate result below.

\begin{proposition}
\label{prop:contain}
Fix some confidence threshold $\conf>0$.
If \eqref{eq:RM} is run with sufficiently small $\curr[\step]$ satisfying the conditions of \cref{prop:APT}, then
\begin{equation}
\label{eq:contain}
\probof{\curr[\eventalt] \given \init\in\nhd}
	\geq 1-\conf
	\quad
	\text{for all $\run=\running$}
\end{equation}
\ie $\state$ remains within the basin of attraction $\cpt$ of $\set$ with probability at least $1-\conf$.
\end{proposition}

\begin{proof}
We begin by bounding the probability of the ``bad realization'' event $\curr[\tilde\eventalt] = \prev[\eventalt] \setminus \curr[\eventalt]$.
Indeed, if $\init\in\nhd$, we have:
\begin{align}
\label{eq:prob-large1}
\probof{\curr[\tilde\eventalt]}
	&= \probof{\prev[\eventalt] \setminus \curr[\eventalt]}
	= \probof{\prev[\eventalt] \cap \{\curr[R] > \eps\}}
	\notag\\
	&= \exof{\one_{\prev[\eventalt]} \times \oneof{\curr[R] > \eps}}
	\notag\\
	&\leq \exof{\one_{\prev[\eventalt]} \times (\curr[R] / \eps)}
	\notag\\
	&= \exof{\curr[\tilde R]} / \eps
\end{align}
where, in the second-to-last line, we used the fact that $\curr[R] \geq 0$ (so $\oneof{\curr[R]>\eps} \leq \curr[R]/\eps$).
Telescoping \eqref{eq:noise-tot-cond} yields
\begin{equation}
\label{eq:prob-large2}
\exof{\curr[\tilde R]}
	\leq \exof*{\tilde R_{0}
		+ \gbound \sum_{\runalt=\start}^{\run} \iter[\step] \iter[\bbound]
		+ \sum_{\runalt=\start}^{\run} \iter[\step]^{2} \iter[\varrho]^{2} }
		- \eps \sum_{\runalt=\start}^{\run} \probof{\preiter[\tilde\eventalt]}
\end{equation}
where we set $\curr[\varrho]^{2} = 2\smooth \gbound^{2} + (2\smooth + \gbound^{2})\curr[\sdev]^{2} + \smooth\curr[\bbound]^{2}$.
Hence, combining \eqref{eq:prob-large1} and \eqref{eq:prob-large2} and invoking \cref{asm:bias,asm:noise}, we get
\(
\sum_{\runalt=\start}^{\run} \probof{\iter[\tilde\eventalt]}
	\leq \frac{1}{\eps} \exof*{
		\sum_{\runalt=\start}^{\run}
			\bracks{ \iter[\step] \gbound \iter[\bbound]
				+ \iter[\step]^{2} \iter[\varrho]^{2}}
				}
	\leq \Gamma/\eps
\)
for some $\Gamma > 0$.
Now, by choosing $\curr[\step]$ sufficiently small, we can ensure that $\Gamma/\eps < \conf$;
therefore, given that the events $\iter[\tilde\eventalt]$ are disjoint for all $\runalt=\running$, we get
\begin{equation}
\probof*{\union_{\runalt=\start}^{\run} \iter[\tilde\eventalt]}
	= \sum_{\runalt=\start}^{\run} \probof{\iter[\tilde\eventalt]}
	\leq \conf
\end{equation}
and hence:
\begin{equation}
\probof{\curr[\eventalt]}
	= \probof*{\intersect_{\runalt=\start}^{\run} \comp{\iter[\tilde\eventalt]}}
	\geq 1 - \conf,
\end{equation}
as claimed.
\end{proof}

\para{Convergence with high probability}

We are finally in a position to prove the convergence of generalized \ac{RM} algorithms:

\begin{proof}[Proof of \cref{thm:attract}]
By \cref{prop:contain}, if $\curr$ is initialized within the neighborhood $\nhd$ defined in \eqref{eq:nhd-init}, we have $\probof{\curr\in\cpt \given \init\in\nhd} \geq 1-\conf$ (note also that the neighborhood $\nhd$ is independent of the required confidence level $\conf$).
Since $\cpt$ is compact, if $\curr\in\cpt$ for all $\run$, we conclude by \cref{thm:APT} that the continuous-time interpoloation $\apt{\ctime}$ of $\curr$ is an \ac{APT} of \eqref{eq:MD}.


Now, if we write $\limset = \intersect_{\ctime\geq\cstart} \cl(\apt{t,\infty})$ for the limit set of $\apt{\ctime}$, we will have $\cpt\cap\limset \neq \varnothing$ by the compactness of $\cpt$ and the fact that $\curr\in\cpt$ for all $\run\geq\start$;
moreover, $\limset$ is itself compact as a closed subset of the compact set $\setdef{\flow{\ctime}{\point}}{0\leq\ctime\leq\horizon,\point\in\cpt}$.
Since points in $\limset\cap\cpt$ are attracted to $\set$ under \eqref{eq:MD} and $\limset$ is invariant under \eqref{eq:MD}, we conclude that $\limset\cap\set \neq \varnothing$.
However, since $\limset$ is \acl{ICT} (by \cref{thm:ICT}) and \acl{ICT} sets do not contain any proper attractors, we conclude that $\limset \subseteq \set$.
This shows that $\apt{\ctime}$ \textendash\ and hence $\curr$ \textendash\ converges to $\set$, and our proof is complete.
\end{proof}

%

\section{Omitted details for \cref{sec:spurious}}
\label{app:spurious}

\subsection{A general criterion for spurious ICT sets in almost bilinear games}

We first provide a generic criterion for the existence of spurious ICT sets in almost bilinear games \eqref{eq:perturbed-bilinear}; \cf \cref{lem:cycles-in-almost-bilinear}. We then verify that the perturbation $\perturb(\maxvar) = \frac{1}{2}\maxvar^2 - \frac{1}{4}\maxvar^4$ employed in \cref{ex:bilinear} indeed satisfies the required conditions.

\begin{lemma}
\label{lem:cycles-in-almost-bilinear}
Let $\perturb(\maxvar) = \sum_{k} a_k \maxvar^k$ be an analytic function such that
\begin{equation}
\label{eq:poly-sufficient-condition}
\sum_{k} a_{2k}k h^{2k} \prod_{i=1}^k \frac{2i-1}{2i} = 0
\end{equation}has a solution with $h>0$. Then, for small enough $\eps$, there is an ICT set of mean dynamics \eqref{eq:MD} with objective $\minmax(\minvar,\maxvar) = \minvar\maxvar + \eps \perturb(\maxvar)$ such that it does not contain any critical point.
\end{lemma}
\begin{proof}
Recall the mean dynamics \eqref{eq:MD}:
\begin{equation}
    \dotorbit{\ctime}= \vecfield(\orbit{\ctime}). \nonumber
\end{equation}In the case of $\minmax(\minvar,\maxvar) = \minvar\maxvar + \eps \perturb(\maxvar)$, \eqref{eq:MD} reads:
\begin{equation} \label{eq:vec-field-almost-bilinear}
    \left\{
                \begin{array}{ll}
                  \dot{\minvar} = -\maxvar \\
                  \dot{\maxvar} =  \minvar + \eps \perturb'(\maxvar)
                \end{array}.
              \right.
\end{equation}
The most important tool of the proof is the \emph{Abelian integral} \citep{christopher2007limit}:
\begin{equation}\tag{AI}
\label{eq:abelian-integral}
I(h) \coloneqq -\oint\limits_{\gamma_h}  \perturb'  \mathrm{d}\minvar
\end{equation}where $h>0$ is a parameter and $\gamma_h$ is a family of ovals defined as in (2.3) of \citep{christopher2007limit}. 

Suppose $\perturb(\maxvar) = a_k\maxvar^k$, so that $\perturb'(\maxvar) = ka_k \maxvar^{k-1}$. We choose $\gamma_h = \setdef{\point}{ \norm{\point} = h }.$ Then, using the polar coordinate representation, we get
\begin{align}
I(h) &= -\oint\limits_{\gamma_h}  \perturb'  \mathrm{d}\minvar  \nonumber \\
&=  ka_k\int_{0}^{2\pi} h^{k} \sin^{k}(\theta) \mathrm{d}\theta \nonumber \\
&= ka_k \cdot \left\{
                \begin{array}{ll}
                  0,  & k \text{ odd} \\
                  2\pi h^{k} \prod_{i=1}^{\frac{k}{2}} \frac{2i-1}{2i}, & k \text{ even}
                \end{array}
              \right. .
\end{align}
Since contour integrals are linear in the integrands, when $\perturb(\maxvar) = \sum_{k} a_k \maxvar^k$ in \eqref{eq:abelian-integral}, we have 
\begin{equation} \nonumber
I(h) = 4\pi \sum_{k} a_{2k}k h^{2k} \prod_{i=1}^k \frac{2i-1}{2i}.
\end{equation}
Therefore, $I(h) = 0$ if and only if \eqref{eq:poly-sufficient-condition} holds. By Theorem 2.4 in \citep{christopher2007limit}, the solution $h^*$ of $I(h^*) = 0$ then implies the existence of a limit cycle in a neighborhood of the oval $\gamma_{h^*} \coloneqq \setdef{\point}{ \norm{\point} = h^* }.$
\end{proof}

Finally, it is easy to verify that for $\perturb(\maxvar) = \frac{1}{2}\maxvar^2 - \frac{1}{4}\maxvar^4$, the condition \eqref{eq:poly-sufficient-condition} is satisfied with $h^* = \sqrt{\frac{4}{3}}$, thus implying the existence of a spurious ICT set near the neighborhood of $\setdef{\point}{ \norm{\point} = \sqrt{\frac{4}{3}} }.$

\subsection{Proof of spurious ICT sets in \cref{ex:forsaken-stable}} \label{app:spurious-forsaken}
We show the existence of two spurious ICT sets in \cref{ex:forsaken-stable}.

The mean dynamics \eqref{eq:MD} for \eqref{eq:forsaken-stable} reads:
\begin{equation} \label{eq:vec-field-forsaken}
    \left\{
                \begin{array}{ll}
                  \dot{\minvar} = -(\maxvar-0.5) - \frac{1}{2}\minvar + 2\minvar^3 - \minvar^5 \\
                  \dot{\maxvar} =  \minvar -  \frac{1}{2}\maxvar + 2 \maxvar^3 - \maxvar^5
                \end{array}.
              \right.
\end{equation}Define $r^2 \coloneqq \minvar^2 + \maxvar^2$. Then straightforward calculations show that:
\begin{align}
\frac{1}{2}\frac{\mathrm{d} }{\mathrm{d} t} r^2&= \minvar \dot{\minvar} + \maxvar \dot{\maxvar} \nonumber \\
&= -\minvar(\maxvar-0.5) - \frac{1}{2}\minvar^2 + 2\minvar^4 - \minvar^6 + \minvar \maxvar - \frac{1}{2}\maxvar^2 + 2\maxvar^4 - \maxvar^6 \nonumber \\
&= 0.5\minvar - \frac{1}{2} r^2 + 2r^4 - r^6 + 3\minvar^4\maxvar^2 + 3\minvar^2\maxvar^4 - 4\minvar^2\maxvar^2 \nonumber \\
&= 0.5\minvar - \frac{1}{2}r^2 + 2r^4 - r^6 + \minvar^2\maxvar^2\left( 3r^2 - 4 \right). \label{eq:dot-r2}
\end{align}

Substituting the value $r^2 = \frac{4}{3}$ into \eqref{eq:dot-r2}, we get
\begin{align}
\frac{1}{2}\frac{\mathrm{d} }{\mathrm{d} t} r^2 &= 0.5\minvar + \frac{1}{2}\cdot \frac{4}{3} + 2\cdot \frac{16}{9} - \frac{64}{27} \nonumber \\
&= 0.5\minvar+ \frac{14}{27} \nonumber \\
&>0 \nonumber
\end{align}since $\abs{\minvar} \leq \sqrt{\frac{4}{3}}$ on $\setdef{r \geq 0}{r^2 = \frac{4}{3}}$, whence $\dot{r} > 0$ on $\setdef{r \geq 0}{r^2 = \frac{4}{3}}$. 
Likewise, one can check that $\dot{r} < 0$ on $\setdef{r \geq 0}{r^2 =  2}$, and that there is no stationary point in the region $\set\coloneqq \setdef{r \geq 0}{ \frac{4}{3} \leq r^2 \leq 2  }$. By the Poincar\'{e}-Bendixson theorem \citep{wiggins2003introduction}, there exists at least a limit cycle in $\set$.


Finally, it is easy to see that $(\minsol,\maxsol) \simeq (0,0.49)$ is a stable critical point of \eqref{eq:forsaken-stable}. Since the region $\set$ is trapping, Poincar\'{e}'s index theorem then dictates that there exists at least another unstable limit cycle inside $\set$, establishing the claim.

\subsection{Second-order methods}
\label{app:perturb}

In this section, we discuss how to cast existing second-order methods as an RM scheme with different driving vector fields, and show that their ICT sets are similar to the first-order methods under practical settings.

\begin{example}[Second-order methods]\label{ex:2nd-order}
Thanks to the efficient implementation of Hessian-gradient multiplications \citep{pearlmutter1994fast}, a popular second-order method for min-max optimization in machine learning is the \emph{Hamiltonian descent} method \cite{abernethy2019last}.
The idea is simply to run \ac{SGD} on $\obj = \norm{\nabla\minmax}^{2}/2$, giving
\begin{equation}
\label{eq:HD}
\tag{HD}
\next
	= \curr
		- \curr[\step] \jmat(\curr) \nabla \minmax( \curr).
\end{equation}
As a (discretized) gradient system, our theory in \cref{sec:analysis} shows that \eqref{eq:HD} does not possess \ac{ICT} sets other than critical points of $\obj$.
However, a serious issue of \eqref{eq:HD} is that it ignores the \emph{sign} of gradients, \ie it does not distinguish between minimization and maximization. 
As such, it has mostly been used as a \emph{gradient penalty} scheme by mixing \eqref{eq:HD} (or its variants) with \eqref{eq:SGDA}, giving rise to a number of other second-order methods such as \acdef{SGA} \cite{balduzzi2018mechanics} and \acdef{ConO} \cite{mescheder2017numerics}.
As in \cref{sec:algorithms}, one can cast these algorithms as \ac{RM} schemes with $\vecfield(\curr)$ replaced by $(\eye-\coef\jmat(\curr))\vecfield(\curr)$, where $\coef$ is the regularization parameter.
The analysis can then proceed as in \cref{sec:analysis} by replacing \eqref{eq:MD} with the appropriate continuous-time systems.

\cref{fig:CO}(\itshape a\upshape) shows the spurious convergence of \ac{SGA} with $\coef = 0.2$ applied to \eqref{eq:forsaken-stable}.
The \ac{ICT} sets of SGA are only slightly different from \crefrange{alg:SGDA}{alg:SPSA} and, in a certain precise sense, are perturbations thereof (so they suffer the same symptoms).
\endenv
\end{example}

We now discuss how to model second-order methods as \ac{RM} schemes. We will showcase on the \acdef{ConO}:
\begin{equation} \tag{ConO}
\label{eq:ConO}
\next
	= \curr
		+ \curr[\step] (\eye-\coef\jmat(\curr))\vecfield(\curr)
\end{equation}
where $\coef >0$ is the regularization parameter. Recalling the efficient implementation scheme of Hessian-gradient multiplication \citep{pearlmutter1994fast}, we make the following assumption on the
\emph{stochastic second-order oracles} (SSO): when called at $\point = (\minvar,\maxvar)$ with random seed $\sample'\in\samples$, an SSO returns a random vector $\mathsf{JV}(\point;\seed') $ of the form
\begin{equation}
\label{eq:2nd-SO}
\tag{SSO}
\mathsf{JV}(\point;\seed')
	= \jmat(\point)\vecfield(\point)
		+ \err'(\point;\seed')
\end{equation}where $ \err'(\point;\seed')$ is assumed to be unbiased and sub-Gaussian as in \eqref{eq:err-base}. With these assumptions, one can then proceed exactly as in \cref{pf:prop1} for the \crefrange{alg:SGDA}{alg:PEG} cases to show that ConO, and its alternating version, give rise to \aclp{APT} of the continuous-time dynamics:
\begin{equation}\nonumber
\dotorbit{\ctime}
	= \bigg(\eye-\coef\jmat(\orbit{\ctime})\bigg)\vecfield(\orbit{\ctime}).
\end{equation}
\cref{fig:CO}(\itshape b\upshape) demonstrates that the spurious ICT sets of ConO for \eqref{eq:forsaken-stable} is similar to that of SGA.

Similarly, one can show (under appropriate assumptions of the oracles) the continuous-time dynamics of \acdef{SGA} is
\begin{equation}\nonumber
\dotorbit{\ctime}
	= \left(\eye- \coef  \left( \frac{\jmat(\orbit{\ctime}) - \jmat(\orbit{\ctime})^\top}{2} \right)   \right)\vecfield(\orbit{\ctime}).
\end{equation}

As explained in \cref{ex:2nd-order}, it is undesirable to set a large number of $\coef$, since then we are essentially treating $\min\max$ and $\max \min$ as the same problem. However, if $\coef$ is small, then the structure stability of \emph{hyperbolic} orbits (which holds for any stable/unstable \ac{ICT} sets}) implies that any stable (unstable) ICT set of \eqref{eq:MD} remains stable (unstable) under perturbations \citep{wiggins2003introduction}. We therefore expect the ICT sets of various second-order algorithms in \cref{ex:2nd-order} be to similar to that of first-order RM schemes.

%
%

\begin{figure*}[t]
\centering
\footnotesize
\includegraphics[height=.27\textwidth]{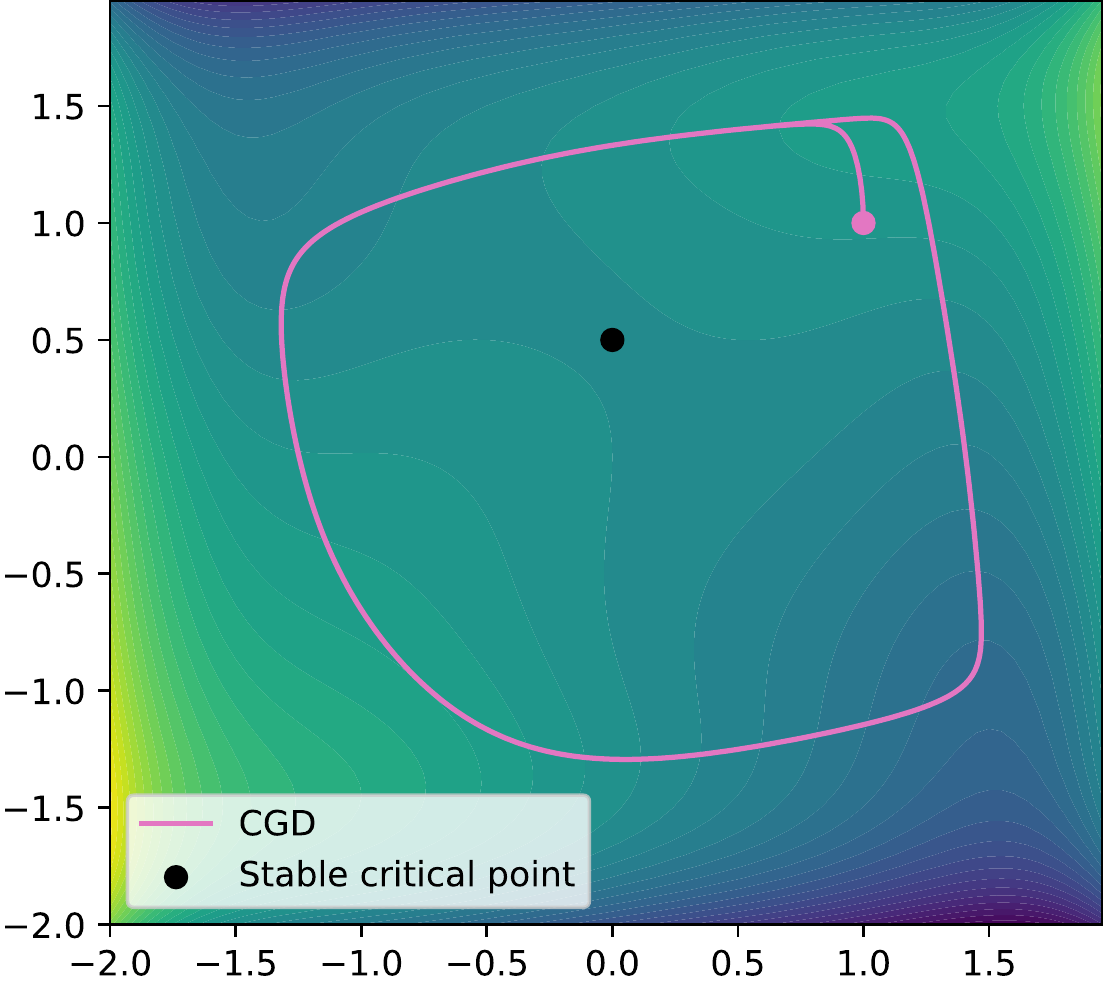}
\hfill
\includegraphics[height=.27\textwidth]{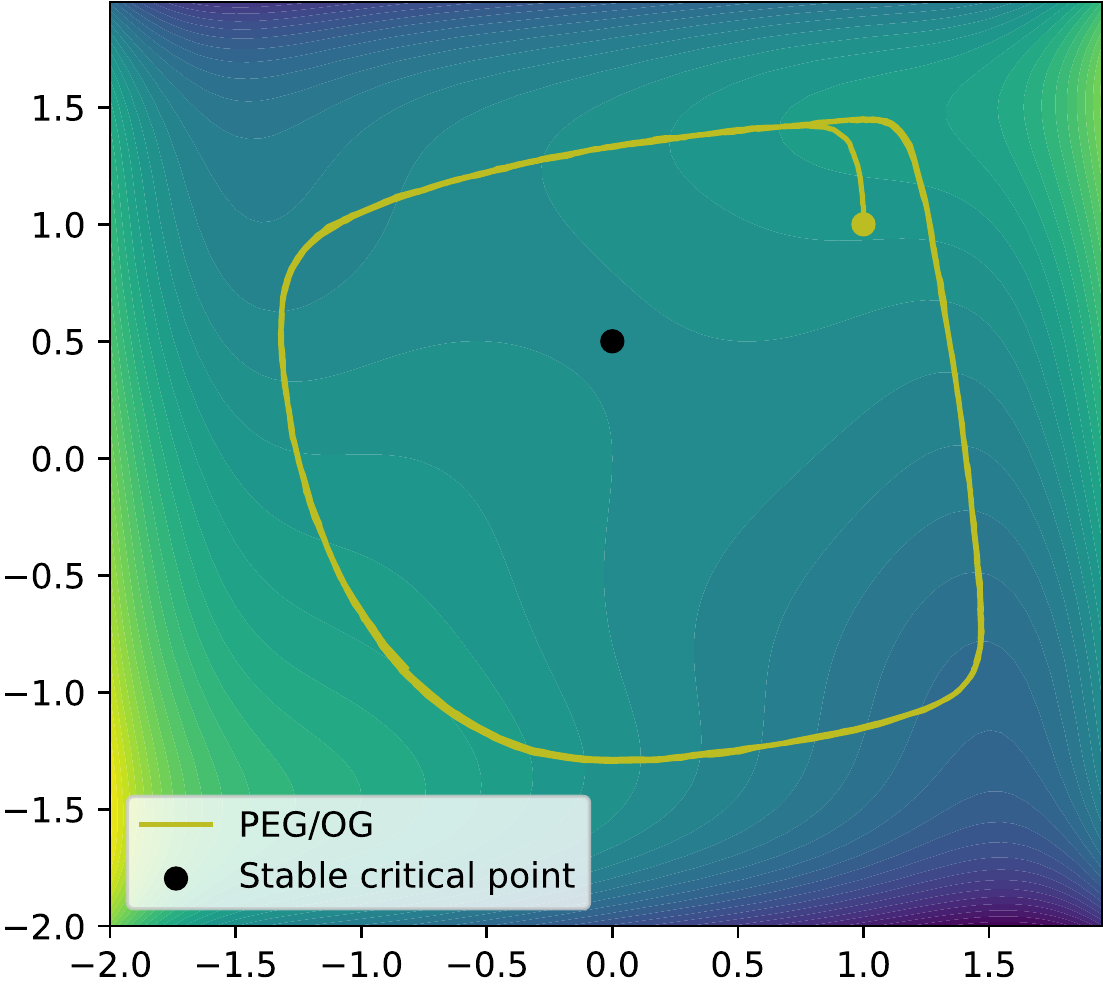}
\hfill
\includegraphics[height=.27\textwidth]{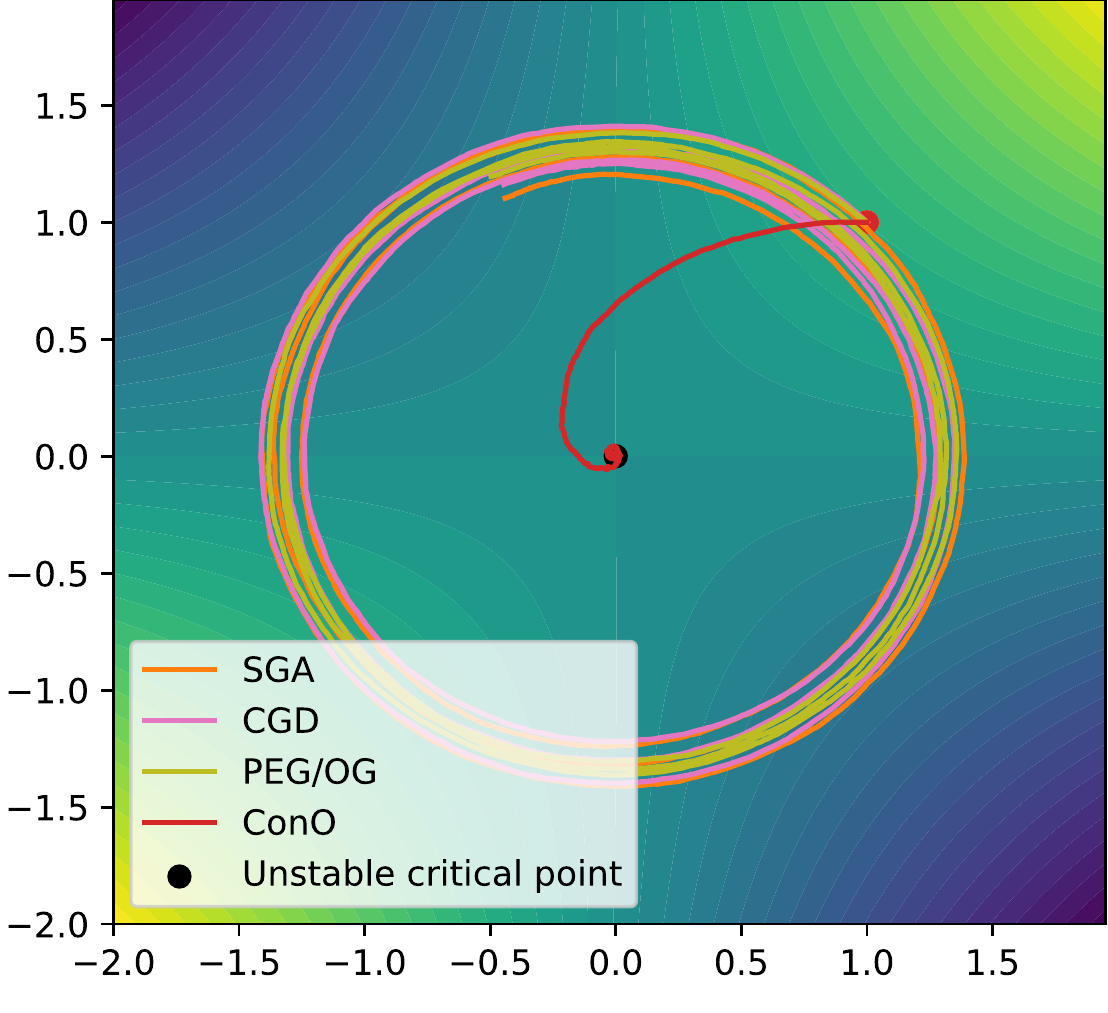}
\caption{Spurious limits of min-max optimization algorithms from the same initialization. 
From left to right:
(\itshape a\upshape)
\ac{CGD} for \eqref{eq:forsaken-stable};
(\itshape b\upshape)
\eqref{eq:PEG} for \eqref{eq:forsaken-stable};
(\itshape c\upshape)
Algorithms for \eqref{eq:perturbed-bilinear}.}
\label{fig:algs}
\end{figure*}

In addition, we have included yet another second-order method, the \acdef{CGD} \citep{schafer2019competitive}, in \cref{fig:algs}(\itshape a\upshape). For ease of comparison, we run \eqref{eq:PEG} with the same initialization in \cref{fig:algs}(\itshape b\upshape). As is evident from the figure, both algorithms perform similarly and converge straight to the spurious ICT set.

Finally, we report the behavior of various algorithms applied to the ``almost bilinear game'' \eqref{eq:perturbed-bilinear} in \cref{fig:algs}(\itshape c\upshape). In this case, all algorithms fail to escape the spurious ICT set, with the sole exception of ConO. Intriguingly, ConO converges to the \emph{unstable} critical point. A plausible explanation of this phenomenon is provided by \citep{abernethy2019last}, where it is shown that the Hamiltonian descent \eqref{eq:HD} converges to critical points for any almost bilinear game. Therefore, it is not surprising that ConO, being a mixture of SGDA and HD, also enjoys similar guarantees. Such a convergence is nonetheless highly undesirable in our example, echoing the concern that gradient penalty schemes cannot distinguish (local) $\min\max$ from $\max\min$.



\begin{figure*}[t]
\centering
\footnotesize
\includegraphics[height=.48\textwidth]{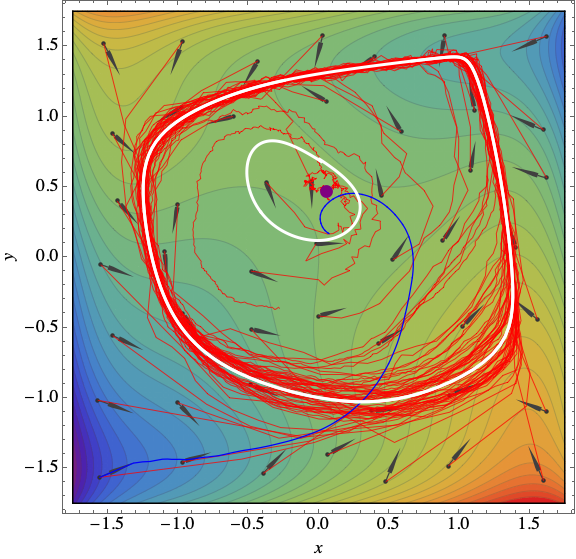}
\hfill
\includegraphics[height=.48\textwidth]{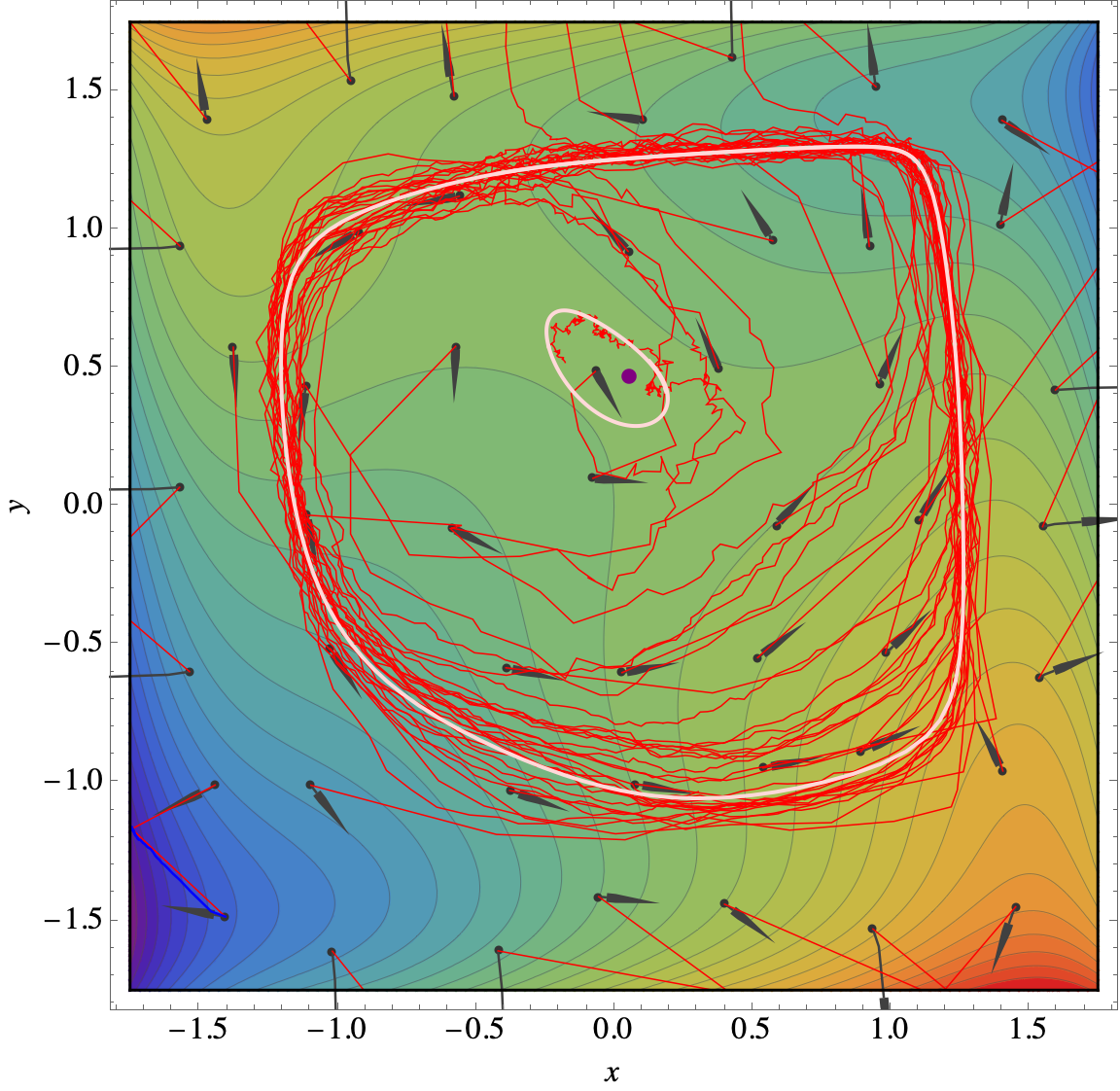}%
\caption{Spurious limits of second-order algorithms. 
From left to right:
(\itshape a\upshape)
\ac{SGA} with $\coef = 0.2$ applied to \eqref{eq:forsaken-stable};
(\itshape b\upshape)
ConO with $\coef = 0.2$ applied to \eqref{eq:forsaken-stable}.
}
\label{fig:CO}
\end{figure*}

\subsection{Constant step-sizes}\label{app:const-step}

\begin{figure*}[t]
\centering
\footnotesize
\includegraphics[width=0.24\textwidth]{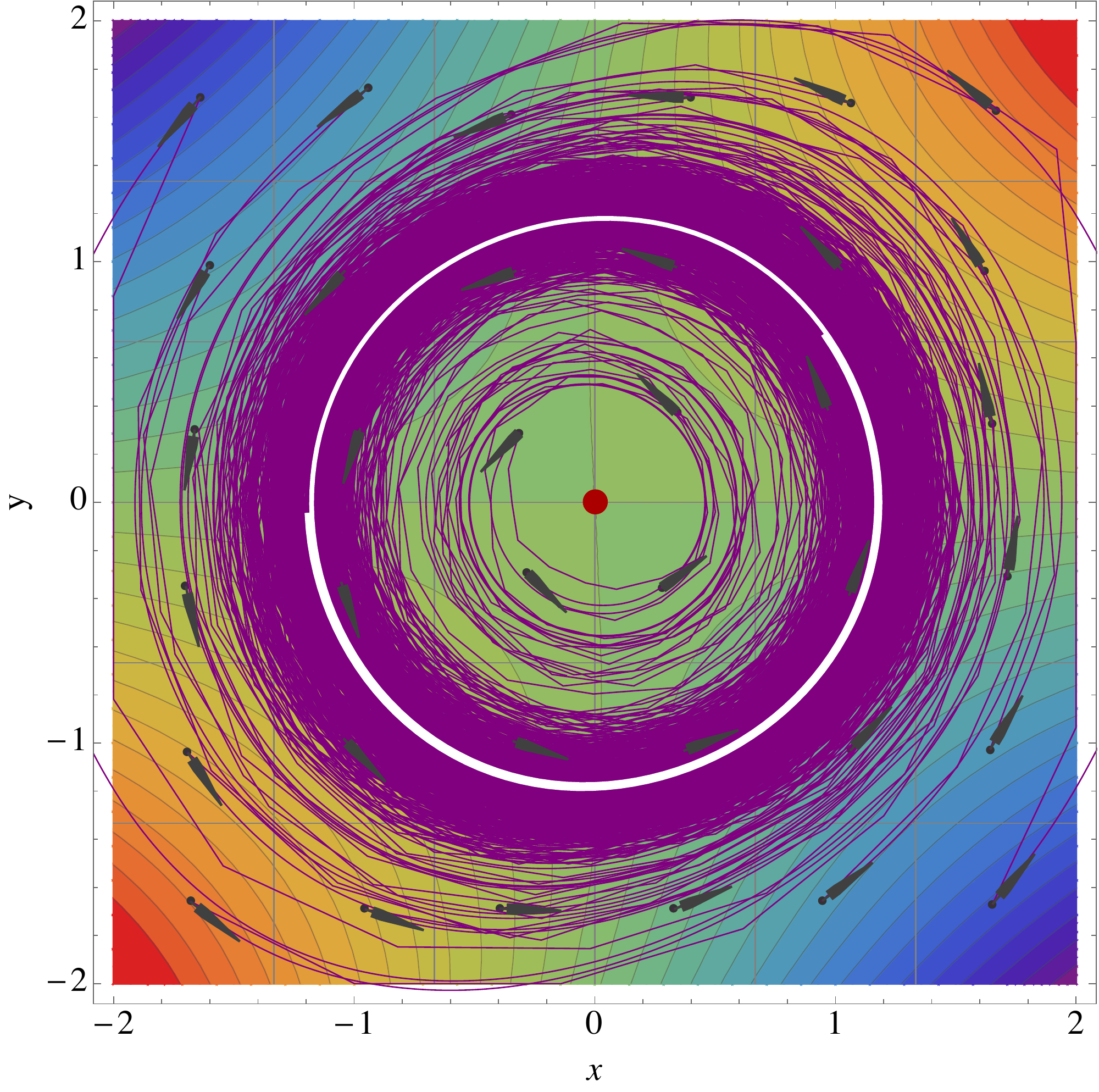} 
\includegraphics[width=0.24\textwidth]{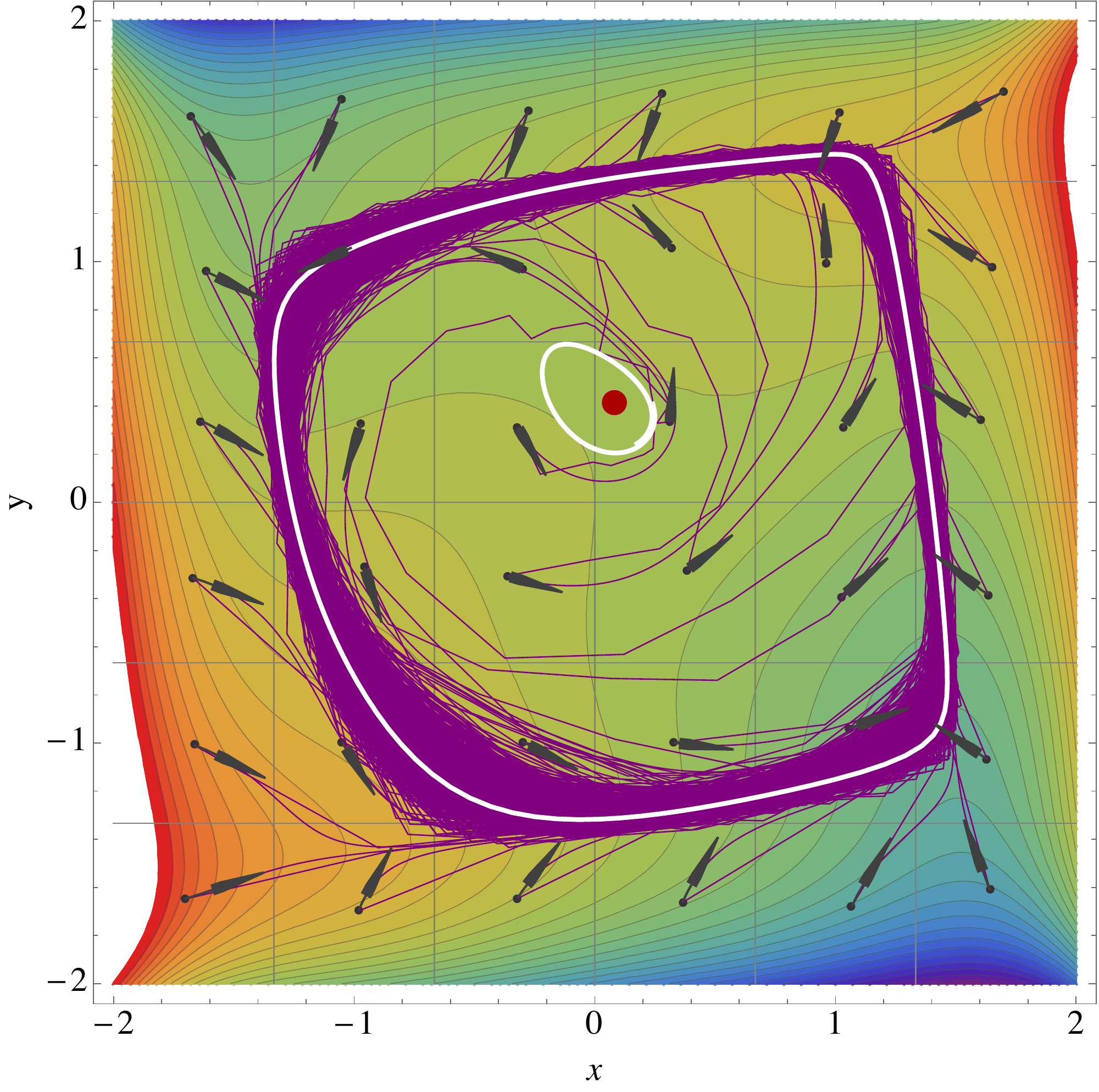}  
\includegraphics[width=0.24\textwidth]{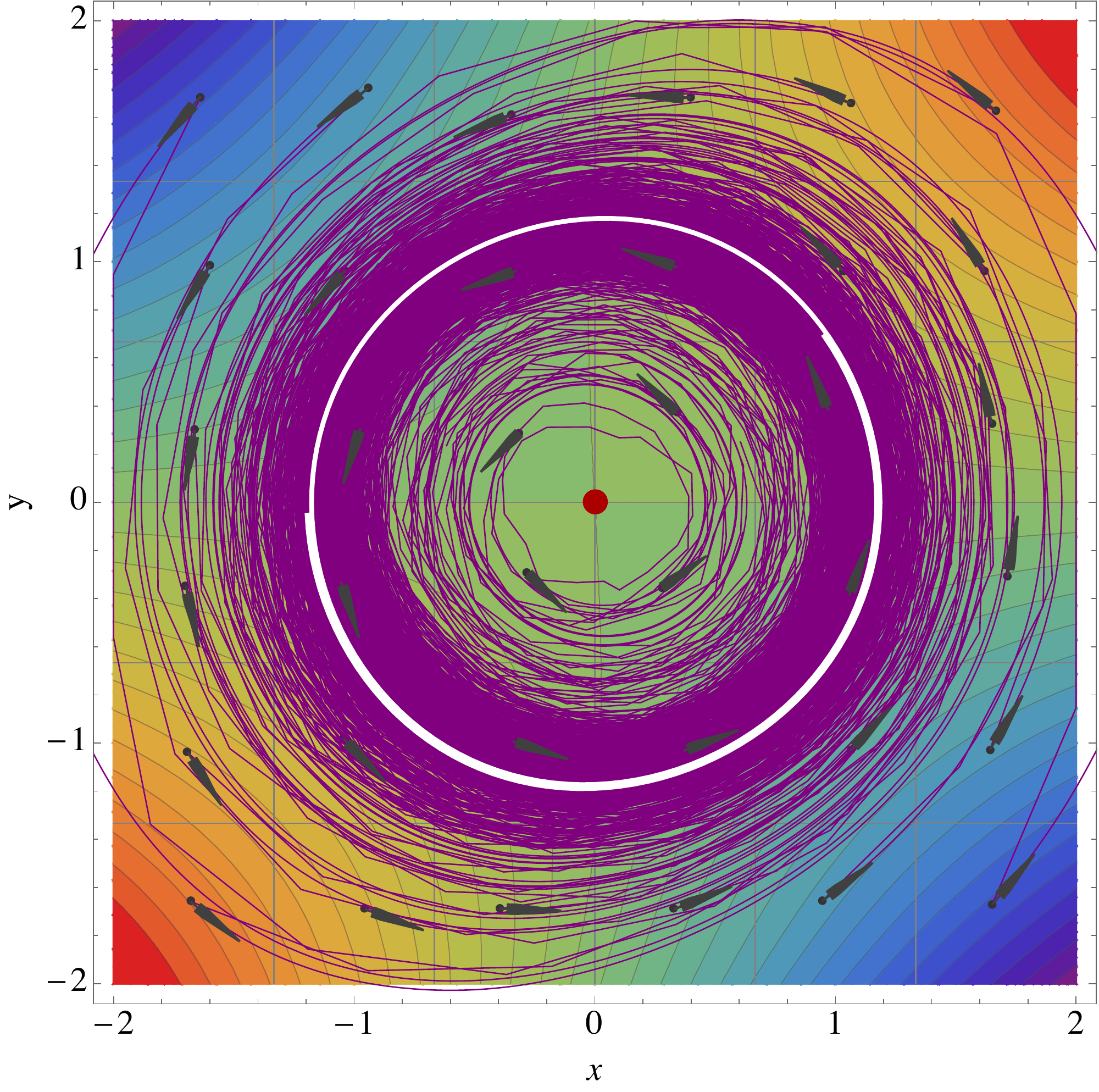} 
\includegraphics[width=0.24\textwidth]{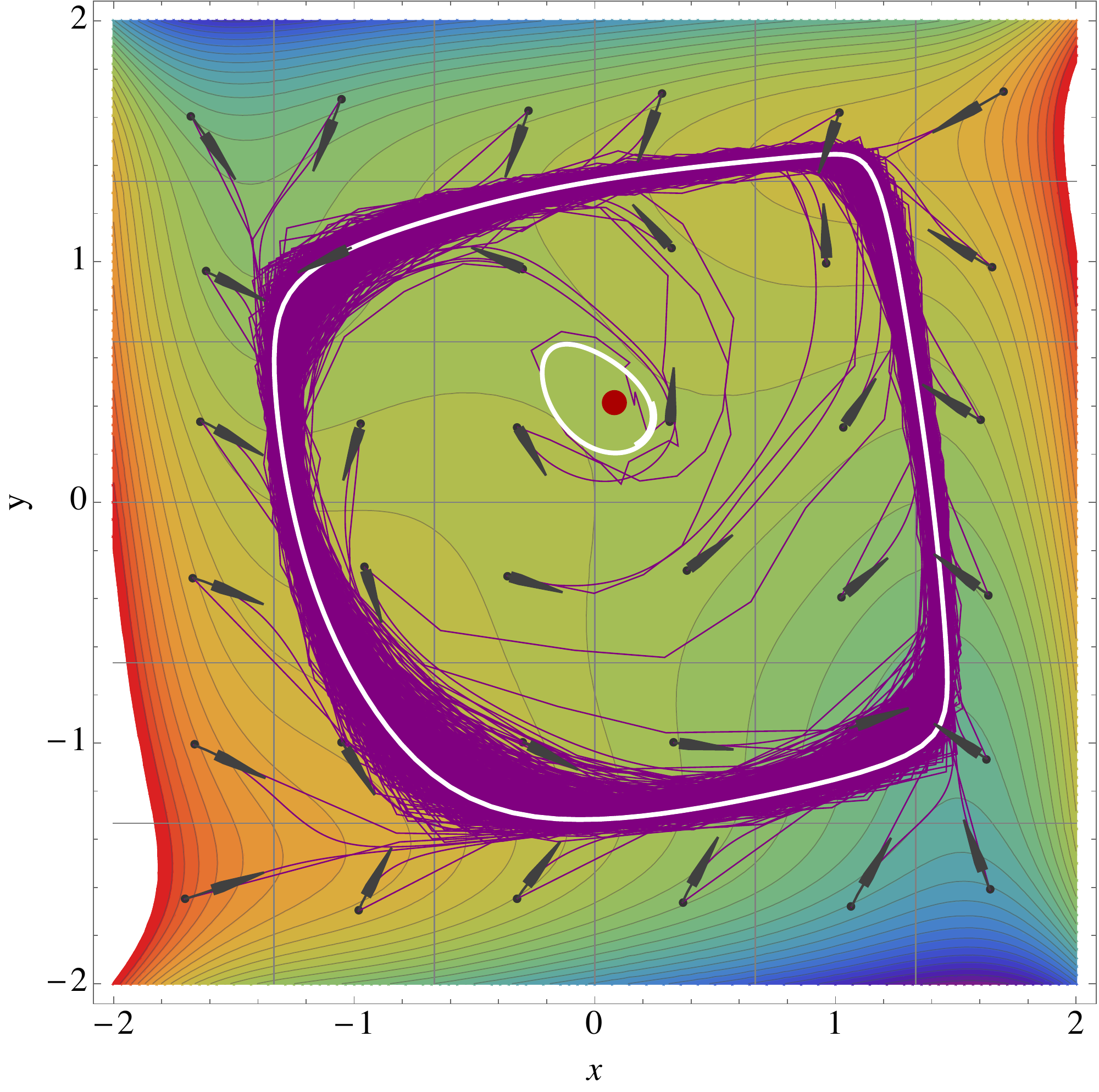}  
\caption{\ac{RM} schemes with constant step-size $\curr[\step] = 0.01$ under the same initializations.
From left to right:
(\itshape a\upshape)
\eqref{eq:SGDA} for \eqref{eq:perturbed-bilinear} with $\eps = .1$;
(\itshape b\upshape)
\eqref{eq:SGDA} for \eqref{eq:forsaken-stable};
(\itshape c\upshape)
\eqref{eq:SEG} for \eqref{eq:perturbed-bilinear} with $\eps = .1$;
(\itshape d\upshape)
\eqref{eq:SEG} for \eqref{eq:forsaken-stable}.}
\label{fig:const-step}
\end{figure*}

We report in \cref{fig:const-step} the behaviors of constant step-size \ac{RM} schemes. In accord with our intuition, these schemes exhibit concentration behaviors around the attractors.


\subsection{Adaptive methods}\label{app:adaptive}

\begin{figure*}[h]
\centering
\footnotesize
\includegraphics[width=0.45\textwidth]{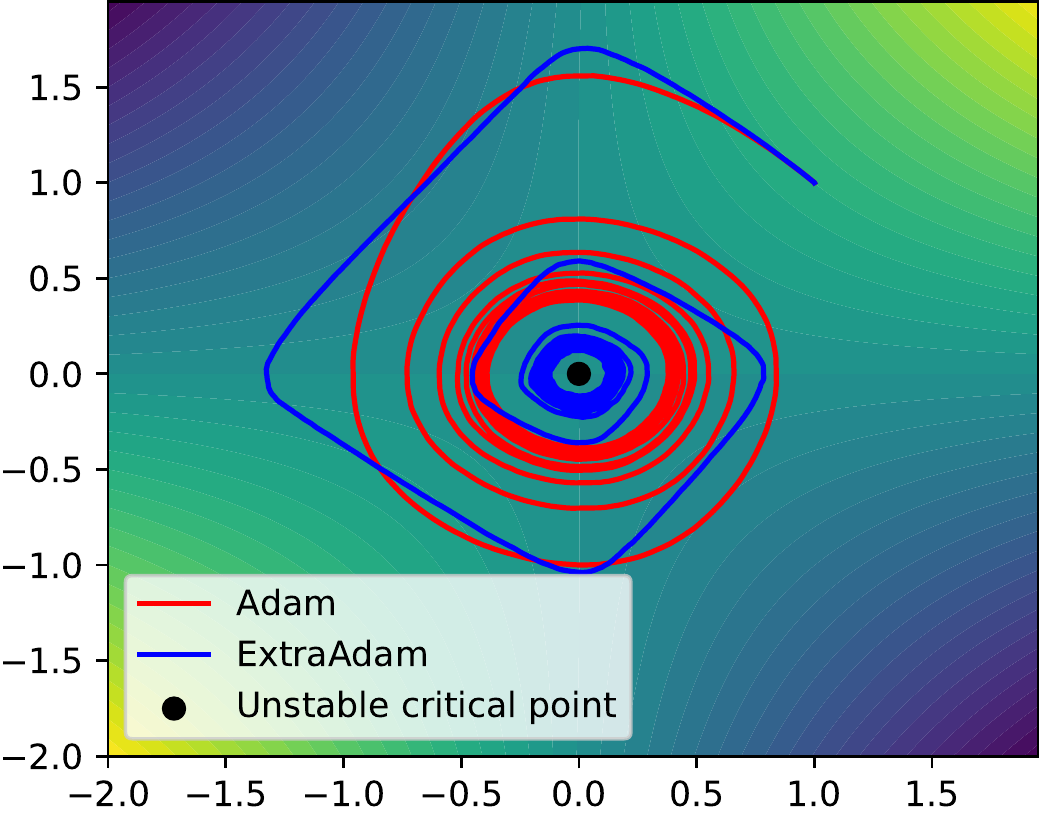} 
\hfill
\includegraphics[width=0.45\textwidth]{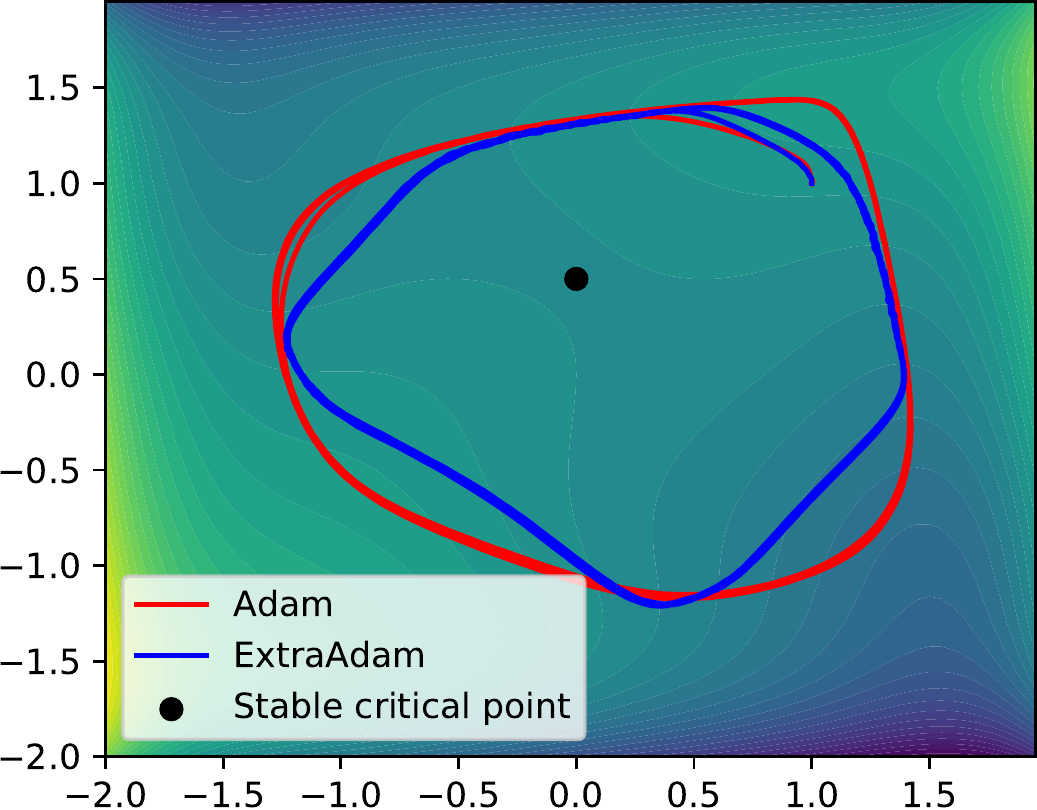}  
\caption{Adaptive algorithm.
From left to right:
(\itshape a\upshape)
Adaptive algorithms for \eqref{eq:perturbed-bilinear};
(\itshape b\upshape)
Adaptive algorithms for \eqref{eq:forsaken-stable}.}
\label{fig:adaptive}
\end{figure*}

We report in \cref{fig:adaptive} the behaviors of popular \emph{adaptive algorithms}
for min-max optimization,
including Adam \cite{kingma2014adam} and its extra-gradient variant \cite{GBVV+19}, both set to default hyperparameter values in PyTorch.
The result reveals a potentially dangerous trend:
while both Adam and ExtraAdam are able to somewhat mitigate cycling phenomena, this comes at the cost of converging to the \emph{max-min} point $(0,0)$ of \eqref{eq:perturbed-bilinear}.
In other words, the algorithm has converged, but to a very bad solution point \textendash\ an observation which, in the terminology of \citet{letcher2020impossibility}, would mean that Adam is not a ``reasonable'' algorithm.
Moreover, as all RM schemes, both adaptive methods fail to reach the ``forsaken'' solutions in \cref{ex:forsaken-stable}.


\bibliographystyle{spbasic}
\bibliography{bibtex/IEEEabrv,bibtex/Bibliography-PM,bibtex/Bibliography-YPH}

\end{document}